\documentclass[5p,times]{elsarticle}

\usepackage{rotating}
\usepackage{etoolbox}
\usepackage{amsthm}
\usepackage{amsfonts}
\usepackage{latexsym}
\usepackage{amsmath}
\usepackage{amsmath,amsfonts}
\usepackage{amssymb}
\usepackage{algorithmic}
\usepackage{algorithm}
\usepackage[linesnumbered,boxed,vlined,algo2e]{algorithm2e}
\usepackage{array}
\usepackage{multirow,booktabs}
\usepackage{textcomp}
\usepackage{url}
\usepackage{bm}
\usepackage{lipsum}
\usepackage{verbatim}
\usepackage{graphicx}
\usepackage{subcaption}
\usepackage{nomencl}

\makenomenclature

\usepackage{cite}
\usepackage{comment}
\usepackage{color}
\usepackage{makecell}
\usepackage{tabularray}
\usepackage{xcolor}

\journal{Renewable Energy (To Appear)}

\begin{document}
\begin{frontmatter}

\title{\centering Dynamic Rolling Horizon Optimization for Network-Constrained V2X Value Stacking of Electric Vehicles Under Uncertainties}

\author[A]{Canchen~Jiang}
\author[A,B]{Ariel~Liebman}
\author[C]{Bo~Jie}
\author[A,B]{Hao~Wang\corref{mycorrespondingauthor}}
\ead{hao.wang2@monash.edu}
\cortext[mycorrespondingauthor]{Corresponding author.}

\address[A]{Department of Data Science and AI, Faculty of Information Technology, Monash University, Melbourne, Victoria, Australia}
\address[B]{Monash Energy Institute, Monash University, Melbourne, Victoria, Australia}
\address[C]{Graduate School of Engineering, The University of Tokyo, Tokyo, Japan}

\begin{abstract}
Electric vehicle (EV) coordination can provide significant benefits through vehicle-to-everything (V2X) by interacting with the grid, buildings, and other EVs. This work aims to develop a V2X value-stacking framework, including vehicle-to-building (V2B), vehicle-to-grid (V2G), and energy trading, to maximize economic benefits for residential communities while maintaining distribution voltage. This work also seeks to quantify the impact of prediction errors related to building load, renewable energy, and EV arrivals. A dynamic rolling-horizon optimization (RHO) method is employed to leverage multiple revenue streams and maximize the potential of EV coordination. To address energy uncertainties, including hourly local building load, local photovoltaic (PV) generation, and EV arrivals, this work develops a Transformer-based forecasting model named Gated Recurrent Units-Encoder-Temporal Fusion Decoder (GRU-EN-TFD). The simulation results, using real data from Australia's National Electricity Market, and the Independent System Operators in New England and New York in the US, reveal that V2X value stacking can significantly reduce energy costs. The proposed GRU-EN-TFD model outperforms the benchmark forecast model. Uncertainties in EV arrivals have a more substantial impact on value-stacking performance, highlighting the significance of its accurate forecast. This work provides new insights into the dynamic interactions among residential communities, unlocking the full potential of EV batteries.
\end{abstract}

\end{frontmatter}

% \begin{table*}[h]
% \begin{framed}
\setlength{\nomitemsep}{0.02cm} % Adjusts the space between entries
\renewcommand\nomgroup[1]{%
  \item[\bfseries
  \ifstrequal{#1}{A}{Abbreviations}{%
  \ifstrequal{#1}{P}{Parameters}{%
  \ifstrequal{#1}{I}{Sets and Indices}{%
  \ifstrequal{#1}{V}{Variables}{}}}}%
]}
% Abbreviations
\nomenclature[A]{\(\text{EV}\)}{Electric Vehicle}
\nomenclature[A]{\(\text{EVPL}\)}{Electric Vehicle Parking Lot}
\nomenclature[A]{\(\text{PV}\)}{Photovoltaics}
\nomenclature[A]{\(\text{RHO}\)}{Rolling Horizon Optimization}
\nomenclature[A]{\(\text{HVAC}\)}{Heating, Ventilation and Air Conditioning}
\nomenclature[A]{\(\text{MIQP}\)}{Mixed-Integer Quadratic programming}
\nomenclature[A]{\(\text{RMSE}\)}{Root Mean Square Error}
\nomenclature[A]{\(\text{RNN}\)}{Recurrent Neural Network}
\nomenclature[A]{\(\text{TFT}\)}{Temporal Fusion Transformer}
\nomenclature[A]{\(\text{GRU}\)}{Gated Recurrent Unit}
\nomenclature[A]{\(\text{ADMM}\)}{Alternating Direction Method of Multipliers}
\nomenclature[A]{\(\text{LSTM}\)}{Long Short-Term Memory}
\nomenclature[A]{\(\text{GRU-EN-TFD}\)}{Gated Recurrent Units-Encoder-Temporal Fusion Decoder}
\nomenclature[A]{\(\text{ISO-NE}\)}{Independent System Operator New England}
\nomenclature[A]{\(\text{NYISO}\)}{New York Independent System Operator}
\nomenclature[A]{\(\text{TPT}\)}{Two-Part Tariff Pricing}
\nomenclature[A]{\(\text{TOU}\)}{Time-of-Use Pricing}
\nomenclature[A]{\(\text{NEM}\)}{National Electricity Market}
\nomenclature[A]{\(\text{V2X}\)}{Vehicle-to-Everything}
\nomenclature[A]{\(\text{V2H}\)}{Vehicle-to-Home}
\nomenclature[A]{\(\text{V2B}\)}{Vehicle-to-Building}
\nomenclature[A]{\(\text{V2G}\)}{Vehicle-to-Grid}
\nomenclature[A]{\(\text{V2V}\)}{Vehicle-to-Vehicle}
% Parameters 
\nomenclature[P]{\(\mu\)}{Charging efficiency of EV battery}
\nomenclature[P]{\(\eta\)}{Discharging efficiency of EV battery}
\nomenclature[P]{\(\alpha_b\)}{The amortized cost coefficient of the EV battery charging and discharging}
\nomenclature[P]{\(P_{\text{Bd,load}}^{i,t}\)}{The inflexible load in Community $i$ in time slot $t$}
\nomenclature[P]{\(T_{\text{OutB}}^{i,t}\)}{The temperature outside building $i$ in time slot $t$}
\nomenclature[P]{\(T_{\text{InB}}^{i,\text{Pref}}\)}{The preferred indoor temperature for building $i$ in time slot $t$}
\nomenclature[P]{\(\overline{P}_{\text{AC}}^{i,t}\)}{The upper bound of central HVAC unit power of building $i$ in time slot $t$}
\nomenclature[P]{\(\underline{P}_{\text{AC}}^{i,t}\)}{The lower bound of central HVAC unit power of building $i$ in time slot $t$}
\nomenclature[P]{\(\overline{T}^{i}\)}{The upper bound of indoor temperature of building $i$}
\nomenclature[P]{\(\underline{T}^{i}\)}{The lower bound of indoor temperature of building $i$}
\nomenclature[P]{\(\overline{G}^{i,t}\)}{The maximum allowable amount of power purchased from the grid for Community $i$ in time slot $t$}
\nomenclature[P]{\(\overline{RE}^{i,t}\)}{The maximum available amount of renewable power of Community $i$ in time slot $t$}
\nomenclature[P]{\(\overline{P}_{\text{ETs}}^{i}\)}{The upper bound for selling power from Community $i$ in the local market}
\nomenclature[P]{\(\overline{P}_{\text{ETb}}^{i}\)}{The upper bound for purchasing power by Community $i$ in the local market}
\nomenclature[P]{\(\overline{P}_{\text{EVc}}^{u,t}\)}{The maximum limit of power for charging EV $u$ in time slot $t$}
\nomenclature[P]{\(\overline{P}_{\text{EVd}}^{u,t}\)}{The maximum limit of power discharging for EV $u$ in time slot $t$}
\nomenclature[P]{\(\overline{B}_{\text{EV}}\)}{The upper bound of stored energy in EV batteries}
\nomenclature[P]{\(\underline{B}_{\text{EV}}\)}{The lower bound of stored energy in EV batteries}
\nomenclature[P]{\(R^{i}\)}{The resistance of the branch between nodes $i-1$ and $i$}
\nomenclature[P]{\(X^{i}\)}{The reactance of the branch between nodes $i-1$ and $i$}
\nomenclature[P]{\(Q_{\text{Load}}^{i,t}\)}{Reactive load at node $i$ in time slot $t$}
\nomenclature[P]{\(\overline{V}^{i}\)}{The upper bound of voltage at node $i$}
\nomenclature[P]{\(\underline{V}^{i}\)}{The lower bound of voltage at node $i$}
\nomenclature[P]{\(\overline{P}_{\text{DN}}^{i}\)}{The maximum active power flow between nodes $i-1$ and $i$}
\nomenclature[P]{\(\underline{P}_{\text{DN}}^{i}\)}{The minimum active power flow between nodes $i-1$ and $i$}
\nomenclature[P]{\(\overline{Q}_{\text{DN}}^{i}\)}{The maximum reactive power flow between nodes $i-1$ and $i$}
\nomenclature[P]{\(\underline{Q}_{\text{DN}}^{i}\)}{The minimum reactive power flow between nodes $i-1$ and $i$}
% Variables 
\nomenclature[V]{\(p_{\text{EVc}}^{u,t}\)}{Charge power of EV $u$ in time slot $t$}
\nomenclature[V]{\(p_{\text{EVd}}^{u,t}\)}{Discharge power of EV $u$ in time slot $t$}
\nomenclature[V]{\(b_{\text{EV}}^{u,t}\)}{The stored energy of EV $u$'s battery in time slot $t$}
\nomenclature[V]{\(p_{\text{renew}}^{i,t}\)}{Local renewable energy usage in Community $i$ in time slot $t$}
\nomenclature[V]{\(p_{\text{grid}}^{i,t}\)}{Grid power purchased by Community $i$ in time slot $t$}
\nomenclature[V]{\(p_{\text{V2B}}^{i,t}\)}{The V2B power from EVPL $i$ in time slot $t$}
\nomenclature[V]{\(p_{\text{V2G}}^{i,t}\)}{The V2G power from Community $i$ in time slot $t$}
\nomenclature[V]{\(p_{\text{AC}}^{i,t}\)}{Central HVAC unit power of building $i$ in time slot $t$}
\nomenclature[V]{\(T_{\text{InB}}^{i,t}\)}{The temperature inside building $i$ in time slot $t$}
\nomenclature[V]{\(P_{\text{ETb}}^{i,t}\)}{The amount of power purchased by Community $i$ from the local market in time slot $t$}
\nomenclature[V]{\(P_{\text{ETs}}^{i,t}\)}{The amount of power sold by Community $i$ to the local market  in time slot $t$}
\nomenclature[V]{\(p_{\text{Rc,Ex}}^{i,t}\)}{The active power export (or the negative active load) of Community $i$ in time slot $t$}
\nomenclature[V]{\(p_{\text{EV,Ex}}^{i,t}\)}{The aggregate power import/export for EVPL $i$}
\nomenclature[V]{\(p^{i,t}\)}{The active power from node $i-1$ to node $i$ in time slot $t$}
\nomenclature[V]{\(q^{i,t}\)}{The reactive power from node $i-1$ to node $i$ in time slot $t$}
\nomenclature[V]{\(v^{i,t}\)}{The voltage at node $i$ in time slot $t$}
\nomenclature[V]{\(x^{u,t}\)}{Binary variable avoiding EV charging and discharging simultaneously for EV $u$ in time slot $t$}
\nomenclature[V]{\(y^{i,t}\)}{Binary variable for the Big-M method for Community $i$ in time slot $t$}
\nomenclature[I]{\(\mathcal{U}\)}{The set of EVs}
\nomenclature[I]{\(\mathcal{H}\)}{The set of time slots}
\nomenclature[I]{\(\mathcal{I}\)}{The set of residential communities/buildings/EVPL/nodes}
\nomenclature[I]{\(u\)}{The $u$-th EV}
\nomenclature[I]{\(i\)}{The $i$-th residential community/building/EVPL/node}
\nomenclature[I]{\(t\)}{The $t$-th time slot}

\printnomenclature
% \end{framed}
% \end{table*}

\section{Introduction}
The global popularity of electric vehicles (EVs) has increased dramatically, due to a strong commitment to a sustainable transport sector aimed at achieving net-zero carbon emissions by 2050 \citep{amani2023technology}. Passenger EVs are mostly used for commuting on weekdays and traveling on weekends, accounting for only about 5\% of the usage time. During the remaining 95\% of their parked time, EVs can be effectively employed for additional roles by leveraging their batteries and communication capabilities \citep{mastoi2023study}. This capability forms the basis of the vehicle-to-everything (V2X) concept. V2X technology leverages the stored energy in an EV's battery for external applications \citep{singh2024multi}. Therefore, it can provide various ancillary services to power systems, which not only support power system operations but also offer economic benefits to EV owners \citep{rehman2023comprehensive}. Recent studies have focused considerable efforts on developing V2X technologies to facilitate interactions among homes, buildings, the grid, and other EVs \citep{gumrukcu2022decentralized}. For example, Vehicle-to-Building (V2B) technology is employed for peak shaving and building power demand smoothing \citep{liu2024peak}. Additionally, Vehicle-to-Grid (V2G) technology can enhance grid management by facilitating peak shaving and stabilizing voltage and frequency \citep{sagaria2025vehicle}. Meanwhile, Vehicle-to-Vehicle (V2V) technology facilitates energy transfer among connected EVs, aiming to reduce the charging load on the grid \citep{khele2023fairness}. However, uncoordinated V2X operations, which involve charging and discharging, can lead to issues such as voltage and frequency fluctuations, supply-demand imbalances, and congestion in the distribution network \citep{yu2023v2v}. Moreover, the effectiveness of V2X can be compromised by uncertainties in the energy system, including variations in renewable energy generation, energy demand, and the timing of EV charging. Hence, it is essential to develop effective approaches to optimize the benefits of EVs while considering inherent uncertainties and complying with power network constraints.

Recent research has explored methods for optimizing the scheduling of EV charging and discharging, aiming to maximize their value through interactions with homes, buildings, the grid, and other EVs. For example, Luo et al. \citep{luo2024development} introduced a prediction-based control strategy incorporating V2B technology, designed to deliver substantial economic benefits to both building consumers and EV owners. Alvaro-Hermana et al. \citep{alvaro2016peer} introduced a peer-to-peer (P2P) energy trading system among EV groups to mitigate the charging impact during peak times. Hassija et al. \citep{hassija2020blockchain} developed a blockchain model for transactions in V2G networks, noted for its scalability and efficiency. A variety of value streams of V2X have been explored, and a detailed literature review is provided in Section~\ref{LR}. Many existing studies have focused on the optimization of EV energy management to exploit a combination of multiple value streams, such as residential building demand response, energy trading, and grid services without considering network constraints. However, local network restrictions could limit the effectiveness of these multiple value streams from EVs, and neglecting these constraints could undermine the credibility of the findings.

In existing studies, the distribution network constraints have been considered for EV charging/discharging scheduling and coordination, which is commonly represented by the DistFlow Equations \citep{kiani2023admm}. For example, A topology-aware V2G energy trading strategy was introduced in \citep{zhong2018topology}, aiming at managing distribution network voltage through EV charging and discharging activities. Mazumder et al. \citep{mazumder2020ev} explored the optimal scheduling of EV charging and discharging to facilitate V2G services and assist with voltage regulation in distribution networks. Turkouglu et al. \citep{turkouglu2024maximizing} introduced a V2G-enabled model that takes into account both grid constraints and EV user preferences, aiming to achieve diverse objectives ranging from minimizing total active power loss to maximizing EV profit. However, the energy system faces significant challenges due to uncertainties such as renewable energy generation, fluctuating loads, and varied EV arrival times. These uncertainties can result in EV scheduling decisions that may not be feasible or could breach network constraints. Despite advances in forecasting renewable energy output and load demands, prediction inaccuracies remain, and their impact on the efficacy of EV scheduling and coordination is insufficiently explored.

Effective management of uncertainties can be achieved through real-time scheduling and coordination of EV charging and discharging, employing rolling horizon optimization (RHO) as a proven approach \citep{kong2024variable}. RHO is particularly useful in dynamic environments, allowing for sequential decision-making over an operational horizon based on the latest data forecasts. Recent studies have extensively applied RHO to address uncertainties in energy system management\citep{muttaqi2020adaptive}. For example, Li et al. \citep{li2019design} proposed a tiered system for managing energy demand that supports direct energy trading in a real-time market, employing RHO to cope with the stochasticity of renewable energy production. To address the uncertainty of renewable energy, Wang et al. \citep{wang2023rolling} developed a RHO-based real-time operation strategy, aimed at minimizing operational costs for all prosumers. Trinh et al. \citep{trinh2023optimal} presented an adaptive rolling horizon optimization framework to handle the uncertainties of renewable energy generation, loads, and electric pricing within the grid-connected microgrid. A detailed literature review on RHO is presented in Section~\ref{LR}. However, prediction inaccuracies are inevitable and often adversely affect the performance of V2X services. These studies did not specifically analyze how various uncertainties impact system operations and the benefits of V2X value stacking.

This work aims to bridge the research gap by exploring V2X value stacking, including V2B, V2G, and energy trading, to maximize the economic benefits while satisfying power network constraints and by quantifying the impact of prediction errors related to building load, renewable energy, and EV arrivals in residential communities. We develop a dynamic RHO method to solve the optimal charging and discharging of EVs for V2X value-stacking in real time. We develop a transformer-based forecasting model to predict solar PV output, building energy consumption, and EV arrivals. We then quantify the effects of forecasting errors in the above three variables on the value stacking. Our proposed V2X value-stacking optimization is validated using real data from Australia's National Electricity Market (NEM), Independent System Operator New England (ISO-NE) and New York Independent System Operator (NYISO) in the US. Compared to existing studies, the main contributions and findings of our work are summarized below.
\begin{itemize}
    \item This work introduces a V2X value-stacking framework that optimally integrates V2B, V2G, and EV energy trading within a distribution network, while considering network voltage limits. V2X value-stacking with a substantial number of EVs is assessed under two retail tariff structures: time-of-use (TOU) pricing and two-part tariff (TPT) pricing. The findings reveal that among the value streams evaluated, V2B consistently delivers the greatest cost savings in all studied markets: NEM, ISO-NE, and NYISO, outperforming V2G and energy trading.
    \item This work formulates a dynamic rolling-horizon optimization problem for V2X value-stacking using a shrinking window up to the end of a day, taking into account updated forecasting results as the rolling horizon advances. 
    \item This work develops a Transformer-based forecasting model, called GRU-EN-TFD (Gated Recurrent Units-Encoder-Temporal Fusion Decoder), to predict hourly building load, local PV generation, and EV arrivals. The GRU-EN-TFD model consistently outperforms the Long Short-Term Memory (LSTM) model across all three tasks, including PV generation, residential load demand, and EV arrivals.
    \item This work quantifies the impact of forecasting errors on V2X value-stacking optimization. The findings show that, compared to PV generation and residential load demand, EV arrivals have the greatest impact on additional costs incurred, highlighting the importance of accurate predictions for EV arrival in V2X value stacking.
\end{itemize}

The remainder of this paper is organized as follows. Section \ref{LR} presents related works in detail. Section \ref{Model system} outlines the system models of residential communities, EVs, and the distribution network. Section \ref{Optimal Problem Formulation} presents the dynamic rolling horizon value-stacking optimization problem and the Transformer-based forecast model for PV generation, residential community electricity consumption, and EV arrivals. Section \ref{evaluation} introduces three benchmark optimization problems along with the criteria for value-stacking performance assessment. Section \ref{Simulation and Discussion} provides a performance evaluation of V2X value stacking using real-world data and insights into the benefits of different value streams and the impact of prediction errors. Finally, this paper is concluded in Section \ref{Conclusion and Future Work}.

\section{Literature Review} \label{LR}

Recent research has explored methods for optimizing the scheduling of EV charging and discharging to enhance their value through interactions with households, buildings, local energy markets, and the grid. For example, Li et al. \citep{li2024cooperative} presented an electric-hydrogen integrated energy system that utilize V2G to alleviate peak demand and bridge energy gaps, thus reducing dependency on the utility grid. In \citep{bibak2021influences}, Bijan et al. investigated the impact of V2G operations and energy transactions on the local grid's reliability, costs, and emissions, taking into account the variability of renewable energy sources. Shurrab et al. \citep{shurrab2021efficient} proposed a V2V energy-sharing model designed to optimize the social welfare while meeting charging requirements. Al-Obaidi et al. \citep{al2020electric} proposed an EV scheduling strategy for P2P energy trading and ancillary services. Thompson et al. \citep{thompson2020vehicle} developed a comprehensive framework to assess the economic benefits of V2X, including V2B, V2V, and V2G interactions. Nonetheless, the effectiveness of these multiple value streams from EVs could be hindered by local network constraints. Without considering these constraints, the findings might lack credibility. 

Recent investigations into EV charging scheduling and coordination have increasingly taken into account the constraints of distribution networks. For example, a framework for managing residential EVs was proposed in \citep{nizami2020coordinated} to leverage V2G capabilities in mitigating grid overloads and maintaining grid voltage. Feng et al. \citep{feng2022peer} presented a P2P energy trading framework using a generalized fast dual ascent method, incorporating power network constraints to protect the distribution network's security. Affolabi et al. \citep{affolabi2021optimal} designed a dual-level market structure for facilitating energy transactions among EVs in charging stations, taking into account the voltage constraints of the power network. Hoque et al. \citep{hoque2024framework} proposed an EV energy trading framework that incorporates network export-import limits, known as the dynamic operating envelope (DOE), considering both network constraints and prosumer preferences. However, the energy system's uncertainties, like renewable energy generation, fluctuating loads, and diverse EV arrival times, pose significant challenges to scheduling EV charging. Such uncertainties can result in EV scheduling decisions that may not be viable or violate network constraints. Despite improved forecasts for renewable energy output and load demands, prediction inaccuracies remain, and the impact on the effectiveness of EV scheduling and coordination is underexplored. 

Real-time scheduling and coordination of EV charging can mitigate the adverse impacts of uncertainties, with RHO emerging as an effective method. RHO has been applied to handle uncertainties in energy system management. For example, the study in \citep{Su2020} introduced an RHO strategy for EVs in electricity exchange markets to handle volatile electricity prices. Nimalsiri et al. \citep{nimalsiri2021coordinated} developed an RHO-based approach for managing the uncertainties related to EV charging and discharging schedules, ensuring voltage levels remain safe. Hou et al. \citep{hou2020real} presented a chance-constrained programming-based rolling horizon method to address the uncertainties of EV charging demand in the energy hub system. Zaneti et al. \citep{zaneti2022sustainable} utilized RHO to minimize the operational costs of a charging station considering uncertainties in the electric bus arrival time and solar PV generation. Mohy-ud-din et al. \citep{muttaqi2021adaptive} presented a forward-looking energy management approach using RHO for virtual power plant operations, effectively countering the variability in renewable energy sources and demand. Wang et al. \citep{wang2023tri} employed RHO to adjust power usage, taking into account the variability of renewables, temperature predictions, and consumer behavior in a tri-level distribution market model. Wu et al. \citep{wu2021model} developed a two-layer control framework that incorporates RHO and a multi-uncertainty sampling method aiming to reduce forecasting errors for EV arrivals and minimize power exchange between the main network and the microgrid. Nevertheless, prediction inaccuracies are unavoidable and often adversely affect the performance of V2X services. These studies did not specifically analyze how various uncertainties impact system operations and the benefits of V2X value stacking. 

A preliminary investigation into the impact of uncertainties on the EV value-stacking model was conducted by \citep{jiang2023network}, considering factors like residential household electricity consumption and rooftop PV generation. However, many existing studies above treated EV arrivals as deterministic elements in the system, neglecting the impact of uncertainty in EV arrivals on system operations. In contrast, this work aims to bridge the research gap by exploring V2X value stacking via EV coordination that satisfies power network constraints and quantifying the impact of prediction errors in building load, renewable energy, and EV arrivals within residential communities. Through a comprehensive evaluation, our work seeks to advance the understanding of V2X value stacking and highlight the significance of various uncertainties in V2X value-stacking.
\section{System Model}\label{Model system}
\begin{figure}[!t] %!htbp
  \centering 
    \includegraphics[width=1.0\linewidth]{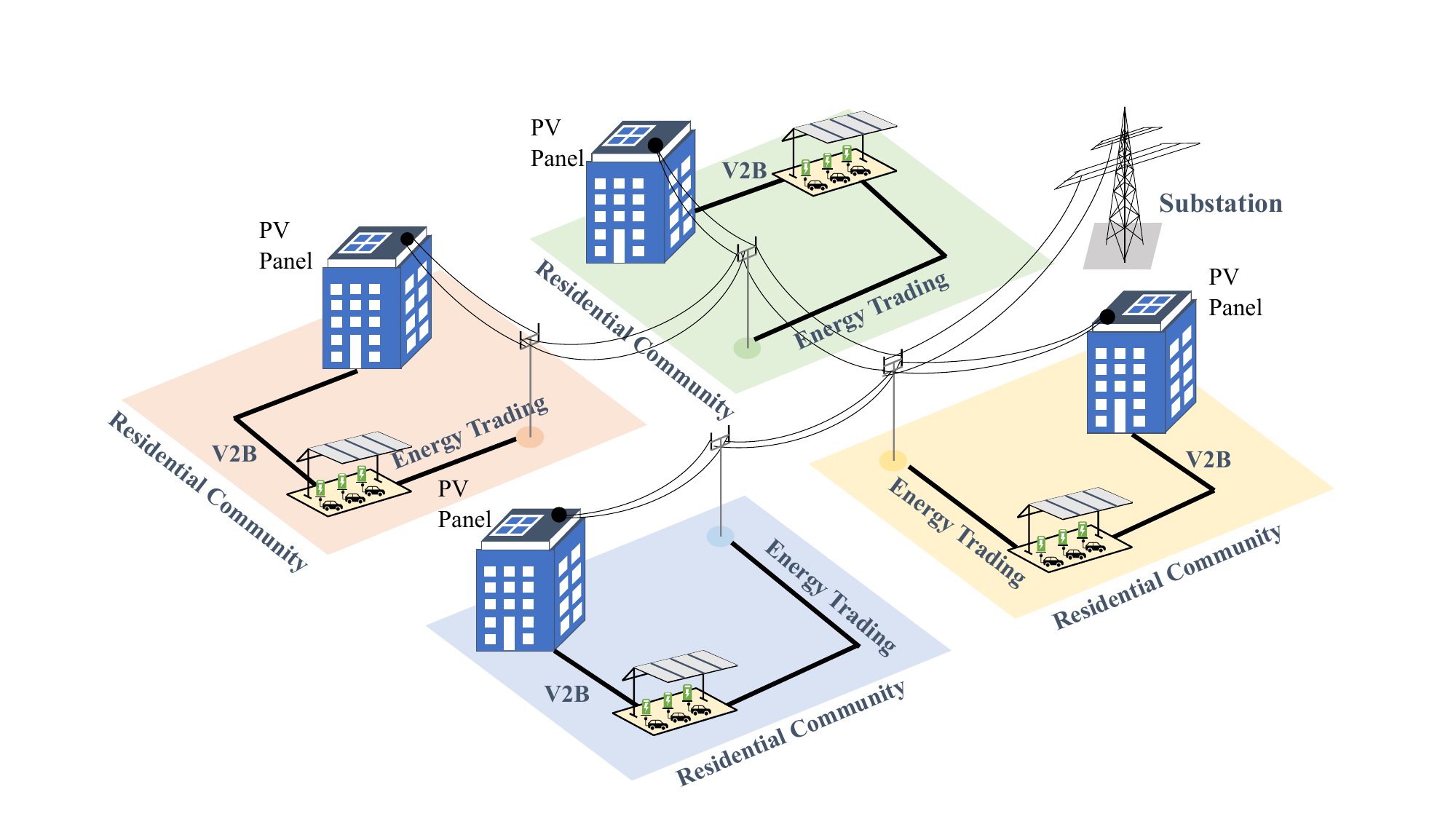} 
  \caption{Illustration of EV coordination across multiple residential communities in distribution network.} 
  \label{fig:system model}  
\end{figure}

An energy system with multiple residential communities in a distribution network is illustrated in Fig. \ref{fig:system model}. Each residential community comprises a building complex with rooftop solar PV and EVs parked in an EV parking lot (EVPL). EVs in the EVPL can supply power to the building, other communities and the grid, and also be recharged from the local energy market, and the external grid. The set of nodes in the distribution network is defined as $\mathcal{I}=\left\{1,\dots, I\right\}$, and each node is indexed by $i$. Additionally, our work uses $\mathcal{U}=\left\{1,\dots, N\right\}$ to denote the set of EVs, where $N$ is the total number of EVs. This work uses the same index $i$ for the residential community connected to node $i$. For the group of EVs parked in Community $i$'s EVPL, this work lets $\mathcal{B}^{i} \subseteq \mathcal{U}$ denote the subset of EVs in Community $i$. It is noted that each EV belongs to one group, and the EV groups $\mathcal{B}^{1}, \mathcal{B}^{2}, \dots, \mathcal{B}^{I}$ do not overlap with each other, i.e., $\mathcal{B}^{1} \cup \mathcal{B}^{2} \cup \dots \cup \mathcal{B}^{I} = \mathcal{U}$ and $\mathcal{B}^{1} \cap \mathcal{B}^{2} \cap \dots \cap \mathcal{B}^{I} = \emptyset$. The operational horizon is denoted as $\mathcal{H} \triangleq \left\{1,\dots, H\right\}$ consisting of evenly-spaced time slots. For daily operations in this study, this work sets $H=24$ hours. To operate the V2X value-stacking system, this work will develop a dynamic rolling-horizon method, which solves the problem at $t$ up to the end of the horizon $H$ with forecast information. The dynamic rolling-horizon is denoted as $\mathcal{H}_{t} = \left\{t,\dots, H\right\} \subseteq \mathcal{H}$ for $t \in \left\{1,\dots, H\right\}$.
 
\subsection{Residential Community Model } %(EVPL, Building, and Power Supply)
\subsubsection{EVPL Model}
EVPL $i$ manages the charge/discharge of EVs parked in it and can export power from EV discharging or import power to meet the EV charging. For EV $u$, its parking time window are denoted as $[t^{u}_m,t^{u}_n]\subseteq {\mathcal{H}}$, in which EV $u$ needs to be charged to the desired amount $B_{\text{EV,des}}^{u}$ before departure. Also, this work defines $p_{\text{EVc}}^{u,t}$ and $p_{\text{EVd}}^{u,t}$ as charge power and discharge power, respectively. The energy stored in EV $u$ at $t$ is represented by $b_{\text{EV}}^{u,t}$. The EV battery dynamics of EV $u$ at $t$ is
\begin{align}\label{eq1}
    &b_{\text{EV}}^{u,t} = b_{\text{EV}}^{u,t-1} + \mu p_{\text{EVc}}^{u,t} \Delta t - \frac{1}{\eta}p_{\text{EVd}}^{u,t} \Delta t,~t\in {[t^{u}_m,t^{u}_n]},\\
    &b_{\text{EV}}^{u,t} = B_{\text{EV,des}}^{u},~t= t^{u}_n,
\end{align}
where $\mu\in [0,1]$ and $\eta\in (0,1]$ represent the charging and discharging efficiencies of EV $u$, respectively. The term $\Delta t$ refers to the duration of the time interval, which is set to 60 minutes, making $\Delta t = 1$ in this study. 

The EV energy $b_{\text{EV}}^{u,t}$ is constrained as
\begin{equation}\label{eq2}
    \underline{B}_{\text{EV}}^{u} \leq b_{\text{EV}}^{u,t} \leq \overline{B}_{\text{EV}}^{u},
\end{equation}
where $\overline{B}_{\text{EV}}$ and $\underline{B}_{\text{EV}}$ are the upper and lower bounds. Note that the terminal energy level $B_{\text{EV,des}}^{u}$ is no greater than the upper bound, i.e., $B_{\text{EV,des}}^{u} \leq \overline{B}_{\text{EV}}$.

The constraints governing the charging and discharging of EVs are given by
\begin{align}
    0 \leq p_{\text{EVc}}^{u,t} \leq \overline{P}_{\text{EVc}}^{u} (1-x^{u,t}), \label{eq3}\\
    0 \leq p_{\text{EVd}}^{u,t} \leq \overline{P}_{\text{EVd}}^{u} x^{u,t}, \label{eq4}
\end{align}
where $\overline{P}_{\text{EVc}}^{u}$ and $\overline{P}_{\text{EVd}}^{u}$ represent the maximum limits for charging and discharging, respectively. To enforce that EVs do not charge and discharge simultaneously, this work introduces a binary variable $x^{u,t}\in \{ 0, 1\}$. When $x^{u,t} = 0$, EV $u$ can charge but cannot discharge. On the contrary, when $x^{u,t} = 1$, EV $u$ can discharge but cannot charge. Moreover, if EVs are not parked at EVPL, i.e., $t \notin {[t^{u}_m,t^{u}_n]}$, then both $p_{\text{EVc}}^{u,t}$ and $p_{\text{EVd}}^{u,t}$ are set to 0.

EV charging and discharging will lead to battery degradation, and the associated cost of EV battery degradation over $\mathcal{H}_{t}$ can be approximated by
\begin{equation}
   \mathbf{C}_{\text{battery}}^{t} = \alpha_b \sum \nolimits_{u\in \mathcal{U}} \left( (p_{\text{EVd}}^{u,t} \Delta t)^2+(p_{\text{EVc}}^{u,t} \Delta t)^2 \right), \notag 
\end{equation}
where $\alpha_b$ represents the amortized cost coefficient for battery usage, as outlined in study \citep{wang2016incentivizing}.

From the perspective of the EVPL, the aggregate power import/export is denoted as $p_{\text{EV,Ex}}^{i,t}$ for EVPL $i$, presented as
\begin{equation}\label{eq31}
p_{\text{EV,Ex}}^{i,t} = \sum_{u\in \mathbf{B}^{i}} (p_{\text{EVd}}^{u,t} - p_{\text{EVc}}^{u,t}).
\end{equation}
Note that if $p_{\text{EV,Ex}}^{i,t}$ takes a positive value, it indicates that the net charge and discharge of EVs lead to extra power available for export in time slot $t$. On the other hand, a negative value of $p_{\text{EV,Ex}}^{i,t}$ represents that EVPL $i$ requires additional power from outside to meet the charging demand at $t$.

\subsubsection{Residential Building Load Model}
For the residential building in Community $i$, its demand consists of inflexible load $ P_{\text{Bd,load}}^{i,t}$ and central HVAC unit load $p_{\text{AC}}^{i,t}$ in time slot $t$. This work only considers the central heating, ventilation and air conditioning (HVAC) unit as the controllable load for energy management, which controls the temperature inside the whole building denoted as $T_{\text{InB}}^{i,t}$. The dynamic of the temperature inside building \citep{cui2019peer} is as
\begin{equation}
     T_{\text{InB}}^{i,t} = T_{\text{InB}}^{i,t-1} - \frac{1}{C^{i}_{\text{AC}}R^{i}_{\text{AC}}}(T_{\text{InB}}^{i,t-1} - T_{\text{OutB}}^{i,t} + \mu^{i} R^{i}_{\text{AC}} p_{\text{AC}}^{i,t} \Delta t),
\end{equation}
where $p_{\text{AC}}^{i,t}$ represents the power of electricity consumption of HVAC unit in building $i$ in time slot $t$, $T_{\text{OutB}}^{i,t}$ denotes the outdoor temperature of building $i$ in time slot $t$. The coefficient $\mu^{i}$ denotes the operating mode of the HVAC system, in which positive values indicate cooling and negative values signify heating. Additionally, $C^{i}_{\text{AC}}$ and $R^{i}_{\text{AC}}$ are the working parameters of the HVAC unit. The upper bound and lower bound of $p_{\text{AC}}^{i,t}$ satisfy the following constraints:
\begin{equation}
     \underline{P}_{\text{AC}}^{i} \leq p_{\text{AC}}^{i,t} \leq \overline{P}_{\text{AC}}^{i},
\end{equation}
where $\underline{P}_{\text{AC}}^{i}$ and $\overline{P}_{\text{AC}}^{i}$ represent the lower bound and upper bound for the HVAC power demand of building $i$. 

The concept of discomfort cost is introduced to measure the dissatisfaction with the thermal environment inside buildings over $\mathcal{H}_{t}$, presented as 
\begin{equation}
     \mathbf{C}_{\text{AC}}^{t} = \beta^{i} \sum_{i \in \mathcal{I}} (T_{\text{InB}}^{i,t} - T_{\text{InB}}^{i,\text{Pref}})^{2},\notag
\end{equation}
where $\beta^{i}$ is the sensitive coefficient and $T_{\text{InB}}^{i,\text{Pref}}$ is the preferred indoor temperature for residential building $i$. To maintain occupant comfort, the indoor temperature should be kept within the range of $[\underline{T}^{i}, \overline{T}^{i}]$ throughout all time slots, which can be presented in the following constraint 
\begin{equation}
     \underline{T}^{i} \leq T_{\text{InB}}^{i,t} \leq \overline{T}^{i}.
\end{equation}

\subsubsection{Power Supply}
The power supply to Community $i$ comprises the local renewable power (such as solar PV) $p_{\text{renew}}^{i,t}$, power export from the local market through energy trading $p_{\text{ETb}}^{i,t} \geq 0$, and the grid power $p_{\text{grid}}^{i,t}$. These sources of power supply satisfy
\begin{align}\label{eq5}
    0 \leq p_{\text{renew}}^{i,t} \leq \overline{RE}^{i,t}, \\
    0 \leq p_{\text{ETb}}^{i,t} \leq \overline{P}_{\text{ETb}}^{i}, \label{eq5_2}\\
    0 \leq p_{\text{grid}}^{i,t} \leq \overline{G}^{i},
\end{align}
where $\overline{RE}^{i,t}$, $\overline{P}_{\text{ETb}}^{i}$, and $\overline{G}^{i}$ are the maximum available amount of renewable power, power purchased from the local market, and power purchased from the grid for Community $i$ in time slot $t$. This study assume that the local renewable can only supply to the local residential building $i$, leading to the following constraint 
\begin{equation}
    P_{\text{Bd,load}}^{i,t} + p_{\text{AC}}^{i,t} - p_{\text{renew}}^{i,t} \geq 0. \label{eq_pv}
\end{equation}

This work assesses two types of tariffs for all residential communities as described below.
\begin{enumerate}
    \item Under the two-part tariff pricing (TPT), the electricity cost of all residential communities over $\mathcal{H}_{t}$ can be written as 
    $$
    \mathbf{C}^{t}_{\text{grid}} = 
    \begin{cases} 
        \pi_g\sum_{i\in \mathcal{I}} p_{\text{grid}}^{i,t} \Delta t & \\
        \pi_g\sum_{i\in \mathcal{I}} p_{\text{grid}}^{i,t} \Delta t + \pi_{\text{peak}}\max p_{\text{grid}}^{i,t}, & \text{if } p_{\text{grid}}^{i,t} = \max p_{\text{grid}}^{i,t},
    \end{cases}
    $$
    where $\pi_\text{g}$ (in $\$$kWh) and $\pi_{\text{peak}}$ (in $\$$kW) are energy price and peak price, respectively. TPT aims to promote peak shaving among consumers.
    \item Under the time-of-use pricing (TOU), the electricity cost of all residential communities over $\mathcal{H}_{t}$ can be written as 
    $$
    \mathbf{C}^{t}_{\text{grid}} = \sum_{i\in \mathcal{I}} \pi_g^t p_{\text{grid}}^{i,t} \Delta t,
    $$
    where $\pi_g^t$ is the TOU price, usually specifying electricity rates during peak, off-peak, and shoulder hours \citep{ausgrid}.
\end{enumerate}

\subsubsection{Power Export and Import of Residential Community}
After introducing the components, such as the EVPL, building, and power supply, in the residential community, this work takes Community $i$'s perspective and denote the export power from Community $i$ as $p_{\text{Rc,Ex}}^{i,t}$, which satisfies the following constraints:
\begin{equation}\label{eq23}
    p_{\text{Rc,Ex}}^{i,t} = p_{\text{EV,Ex}}^{i,t} - (P_{\text{Bd,load}}^{i,t} + p_{\text{AC}}^{i,t} - p_{\text{renew}}^{i,t}).
\end{equation}

Note that a positive value of $p_{\text{Rc,Ex}}^{i,t}$ indicates that the power exported by EVPL and PV supply in Community $i$ are sufficient to meet the building's demand, such that no power is purchased from the grid and local energy market, i.e., $p_{\text{grid}}^{i,t} = 0$ and $p_{\text{ETb}}^{i,t} = 0$. Moreover, the surplus power available in Community $i$ can be exported to the main grid by performing V2G denoted as $p_{\text{V2G}}^{i,t}$ and the local energy market by performing energy trading, denoted as $p_{\text{ETs}}^{i,t}$. Further details regarding $p_{\text{V2G}}^{i,t}$ and $p_{\text{ETs}}^{i,t}$ will be provided in Section~\ref{valuestream}. On the contrary, a negative value of $p_{\text{Rc,Ex}}^{i,t}$ indicates that Community $i$ needs to import power, e.g., purchase power from the grid denoted as $p_{\text{grid}}^{i,t}$ and from the local energy market denoted as $p_{\text{ETb}}^{i,t}$ in order to meet the demand in time slot $t$.

\subsection{Distribution Network Model}
The residential communities are interconnected with each other and also connected to the grid on a distribution network. Given the widespread recognition and application of linearized models for distribution networks, such as in \citep{zhong2018topology} and \citep{wang2015decentralized}, this work also adopts a linearized model, as shown below.
\begin{align}
    & p^{i+1,t} = p^{i,t} + p_{\text{Rc,Ex}}^{i,t}\label{eq16},\\
    & \underline{P}^{i}_{\text{DN}} \leq p^{i,t} \leq \overline{P}^{i}_{\text{DN}}\label{eq18},\\
    & q^{i+1,t} = q^{i,t} - Q_{\text{Load}}^{i,t}\label{eq19},\\
    & \underline{Q}^{i}_{\text{DN}} \leq q^{i,t} \leq \overline{Q}^{i}_{\text{DN}}\label{eq20},\\
    & v^{i+1,t} = v^{i,t} - \left( R^{i+1} p^{i+1,t} + X^{i+1} q^{i+1,t} \right) / V^{0,t}\label{eq21},\\
    & \underline{V}^{i} \leq v^{i,t} \leq \overline{V}^{i}\label{eq22}.
\end{align}

This model is essential for imposing constraints on the active power, reactive power, and voltage within the distribution network. Here, $p^{i,t}$ represents the active power from node $i-1$ to node $i$ at $t$, and $p_{\text{Rc,Ex}}^{i,t}$ denotes the active power export (or the negative active load) from Community $i$ at node $i$. The maximum and minimum active power flows between nodes are denoted as $\overline{P}^{i}_{\text{DN}}$ and $\underline{P}^{i}_{\text{DN}}$, respectively. The reactive power between nodes $i-1$ and $i$ at $t$ is denoted by $q^{i,t}$, and $Q_{\text{Load}}^{i,t}$ is the reactive load at node $i$ at $t$. The limits for reactive power are specified by $\overline{Q}^{i}_{\text{DN}}$ and $\underline{Q}^{i}_{\text{DN}}$. The voltage at node $i$ is given by $v^{i,t}$, with $V^{0,t}$ representing the voltage at node $0$ at time $t$. The resistance and reactance of the branch between node $i-1$ and node$i$ are denoted by $R^i$ and $X^i$, respectively. Voltage levels must be maintained between the nominal range from $\underline{V}^{i}$ to $\overline{V}^{i}$.

\section{V2X Value-Stacking Optimization Problem}\label{Optimal Problem Formulation}
\begin{figure}[!t]
\centering
\includegraphics[width=1.0\linewidth]{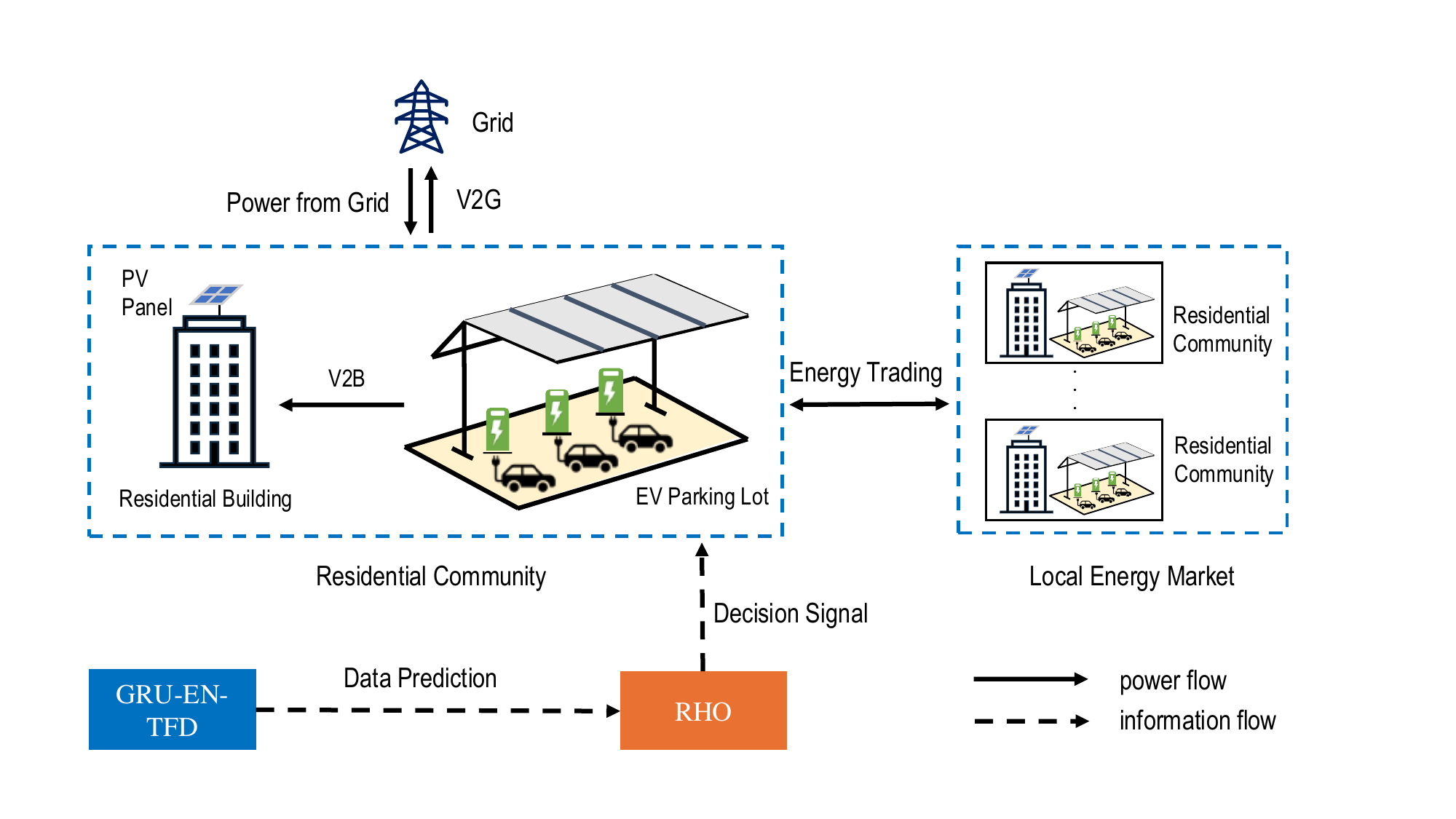} 
\caption{The framework of the proposed V2X value-stacking, including V2B, V2G, and energy trading.}
\label{fig: value-stacking}
\end{figure}

This section introduces our V2X value-stacking optimization problem, including three value streams: V2B, V2G, and energy trading. The V2X value-stacking scenario is depicted in Fig.~\ref{fig: value-stacking}, featuring each residential community equipped with a solar PV system, an EVPL, residential building load, and an energy management system. Here, EVs parked in the EVPL have the capability to discharge their batteries for V2B, V2G, and energy trading activities. Additionally, EVs can recharge using power from the grid or local energy market. This work proposes a dynamic RHO strategy that leverages prediction information. Specifically, in each time slot $t$, the GRU-EN-TFD model forecasts the building load, PV generation, and EV arrivals for the upcoming period from the next time slot $t+1$ up to the end of the day $H$ for Community $i$. Based on these predictions, this study develop a dynamic RHO approach to determine the optimal charge/discharge decisions of EVs for the V2X value-stacking over the operational horizon $[t, H]$. Residential communities implement the real-time decision in each time slot $t$ and continue this process until all decisions are made up to $H$, marking the conclusion of each daily operation. In the following, this work will outline the value streams from V2B, V2G, and energy trading, present forecasting models for building load, PV generation, and EV arrivals, and detail the dynamic RHO formulation.

\subsection{Value-Stacking Model}\label{valuestream}

In this part, three value streams, including V2B, V2G, and energy trading among residential communities, are presented. Then, this work models residential communities' V2X value stacking in multiple value streams.

\subsubsection{Value Stream of V2B}\label{V2B}
In V2B, EVs can charge their batteries when tariffs are low and then supply power to support inflexible loads and HVAC usage in residential buildings, thereby reducing peak demand and energy costs. Specifically, during hot weather, EVs' batteries can supply power to the building's air conditioning during peak demand, lowering HVAC usage and reducing stress on the grid. Conversely, during colder weather, it can help in heating systems. EV discharge from EVPL $i$ for performing V2B in time slot $t$ is denoted as $p_{\text{V2B}}^{i,t}$ with the following constraint
\begin{equation}\label{eq9}
    0 \leq p_{\text{V2B}}^{i,t} \leq \overline{P}_{\text{V2B}}^{i},
\end{equation}
where $\overline{P}_{\text{V2B}}^{i}$ denotes the maximum allowed V2B power from EVPL $i$. It is important to note that this work prioritizes supplying power from the EVPL $p_{\text{EV,Ex}}^{i,t}$ to local residential buildings for V2B operations, when the EVPL can export power, i.e., $p_{\text{EV,Ex}}^{i,t} > 0$. More specifically, if $p_{\text{EV,Ex}}^{i,t} \leq P_{\text{Bd,load}}^{i,t} + p_{\text{AC}}^{i,t} - p_{\text{renew}}^{i,t}$, the V2B power from the EVPL is insufficient to fulfill the deficit of building net demand $P_{\text{Bd,load}}^{i,t} + p_{\text{AC}}^{i,t} - p_{\text{renew}}^{i,t}$. If $p_{\text{EV,Ex}}^{i,t} > P_{\text{Bd,load}}^{i,t} + p_{\text{AC}}^{i,t} - p_{\text{renew}}^{i,t}$, the V2B power can fulfill the deficit of building net demand and the excess power can then be exported from the community to perform V2G and energy trading with other communities in the local market.

\subsubsection{Value Stream of V2G}\label{V2G}
In V2G, the exported power of Community $i$ can be sold to the grid to make a profit. Note that the exported power of Community $i$ comes from EVPL $i$, and $p_{\text{V2G}}^{i,t}$ is denoted as the V2G energy that EVPL $i$ sells to the grid. The constraint of $p_{\text{V2G}}^{i,t}$ is as follows
\begin{equation}\label{eq10}
    0 \leq p_{\text{V2G}}^{i,t} \leq \overline{P}_{\text{V2G}}^{i},
\end{equation}
where $\overline{P}_{\text{V2G}}^{i}$ is the maximum V2G from EVPL $i$ at $t$.

This study let $\pi_{\text{V2G}}^t$ be the time-varying V2G price in time slot $t$ and thus can obtain the revenue of V2G at $t$ as $\pi_{\text{V2G}}^{t} p_{\text{V2G}}^{i,t}$. The total revenue of V2G for all EVPLs over $\mathcal{H}_{t}$ is 
$$
\mathbf{R}_{\text{V2G}}^{t} = \sum \nolimits_{i\in\mathcal{I}} \pi_{\text{V2G}}^{t} p_{\text{V2G}}^{i,t} \Delta t.
$$

\subsubsection{Value Stream of Energy Trading in Local Market}\label{ET}
Residential communities can engage in energy trading using EV batteries in EVPLs, operating within the local market under dynamic pricing $\pi_{\text{ET}}^{t}$, which is calculated based on the Mid-Market Rate. This rate is determined as the midpoint between the retail buying and selling prices \citep{long2017peer}, ensuring that the selling price does not exceed the buying price. This work separates buying and selling in energy trading, denoted as $p_{\text{ETb}}^{i,t}$ and $p_{\text{ETs}}^{i,t}$, respectively. 
The constraints for buying have been specified in Eq. \eqref{eq5_2} and the constraints for selling in energy trading are 
\begin{align}
    0 \leq p_{\text{ETs}}^{i,t} \leq \overline{P}_{\text{ETs}}^{i}, \label{eq14}
\end{align}
where $\overline{P}_{\text{ETs}}^{i}$ denotes the upper bound for selling energy in the local market for Community $i$. Note that buying and selling decisions as part of the energy trading do not happen at the same time for the same community, and this work will discuss the constraints later in Eq. \eqref{eq25}. The energy selling and purchasing to and from the local energy market in time slot $t$ should be balanced and the constraint is 
\begin{equation}\label{eq24}
    \sum_{i \in \mathcal{I}}p_{\text{ETs}}^{i,t} =  \sum_{i \in \mathcal{I}}p_{\text{ETb}}^{i,t}.
\end{equation}
Therefore, the net cost of energy trading, $\mathbf{C}_{\text{ET}}^{t}$ ,for the entire system in time slot $t$ is zero. 

\subsubsection{Multiple Value Streams of V2X}\label{VS}
Various value streams that are enabled by V2X technologies are explored, including V2B, V2G, and energy trading among residential communities in the local market. The residential communities optimize the EV charging/discharging scheduling in their EVPLs to capitalize on these multiple value streams. The exported power from Community $i$ satisfies the following constraints:
\begin{equation}\label{eq25}
   \left \{
    \begin{aligned}
        & p_{\text{V2G}}^{i,t} + p_{\text{ETs}}^{i,t} - p_{\text{Rc,Ex}}^{i,t} = 0,~\text{if}~ p_{\text{Rc,Ex}}^{i,t} > 0, \\
        & p_{\text{grid}}^{i,t} + p_{\text{ETb}}^{i,t} =0,~\text{if}~ p_{\text{Rc,Ex}}^{i,t} > 0,\\
        & p_{\text{V2G}}^{i,t} + p_{\text{ETs}}^{i,t} = 0,~ \text{if}~ p_{\text{Rc,Ex}}^{i,t} \leq 0,\\
        & p_{\text{ETb}}^{i,t} + p_{\text{grid}}^{i,t} + p_{\text{Rc,Ex}}^{i,t} = 0,~ \text{if}~ p_{\text{Rc,Ex}}^{i,t} \leq 0.
    \end{aligned}
   \right.
\end{equation}
To implement conditional constraint in Eq.~\eqref{eq25}, this work reformulates it into Eq.~\eqref{eq26} using the Big-M method via introducing a very large number $M$ and a binary variable $y^{i,t} \in \{0,1\}$ for Community $i$ in time slot $t$. Note that $y^{i,t} = 0$ represents $p_{\text{Rc,Ex}}^{i,t} > 0$, and $y^{i,t} = 1$ indicates $p_{\text{Rc,Ex}}^{i,t} \leq 0$. 
\begin{equation}\label{eq26}
   \left \{
    \begin{aligned}
        & 0 \leq  p_{\text{V2G}}^{i,t} + p_{\text{ETs}}^{i,t} - p_{\text{Rc,Ex}}^{i,t} \leq My^{i,t}, \\
        & 0 \leq p_{\text{grid}}^{i,t} + p_{\text{ETb}}^{i,t} \leq My^{i,t},\\
        & 0 \leq p_{\text{V2G}}^{i,t} + p_{\text{ETs}}^{i,t} \leq M(1-y^{i,t}), \\
        & 0 \leq p_{\text{ETb}}^{i,t} + p_{\text{grid}}^{i,t} + p_{\text{Rc,Ex}}^{i,t} \leq M(1-y^{i,t}).
    \end{aligned}
   \right.
\end{equation}

\subsection{Transformer-based Forecasting Model}\label{Transformer}
\begin{figure}[!t]
\centering
\includegraphics[width=1.0\linewidth]{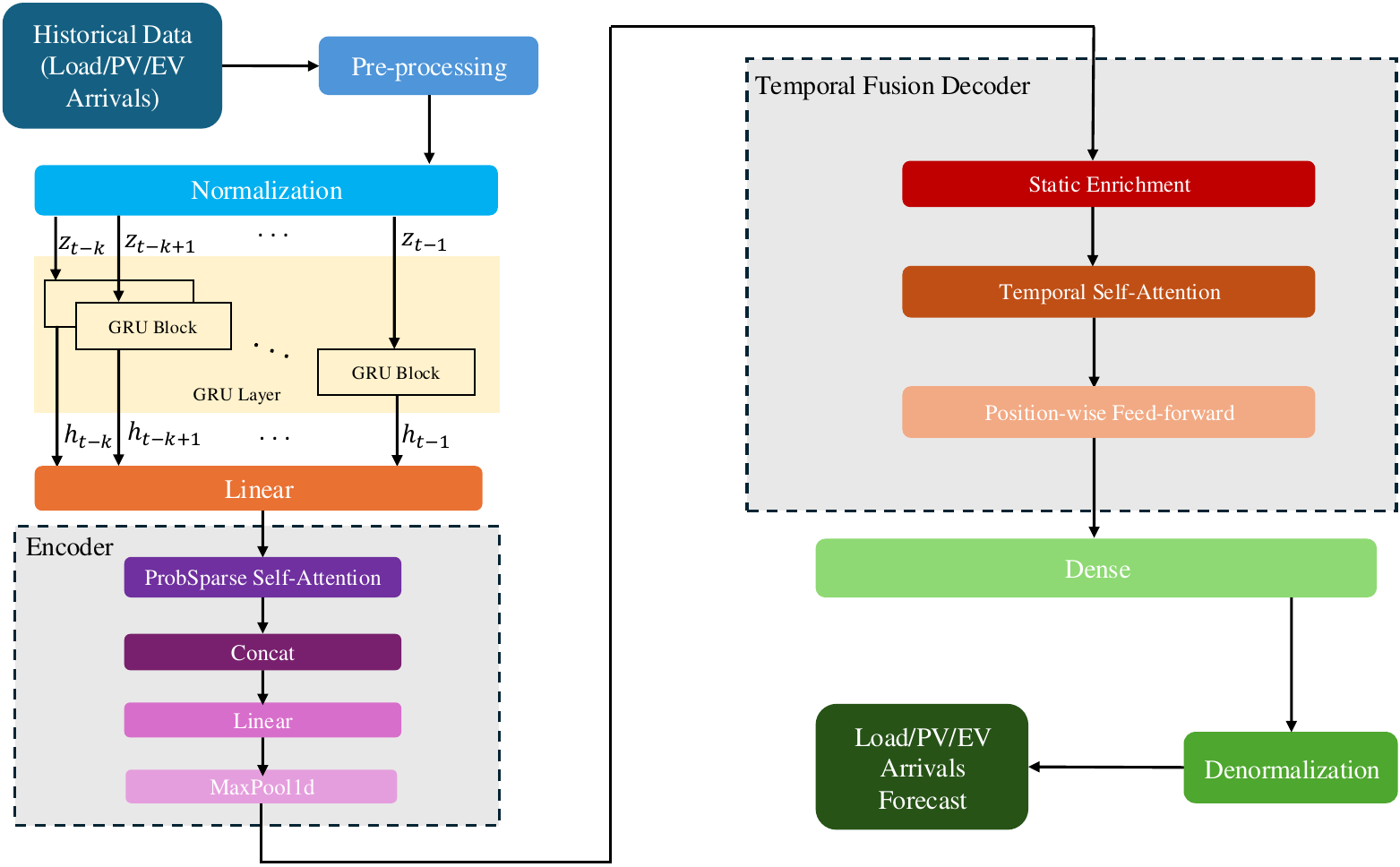} 
\caption{Framework of the GRU-EN-TFD model for forecasting building demand, PV generation, and EV arrivals.}
\label{fig: GRU-EN-TFD}
\end{figure}

To perform the RHO for V2X value-stacking, forecasting critical future information, including residential building load, solar PV generation, and EV arrivals, is essential. This work develops a Transformer-based forecasting model called GRU-EN-TFD, as depicted in Fig.~\ref{fig: GRU-EN-TFD}. Different from the traditional LSTM, Gated Recurrent Unit (GRU), or Transformer models, the proposed GRU-EN-TF model integrates a Gated Recurrent Unit (GRU) layer, Encoder, and Temporal Fusion Decoder. 

\subsubsection{GRU Layer}
The GRU module is a simplified variant of Recurrent Neural Network (RNN) that incorporates gating mechanisms—reset and update gates—to effectively manage time-series input \citep{wang2024high}. In this model, the GRU layer processes input time-series data to extract temporal features, and the hidden layer outputs from each node serve as inputs for the encoder.

\subsubsection{Encoder} 
The encoder consists of Self-Attention, Dense, and Temporal Pooling. Self-Attention is a computationally efficient variant of attention mechanisms and suitable for long sequences \citep{mo2024powerformer}. It selectively focuses on a subset of important elements rather than the entire sequence. The Temporal Pooling can aggregate information over time, reducing the sequence length while retaining important features.

\subsubsection{Temporal Fusion Transformers Model}
The Temporal Fusion Transformer (TFT) is employed for multivariate time-series forecasting for its ability to capture temporal correlations among input features and assign importance weights to each feature \citep{zheng2023interpretable}. The temporal fusion decoder is part of the TFT architecture, processing the outputs from earlier layers and combining them with the learned representations of temporal dependencies to generate accurate forecasts.

\color{black}
\subsubsection{Forecasting of Building Demand, PV Generation, and EV Arrivals}
To forecast residential building demand, PV generation, and EV arrivals, this work develops three separate GRU-EN-TFD networks. Each network is designed to forecast hourly building load, PV generation, and EV arrivals, respectively. The process begins with the preprocessing the input data, which includes historical records of building load, PV generation, and EV arrivals, respectively. The normalization is performed to scale the input features (e.g., building load demand, PV generation, and EV arrivals) to a common scale without distorting the differences in the ranges of values. This step is crucial for promoting faster convergence during training. In this work, the Min-Max scaling technique is applied. 

In the training process, this work uses $\{z_{t-k},z_{t-k+1},...,z_{t-1}\}$ to denote the past $k$ time-slot sequence as the input data, where $z_{t} \in \mathbf{R}^{n} $ is an $n$-dimensional vector of real values in the $t$-th time slot. In the GRU layer, it processes the sequential input data, capturing the temporal dependencies within the previous 48 time-slot input data. In the Encoder, the ProbSparse Self-Attention selectively focuses on the most significant parts of the input sequence, enhancing the model's ability to capture crucial temporal events without being overwhelmed by data scale. Each output from the attention mechanism is processed via a dense (fully connected) layer, which allows for nonlinear processing of the features. A pooling layer aggregates the temporal information, reducing the sequence length while preserving essential features, thus preparing the data for decoding. In the Temporal Fusion Decoder, static data include invariant features, providing contextual grounding, which is essential for accurate forecasting in diverse settings. Temporal Self-Attention and Position-Wise Feed-Forward Network process the enriched features, further refining the model's ability to predict future energy consumption by considering both historical and contextual data. The final output of the decoder is passed through a dense layer that maps the hidden states to the forecasted data for the next 24 time periods, aligning the output dimensions with the prediction horizon. The prediction outputs are in the normalized scale used during training. Denormalization is applied as the process of converting normalized predictions back to the original scale of the data. This step is crucial for making the predictions interpretable and useful in a real-world context.

This work evaluates the performance of the proposed GRU-EN-TFD model. For comparison, we introduce the LSTM network, a widely used type of RNNs that have demonstrated significant promise in time-series forecasting, as noted in previous studies, such as \citep{kong2017short}.

\color{black}
\subsection{Value-Stacking Rolling Horizon Optimization Problem}\label{valuestacking}
\begin{figure}[!t]
    \centering
    \includegraphics[width=1.0\linewidth]{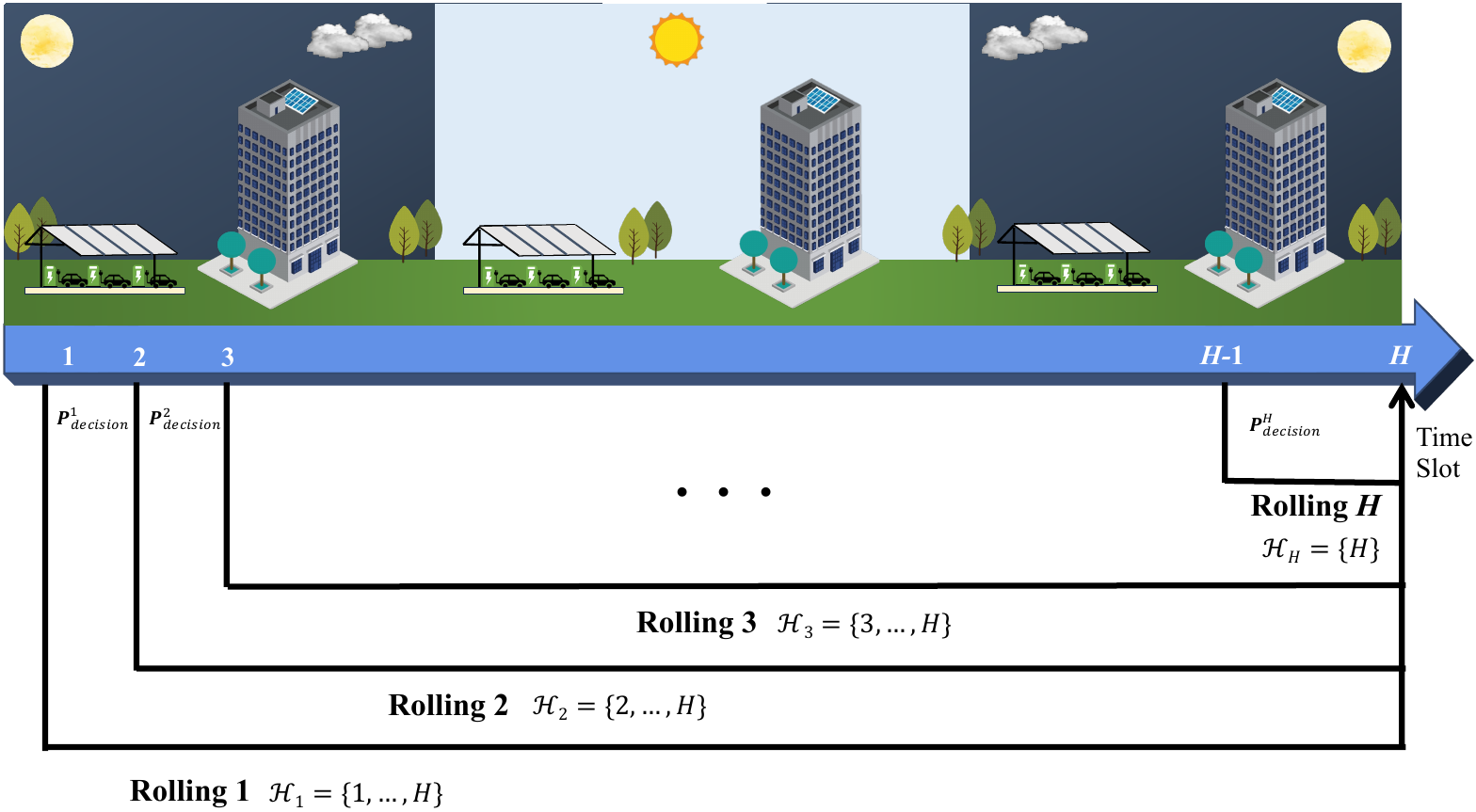}
    \caption{Shrinking time window for the dynamic rolling-horizon EV value-stacking optimization.}
\label{fig: dynamic window}
\end{figure}

After presenting the system model in Section~\ref{Model system}, V2X value stream models in Section~\ref{valuestream}, and the Transformer-based forecast model for building load, PV generation, and EV arrivals in Section~\ref{Transformer}, the dynamic RHO is now formulated for the real-time V2X value-stacking problem in this section. This formulation co-optimizes three value streams, including V2B, V2G, and energy trading among residential communities, to effectively manage the charging and discharging of EVs parked in EVPLs. This work introduces $\mathbf{C}_{\text{total}}^{t}$ to represent the total cost across all residential communities. The objective of the dynamic RHO at time $t$ is to minimize the total cost over $\mathcal{H}_{t}$ by optimizing decision variables such as $p_{\text{grid}}^{i,t}$, $p_{\text{renew}}^{i,t}$, $p_{\text{AC}}^{i,t}$, $T_{\text{InB}}^{i,t}$, $p_{\text{EVc}}^{u,t}$, $p_{\text{EVd}}^{u,t}$, $x^{u,t}$, $b_{\text{EV}}^{u,t}$, $p_{\text{V2B}}^{i,t}$, $p_{\text{V2G}}^{i,t}$, $p_{\text{ETs}}^{i,t}$, $p_{\text{ETb}}^{i,t}$, $p_{\text{EV,Ex}}^{i,t}$, $p_{\text{Rc,Ex}}^{i,t}$, $y^{i,t}$, $p^{i,t}$, $q^{i,t}$, and $v^{i,t}$. It is crucial to note that unlike an offline optimization problem with a fixed time horizon, the value-stacking RHO employs a dynamic operational horizon, i.e., $\mathcal{H}_{t}$, which adapts according to the dynamic time window illustrated in Fig. \ref{fig: dynamic window}. Specifically, on each operational day, the horizon $\mathcal{H}_{t}$ gradually decreases as time $t$ progresses until the end of the day at $t=H$.

Therefore, the optimization problem for V2X value stacking over the dynamic rolling-horizon $\mathcal{H}_{t}$ is formulated as follows:
\begin{equation}
\begin{aligned} \label{S1}
\min~ & \mathbf{C}_{\text{total}}^{t} =  \sum_{t \in \mathcal{H}_{t}} (\mathbf{C}^{t}_{\text{grid}} + \mathbf{C}_{\text{battery}}^{t}+ \mathbf{C}_{\text{AC}}^{t} -\mathbf{R}_{\text{V2G}}^{t})\\
\text{s.t.} ~& \eqref{eq1}-\eqref{eq24}~\text{and}~ \eqref{eq26}.
\end{aligned}
\end{equation}

This formulation represents a mixed-integer quadratic programming (MIQP) problem that operates over a dynamic rolling time window. In time slot $t$, taking the forecasting information shown in Section~\ref{Transformer}, the dynamic RHO problem solves the optimal EV charging/discharging to minimize the total cost over the time-varying horizon $\mathcal{H}_{t}$. Only the decisions at $t$ are executed. The complete process is presented in Algorithm~\ref{Process}, where $\mathbf{P}_{\text{decision}}^{t}$ is defined as the variable vector including all the decision variables in time slot $t$ in \eqref{S1}.

In each time slot $t$, the V2X value-stacking optimization determines the optimal decisions and generates schedules for all residential communities, based on the predicted building demand, PV generation, and EV arrivals from time slot $t+1$ to the end of the horizon $H$. In this procedure, only the decision variables for the current time slot are executed and deployed for cost calculation. As the rolling horizon advances, the size of the dynamic window diminishes by one-time slot with each step until reaching the end of the daily operational horizon. The optimal cost derived from the value-stacking problem with perfect predictions serves as the baseline for evaluating the impact of prediction errors.

\begin{algorithm}[!t]
\DontPrintSemicolon
\caption{Dynamic Rolling-horizon Optimization}\label{Process}
\SetKwInOut{Input}{input}
\SetKwInOut{Output}{output}
\renewcommand{\Input}{\textbf{Initialization:~}}
\Input{$\text{Time slot}~t$;\newline
Dynamic time horizon $\mathcal{H}_{t}$ $\leftarrow$ $\left\{t,\dots, H\right\}$;\newline
Cost of all communities in time slot $t$, i.e., $\mathbf{C}_{\text{Total}}^{t}$;\newline}
\For{$t= 1:H$}{\text{$\bullet$ Input realized data for demand, PV generation, EV arrival}  
\text{up to time $t$, along with predicted data from the GRU-EN-TFD in} \text{time slots $[t+1, H]$ into the RHO formulated in~\eqref{S1};}

\text{$\bullet$ Solve the RHO problem in~\eqref{S1} under} \text{dynamic time horizon $\mathcal{H}_{t}$ and execute the decision } \text{$\mathbf{P}_{\text{decision}}^{t}$ in time slot $t$;}\; 

\text{$\bullet$ Calculate the total cost of all communities in current} \text{time slot $t$ based on $\mathbf{P}_{\text{decision}}^{t}$;}}

\Output{Optimal decisions $\mathbf{P}_{\text{decision}}^{t}$ over dynamic time horizon $\mathcal{H}_{t}$ and minimized RHO cost for all residential communities.}
\end{algorithm}

\section{Model Performance Evaluation}\label{evaluation}
In this section, three baseline optimization problems are introduced for comparison with our V2X value-stacking model, demonstrating its superiority over any single value stream. Additionally, this work presents three further baseline optimization problems to assess the marginal contributions of individual value streams within V2X. To explore the effects of prediction errors, this work proposes performance metrics designed to quantify the additional costs incurred due to prediction errors, utilizing a dynamic rolling-horizon optimization approach.

\subsection{Baseline Optimization Problems}\label{baseline}
To assess the benefits of each value stream, three baseline optimization problems are explored, each focusing on a single aspect: V2G, V2B, and energy trading among residential communities. These baselines are outlined as follows:
\begin{itemize}
    \item V2G Alone: Residential communities optimize V2G interactions by scheduling variables such as $p_{\text{grid}}^{i,t}$, $p_{\text{renew}}^{i,t}$, $p_{\text{AC}}^{i,t}$, $T_{\text{InB}}^{i,t}$,  $p_{\text{EVc}}^{u,t}$, $p_{\text{EVd}}^{u,t}$, $x^{u,t}$, $b_{\text{EV}}^{u,t}$, $p_{\text{EV,Ex}}^{i,t}$, $p_{\text{V2G}}^{u,t}$, $p^{i,t}$, $q^{i,t}$, and $v^{i,t}$, in which if $p_{\text{EV,Ex}}^{i,t} > 0$, $p_{\text{EV,Ex}}^{i,t} = p_{\text{V2G}}^{i,t}$.
    \item V2B Alone: residential community optimizes specific variables to minimize their costs, including $p_{\text{grid}}^{i,t}$, $p_{\text{renew}}^{i,t}$, $p_{\text{AC}}^{i,t}$, $T_{\text{InB}}^{i,t}$, $p_{\text{EVc}}^{u,t}$, $p_{\text{EVd}}^{u,t}$, $x^{u,t}$, $b_{\text{EV}}^{u,t}$, $p_{\text{EV,Ex}}^{i,t}$, and $p_{\text{V2B}}^{i,t}$, in which if $p_{\text{EV,Ex}}^{i,t} > 0$, $p_{\text{EV,Ex}}^{i,t} = p_{\text{V2B}}^{i,t}$. 
    \item Energy Trading Alone: residential communities with EVs parked in EVPL trade energy with others by scheduling variables, including $p_{\text{grid}}^{i,t}$, $p_{\text{renew}}^{i,t}$, $p_{\text{AC}}^{i,t}$, $T_{\text{InB}}^{i,t}$, $p_{\text{EVc}}^{u,t}$,  $p_{\text{EVd}}^{u,t}$, $x^{u,t}$, $b_{\text{EV}}^{u,t}$, $p_{\text{ETs}}^{i,t}$, $p_{\text{EV,Ex}}^{i,t}$, $p_{\text{ETb}}^{i,t}$, $p^{i,t}$, $q^{i,t}$, and $v^{i,t}$, in which if $p_{\text{EV,Ex}}^{i,t} > 0$, $p_{\text{EV,Ex}}^{i,t} = p_{\text{ETs}}^{u,t}$.
\end{itemize}
Additionally, to determine the marginal contribution of each value stream within our V2X value-stacking problem, this work considers three supplementary baseline scenarios. Each scenario omits one value stream from the comprehensive value-stacking model, thereby illustrating the incremental benefit of the excluded stream.

\subsection{Performance Metrics}
To assess the impact of prediction errors, the Relative Extra Cost (REC) is computed, compared to the optimal baseline cost using the following formula:
\begin{equation}\label{eq27}
    \text{REC} = \frac{\sum_{t\in \mathcal{H}} \sum_{i \in \mathcal{I}} |\widetilde{\mathbf{C}}_{\text{cost}}^{i,t} - \mathbf{C}_{\text{cost}}^{i,t}|}{\sum_{t\in \mathcal{H}} \sum_{i \in \mathcal{I}} (\mathbf{C}_{\text{cost}}^{i,t})},
\end{equation}
where $\mathbf{C}_{\text{cost}}^{i,t}$ and $\widetilde{\mathbf{C}}_{\text{cost}}^{i,t}$ represent the actual and predicted costs for Community $i$ in time slot $t$, respectively.

For the forecasting models, this work computes the Relative Errors (RE) for the predicted load, PV generation, and EV arrivals as follows:
\begin{align}
  & \text{RE}_{\text{Bd,load}} = \frac{\sqrt{\sum_{t\in \mathcal{H}} \sum_{i \in \mathcal{I}} (\widetilde{P}_{\text{Bd,load}}^{i,t} - P_{\text{Bd,load}}^{i,t})^{2}}}{\sqrt{\sum_{t\in \mathcal{H}} \sum_{i \in \mathcal{I}} (P_{\text{Load}}^{i,t})^{2}}}\label{eq28},\\
  & \text{RE}_{\text{PV}} = \frac{\sqrt{\sum_{t\in \mathcal{H}} \sum_{i \in \mathcal{I}} (\widetilde{P}_{\text{renew}}^{i,t} - p_{\text{renew}}^{i,t})^{2}}}{\sqrt{\sum_{t\in \mathcal{H}} \sum_{i \in \mathcal{I}} (p_{\text{renew}}^{i,t})^{2}}}\label{eq29}, \\ 
  & \text{RE}_{\text{EV}} = \frac{\sqrt{\sum_{t\in \mathcal{H}} \sum_{i \in \mathcal{I}} (\widetilde{N}_{\text{EV}}^{i,t} - N_{\text{EV}}^{i,t})^{2}}}{\sqrt{\sum_{t\in \mathcal{H}} \sum_{i \in \mathcal{I}} (N_{\text{EV}}^{i,t})^{2}}}\label{eq30},
\end{align}
where $\widetilde{P}_{\text{Bd,load}}^{i,t}$, $\widetilde{P}_{\text{renew}}^{i,t}$, and $\widetilde{N}_{\text{EV}}^{i,t}$ represent predicted building demand, PV generation, and EV arrivals for Community $i$ in time $t$.

\section{Simulations and Discussions}\label{Simulation and Discussion}
\begin{figure}[!b]
    \centering
    \includegraphics[width=1.0\linewidth]{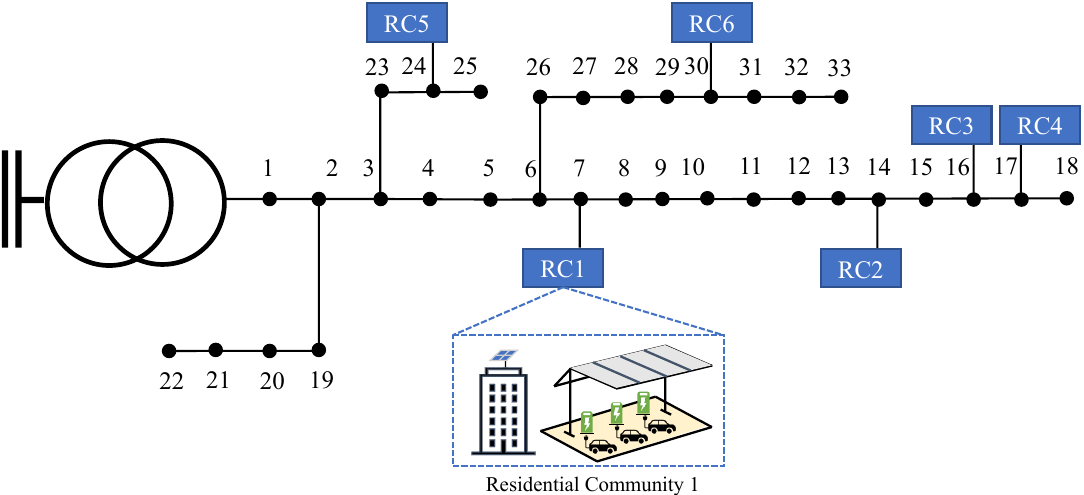}
    \caption{Modified IEEE 33-bus distribution test system with six residential communities.}
\label{fig:IEEE33bus}
\end{figure}

\begin{figure}[t]
\centering
\includegraphics[width=0.8\linewidth]{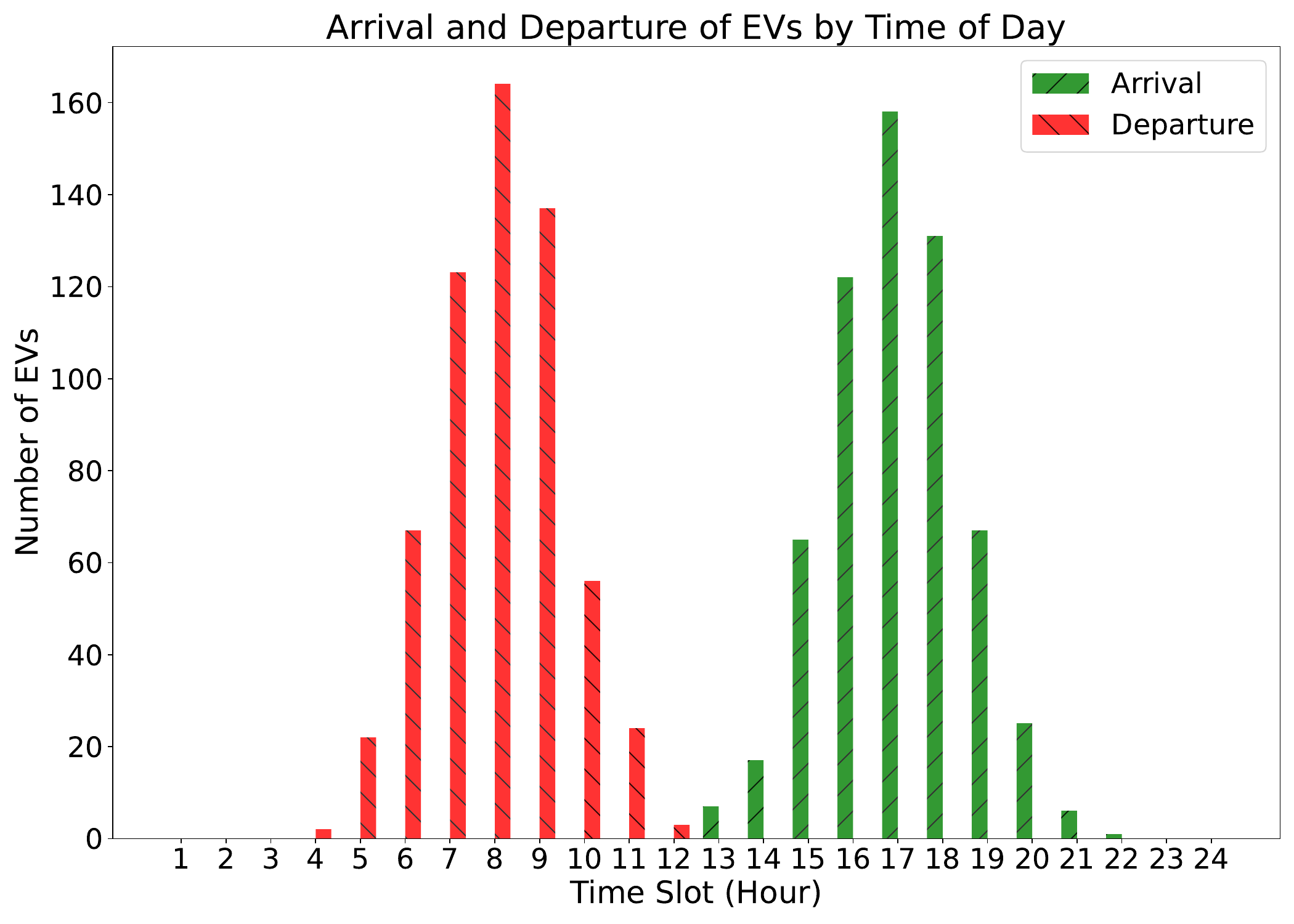} 
\caption{The arrival/departure times distribution with respect to the number of EVs.}
\label{fig: EV arrivals}
\end{figure}

This work assesses the proposed V2X value-stacking using the modified IEEE 33-bus distribution system, depicted in Fig. \ref{fig:IEEE33bus}. This system includes six residential communities connected to nodes 7, 14, 16, 17, 24, and 30, respectively. In this IEEE 33-bus distribution network, the standard values for apparent power and voltage are set at $10$ MVA and $12.66$ kV, respectively. Each community has a residential building and $50$ EVs parked in EVPL. In total, this work considers $300$ EVs. The EVs have a battery capacity of 50 kWh each, with initial stored energy varying between 20 kWh and 30 kWh. EV arrival and departure times are generated following a truncated normal distribution \citep{guldorum2022management}, as depicted in Fig.~\ref{fig: EV arrivals}. It is important to note that the generated scenarios reflect typical weekday patterns. The data reveals that EV owners generally depart from their homes in the early morning hours and return at the end of the workday. Additionally, this work considers a single central HVAC unit capable of controlling the indoor temperature for each building. The HVAC system's parameters, including heat capacity and thermal resistance, are set to $3.3$ and $1.35$ respectively. The preferences on the indoor temperature is set to $25^{\circ}$C.

For the GRU-EN-TFD model, the input data consist of 48-hour sequences of hourly building load, PV generation, and EV arrivals, respectively. The model outputs the subsequent 24-hour sequences of hourly data. An initial learning rate of 0.005 is employed, a maximum of 100 epochs, and training using the root mean square error (RMSE) metric. The Adaptive Moment Estimation (Adam) algorithm is utilized for optimization. Data are batched with a size of 32 sequences per batch, ensuring adequate sampling diversity while maintaining computational efficiency. Given the sequential nature of the data, shuffling is not employed to preserve the temporal ordering, which is crucial for learning accurate temporal dynamics.

\begin{figure}[!t]
\centering
\includegraphics[width=1.0\linewidth]{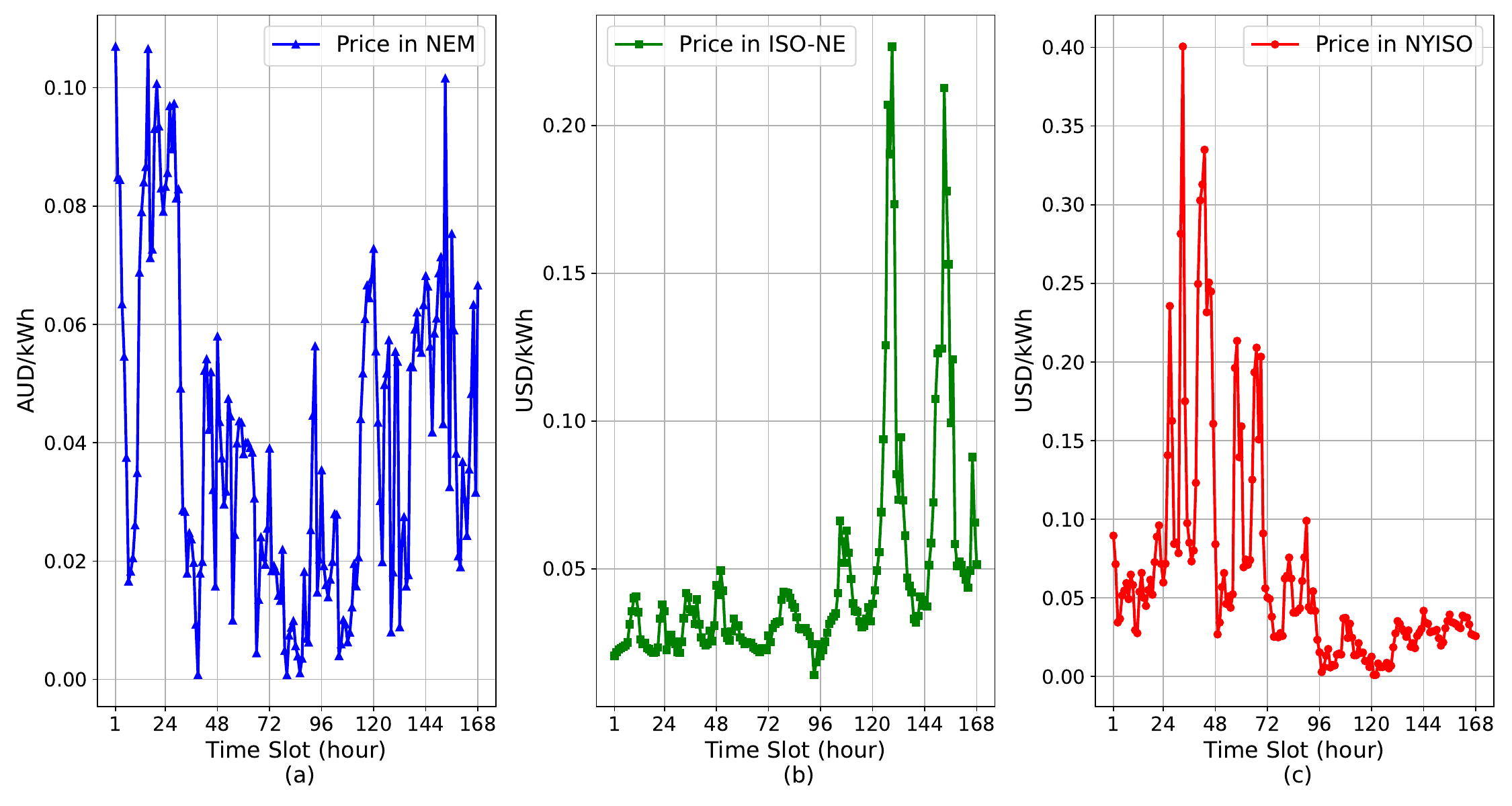} 
\caption{Illustration of electricity prices in NEM, NYISO, and ISO-NE markets over one week from December 24th to 30th, 2022.}
\label{fig: electricity prices}
\end{figure}

The model uses real market data from entities like NEM, ISO-NE, and NYISO. One week of data (from December 24th to 30th, 2022) is taken as an example to illustrate the hourly price dynamics over a week in three different markets, as shown in Fig.~\ref{fig: electricity prices}. It is important to note that we retain the local currencies of the original electricity prices from each market without converting them to a common unit. Additionally, this work evaluates two tariff structures for residential communities: TOU and TPT. Under the TOU tariff, the off-peak price is AU\$ 0.20, and the peak price is AU\$ 0.32 \citep{tariff}. For the TPT plan, the energy price is AU\$ 0.2, with a peak price of AU\$ 0.80 \citep{tariff}. In the TOU plan for NYISO and ISO-NE markets, the off-peak, shoulder, and peak prices are US\$ 0.02, US\$ 0.11, and US\$ 0.32, respectively. The value-stacking problem is efficiently solved using the Gurobi solver with JuMP/Julia on a commodity PC equipped with an M1 Pro Chip and 32GB of memory \citep{Gurobi}, \citep{JuMP}.

In this section, this work first analyzes the value-stacking performance through an analysis based on perfect prediction, as detailed in \ref{VS with perfect para}. Additionally, this work illustrates the decision-making processes in EV value stacking against three baselines under TOU and TPT tariffs, presented in \ref{decision of EV}. Furthermore, the impact of V2B on HVAC usage is elucidated through a comparison between scenarios with and without V2B integration, as discussed in \ref{HVAC with and without V2B}. The accuracy of the GRU-EN-TFD model in forecasting residential building demand, PV generation, and EV arrivals is evaluated in Section \ref{forecast assessment}. Finally, the effects of prediction errors on the performance of value stacking are analyzed in Section \ref{impact of prediction error on VS}.

\color{black}
\subsection{Performance of Value Stacking Under Perfect Prediction}\label{VS with perfect para}
\begin{figure}[!t]
    \centering
    \includegraphics[width=0.8\linewidth]{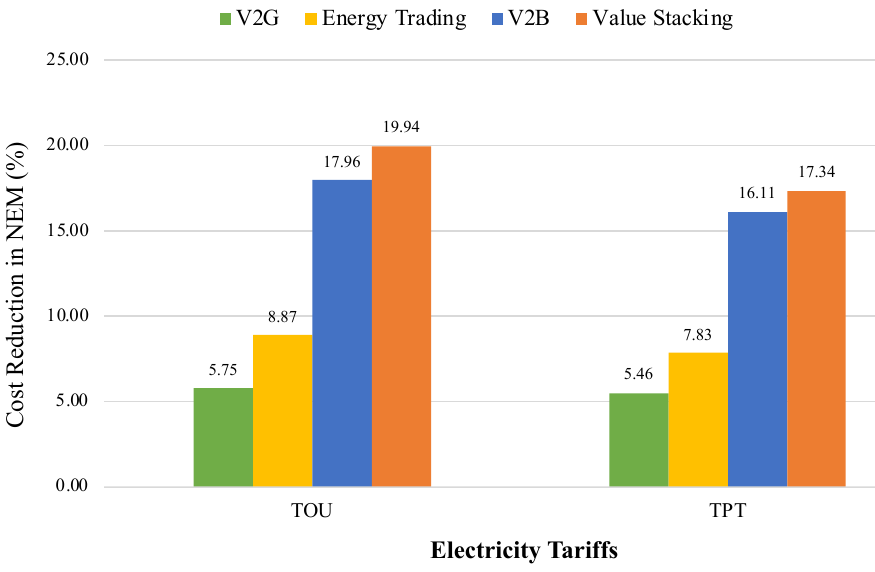}
    \caption{Cost reductions for residential communities under two tariffs across four scenarios.}
\label{fig:cost reduction}
\end{figure}

This work initially assesses the efficacy of V2X value-stacking, focused on cost savings for all residential communities, assuming accurate forecasts of building demand, PV generation, and EV arrivals. The evaluation covers a 24-hour simulation horizon with average daily costs calculated over a week (using real data from December 24th to 30th, 2022), and the findings are depicted in Fig. \ref{fig:cost reduction}. The vertical axis displays the percentage of total cost reductions achieved by value stacking and three baselines—V2G, V2B, and energy trading among communities—relative to a basic baseline where only EV charging is optimized without allowing for discharging. The two groups of bars on the horizontal axis represent two different tariffs: TOU and TPT, applied to the residential communities. Under TOU pricing, value stacking leads to a cost reduction of 19.94\%. Out of the three value streams, V2G alone yields the smallest reduction, at 5.75\%. The advantage of V2B surpasses that of energy trading within the local market, since all residential communities face identical peak load hours, thus limiting opportunities for energy trading but enhancing V2B's role in reducing costs. Under TPT, the cost reduction from value stacking is 17.34\%, which is less than under TOU. For the baselines, V2G achieves a 5.46\% reduction, slightly less than its performance under TOU. V2B still outperforms energy trading among communities, though its effectiveness drops to 16.11\% under TPT. 

Under both TOU and TPT, value stacking achieves the most cost reduction compared to the baselines limited to individual value streams, indicating that no single value stream dominates. The benefit of value stacking stems from the synergy among different value streams, which fully unlocks the value of EV batteries and their operational flexibility (charging and discharging) in interacting with multiple entities, such as the grid, buildings, and other EVs, over time.

Compared to TOU, TPT leads to a lower overall cost reduction from V2X value stacking. This decline can be attributed to the inherent flexibility of TPT, which enables peak shaving and cost reduction even in the baseline cases. In addition, the varying times of peak loads may not align with the cost-reduction opportunities offered by V2X value stacking, which further decreases the cost reduction under TPT. 

Among the three value streams, V2G alone yields the smallest cost reduction. Notably, the advantage of V2B surpasses that of energy trading within the local market. This is primarily because all residential communities experience similar peak load hours, which limits the opportunities for energy trading but enhances V2B's role in cost savings. This work further explores the individual contributions of each value stream within the V2X framework to understand their specific impacts on cost efficiency.

\begin{figure}[!t]
    \centering
    \includegraphics[width=0.8\linewidth]{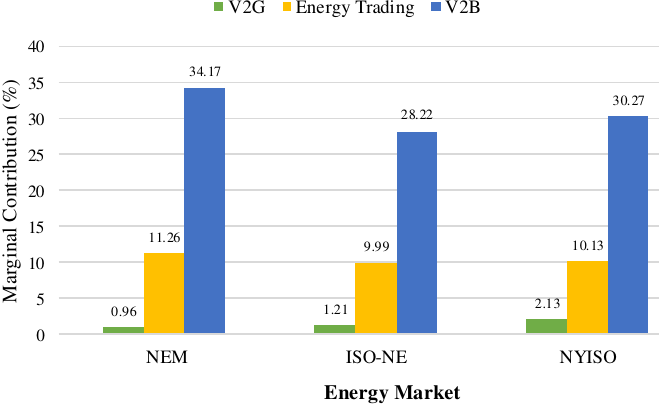}
    \caption{Marginal contribution of various value streams in three energy markets.}
\label{fig:Marginal}
\end{figure}

To determine the marginal contributions of each value stream within our value-stacking model discussed in Section~\ref{valuestacking}, this work assesses their impact across different markets: specifically, NEM, ISO-NE, and NYISO. The marginal contribution is calculated as the difference in cost reductions between comprehensive value stacking and baseline scenarios that omit a specific value stream. As shown in Fig.~\ref{fig:Marginal}, V2B significantly outperforms V2G and energy trading across all three markets, with contributions reaching 34.17\%, 28.22\%, and 30.27\%, respectively. This demonstrates that V2B offers the most substantial cost savings when implemented in these markets. Furthermore, local energy trading enhances the total system cost efficiency by optimizing energy allocations within the value-stacking model, leading communities to prefer selling power locally at peak consumption times over supplying power for V2G, which leads to a relatively modest contribution from V2G across the three markets.

The results reveal that V2B's marginal contribution is consistently greater than those of V2G and energy trading in all three markets. The primary reason is that the high purchase prices incentivize communities to store energy in EV batteries during off-peak periods for use during peak times via V2B. In contrast, V2G and energy trading often cannot offer greater economic benefits at the same time.

\begin{figure}[t]
    \centering
    \includegraphics[width=0.8\linewidth]{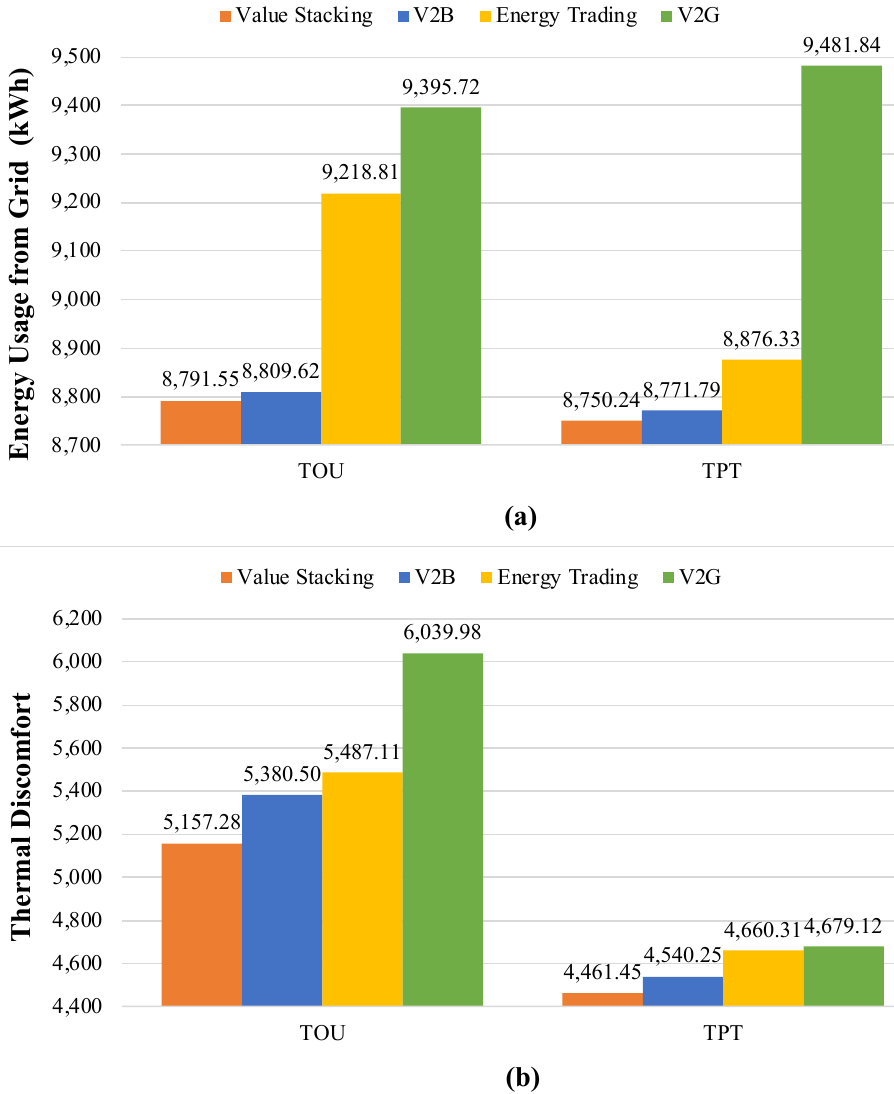}
    \caption{Comparison of two tariffs in energy usage from grid and thermal discomfort in four scenarios.}
\label{fig:Thermal Discomfort}
\end{figure}

This work also assesses the effectiveness of V2X value stacking in managing electricity consumption from the grid and addressing thermal discomfort within the framework of TOU and TPT tariffs, as depicted in Fig.~\ref{fig:Thermal Discomfort}. In Fig.~\ref{fig:Thermal Discomfort} (a), the vertical axis illustrates the total electricity consumption from the grid (kWh) for value stacking and three baseline scenarios: V2G, V2B, and energy trading among communities. Correspondingly, in Fig.~\ref{fig:Thermal Discomfort} (b), the vertical axis represents thermal discomfort across six communities. The two groups of bars on the horizontal axis represent the TOU and TPT tariffs. In the case of TPT, value stacking demonstrates the most favorable outcomes in terms of both energy usage and thermal discomfort. V2B outperforms energy trading among residential communities in both electricity consumption and thermal discomfort. This disparity might be attributed to V2B's direct power supply to buildings, effectively reducing building load usage, whereas energy trading relies on the availability of surplus power within neighboring communities for purchase within the local energy market. V2G contributes less to building load usage and thermal discomfort, as surplus power from EVPL is exported to the grid for profit. Conversely, under TOU pricing, value stacking shows more electricity consumption from the grid compared to TPT, but the HVAC discomfort is higher than under TPT. This could be explained by that, during peak hours in TOU, communities reduce HVAC usage but sacrifice thermal comfort because the electricity cost during peak hours is higher than the cost of dissatisfaction with the thermal environment. Similar to TPT, V2B under TOU performs better than energy trading, excelling in both energy usage reduction and thermal comfort management.

\subsection{Decision-Making Process of EV Batteries under Two Tariffs across Six Communities}\label{decision of EV}
\begin{figure*}
\centering
\begin{subfigure}{0.45\textwidth}
    \includegraphics[width=\textwidth]{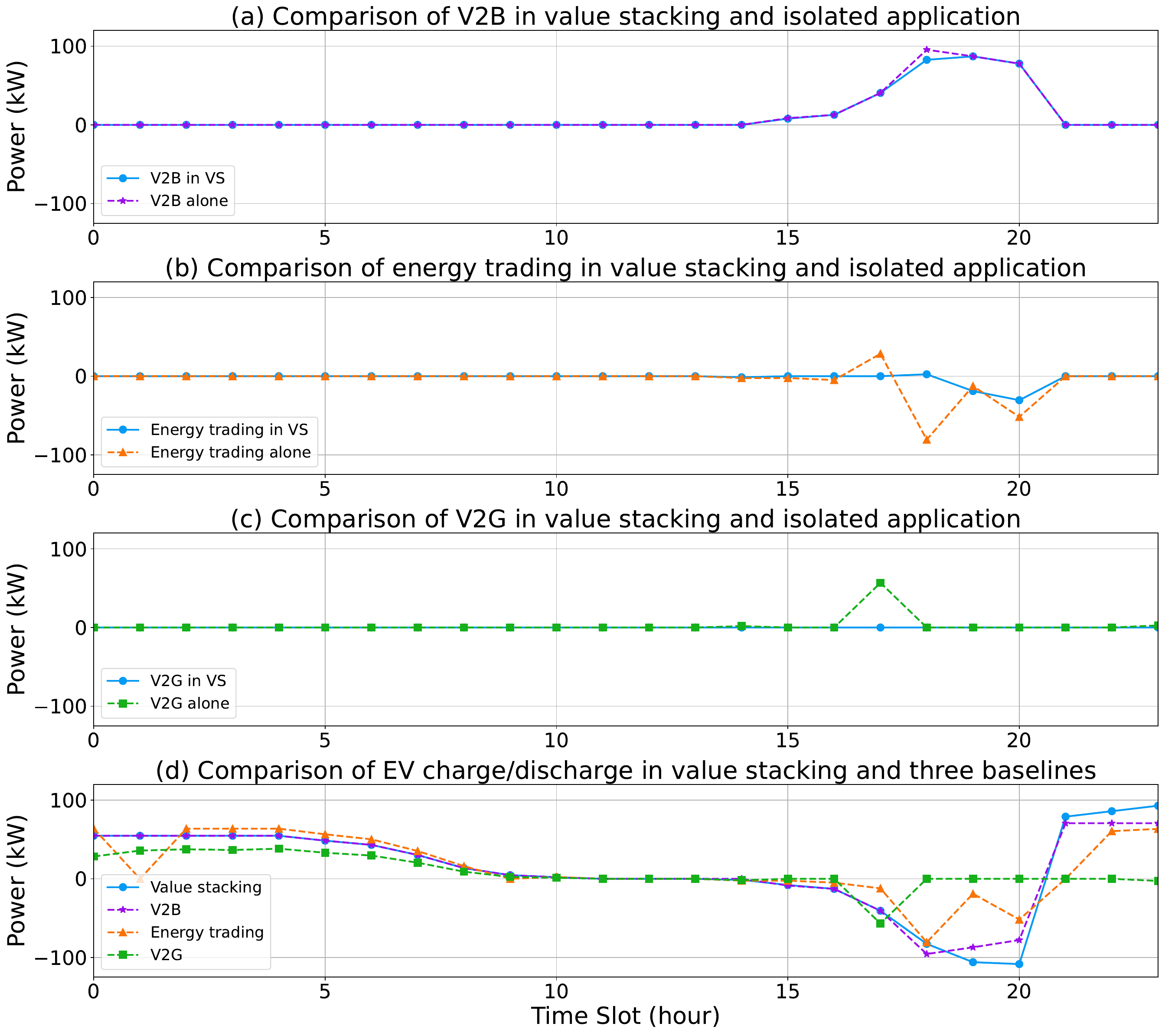}
    \caption{EV batterys' decision of Community 1 under TOU tariff.}
    \label{fig:Community 1 Decision TOU}
\end{subfigure}
\hfill
\begin{subfigure}{0.45\textwidth}
    \includegraphics[width=\textwidth]{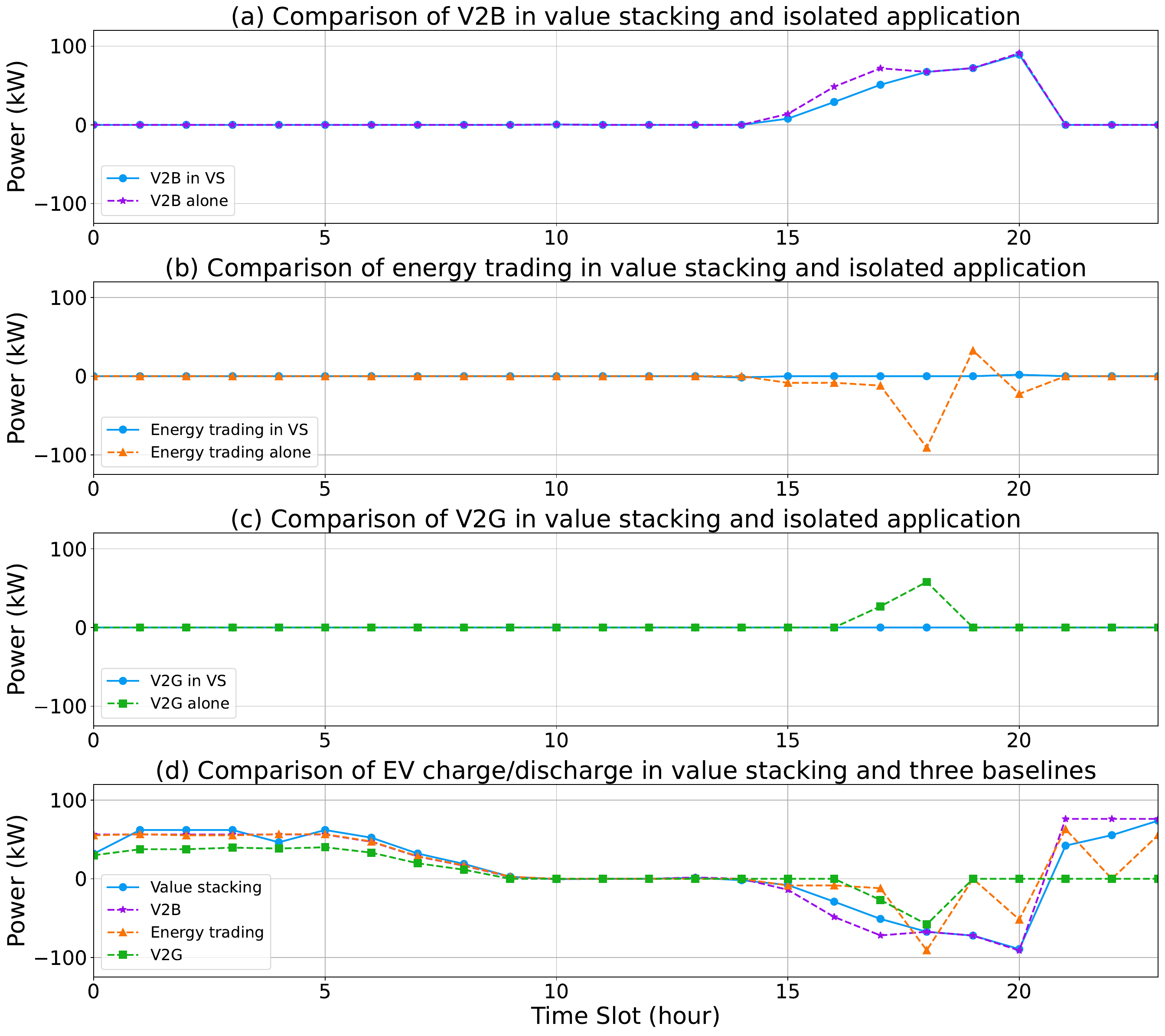}
    \caption{EV batterys' decision of Community 2 under TOU tariff.}
    \label{fig:Community 2 Decision TOU}
\end{subfigure}
\hfill
\begin{subfigure}{0.45\textwidth}
    \includegraphics[width=\textwidth]{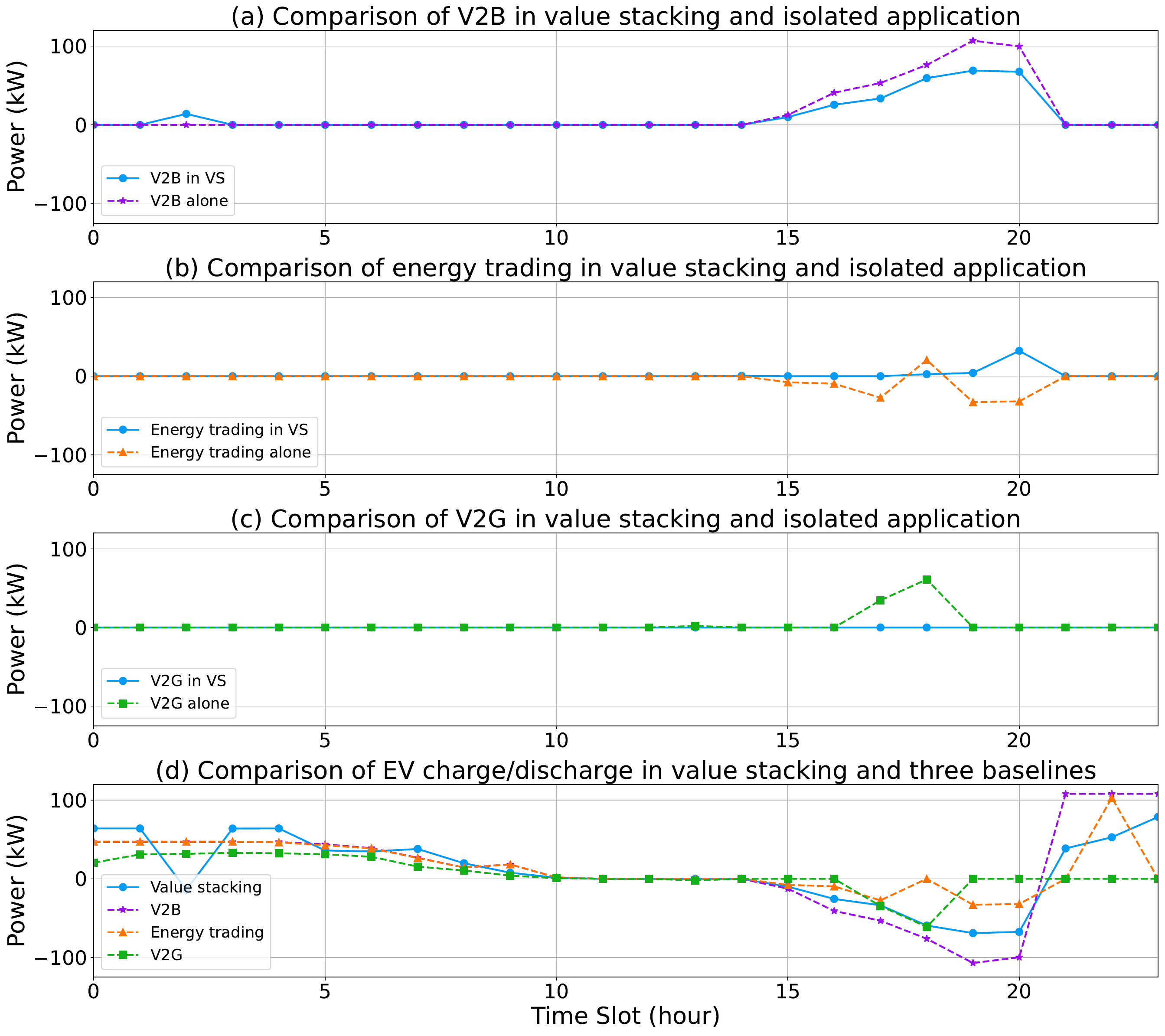}
    \caption{EV batterys' decision of Community 3 under TOU tariff.}
    \label{fig:Community 3 Decision TOU}
\end{subfigure}
\hfill
    \begin{subfigure}{0.45\textwidth}
        \centering
        \includegraphics[width=\textwidth]{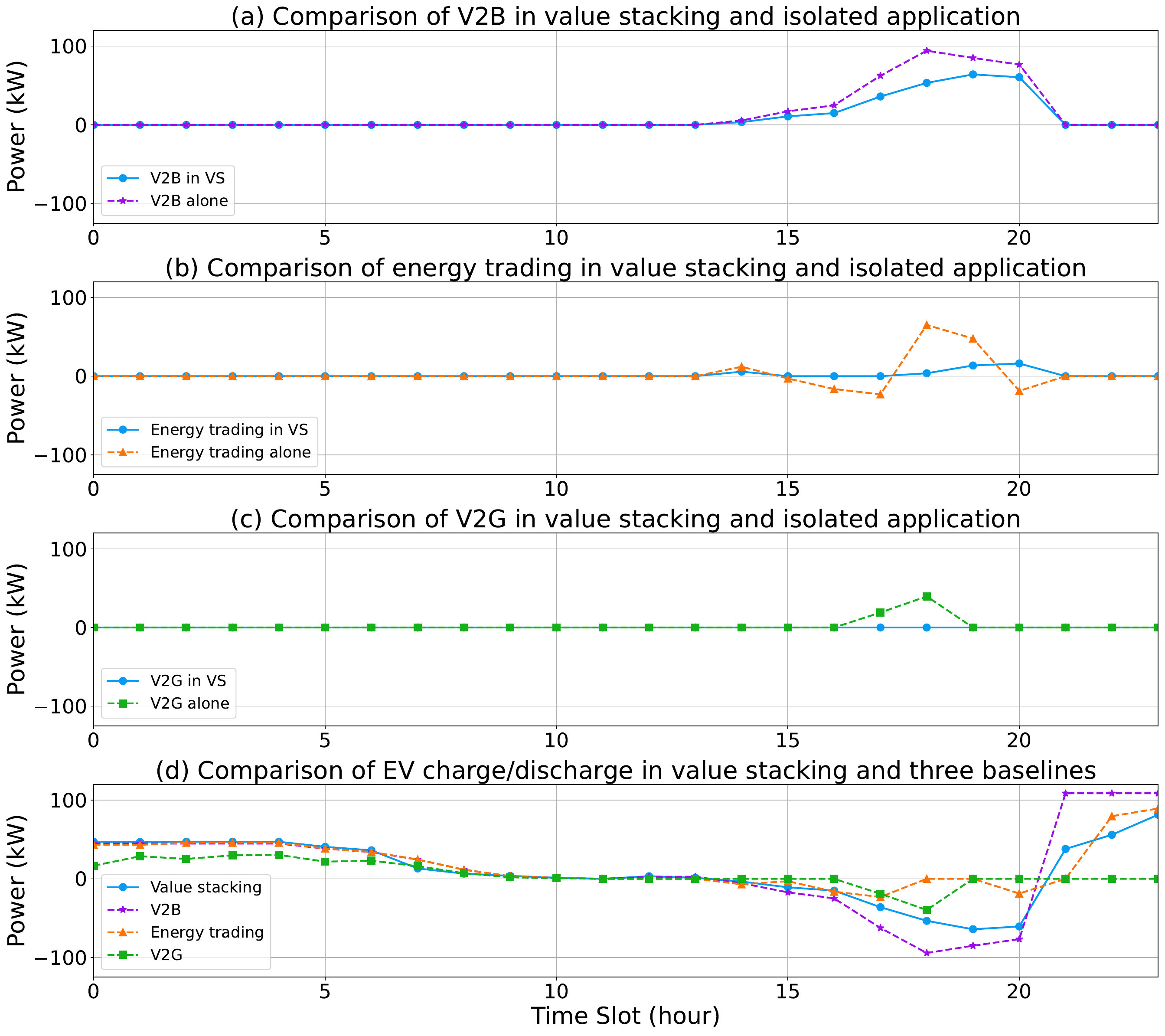}
        \caption{EV batterys' decision of Community 4 under TOU tariff}
        \label{fig:Community 4 Decision TOU}
    \end{subfigure}
\hfill
\begin{subfigure}{0.45\textwidth}
    \includegraphics[width=\textwidth]{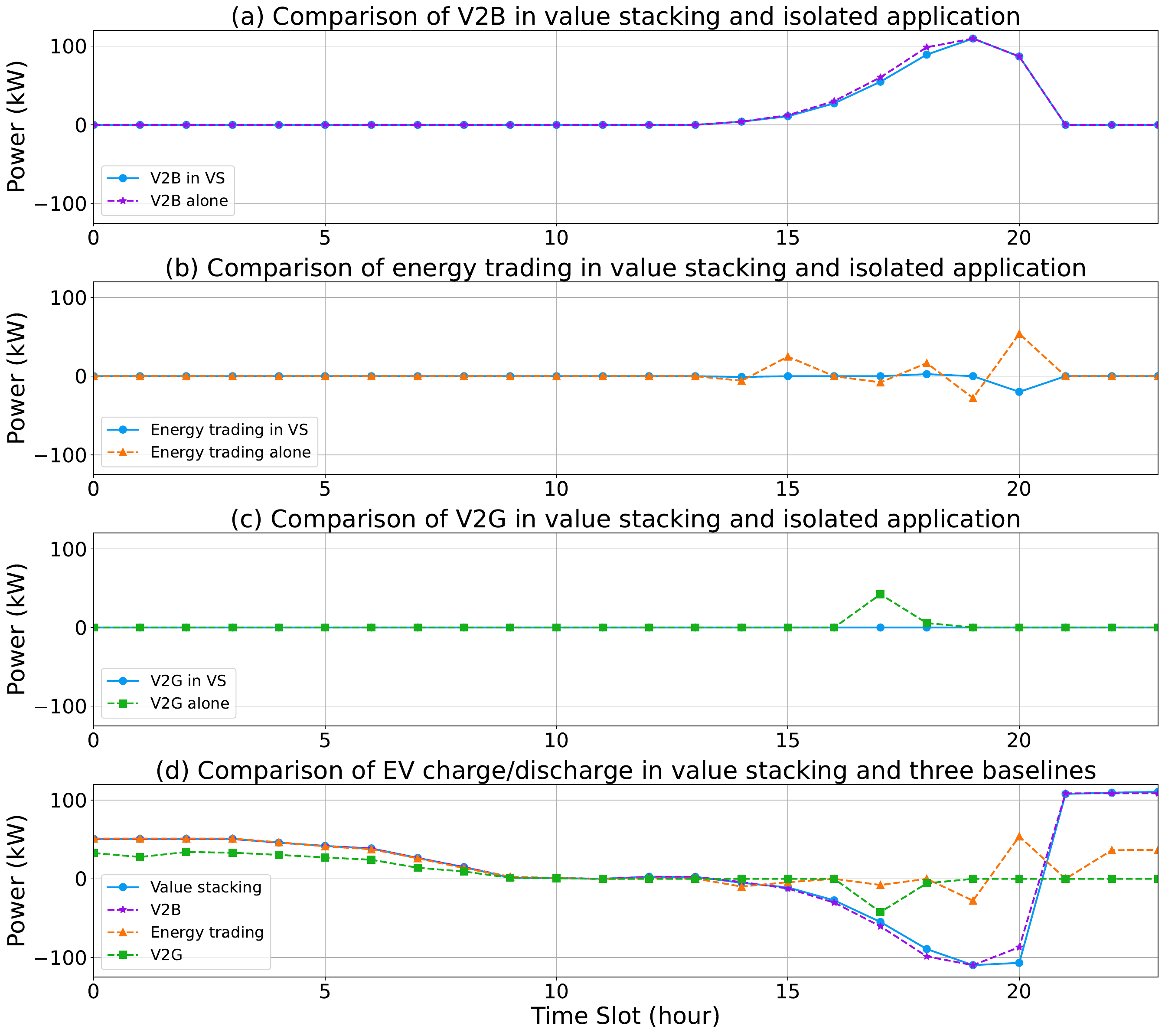}
    \caption{EV batterys' decision of Community 5 under TOU tariff.}
    \label{fig:Community 5 Decision TOU}
\end{subfigure}
\hfill
    \begin{subfigure}{0.45\textwidth}
        \centering
        \includegraphics[width=\textwidth]{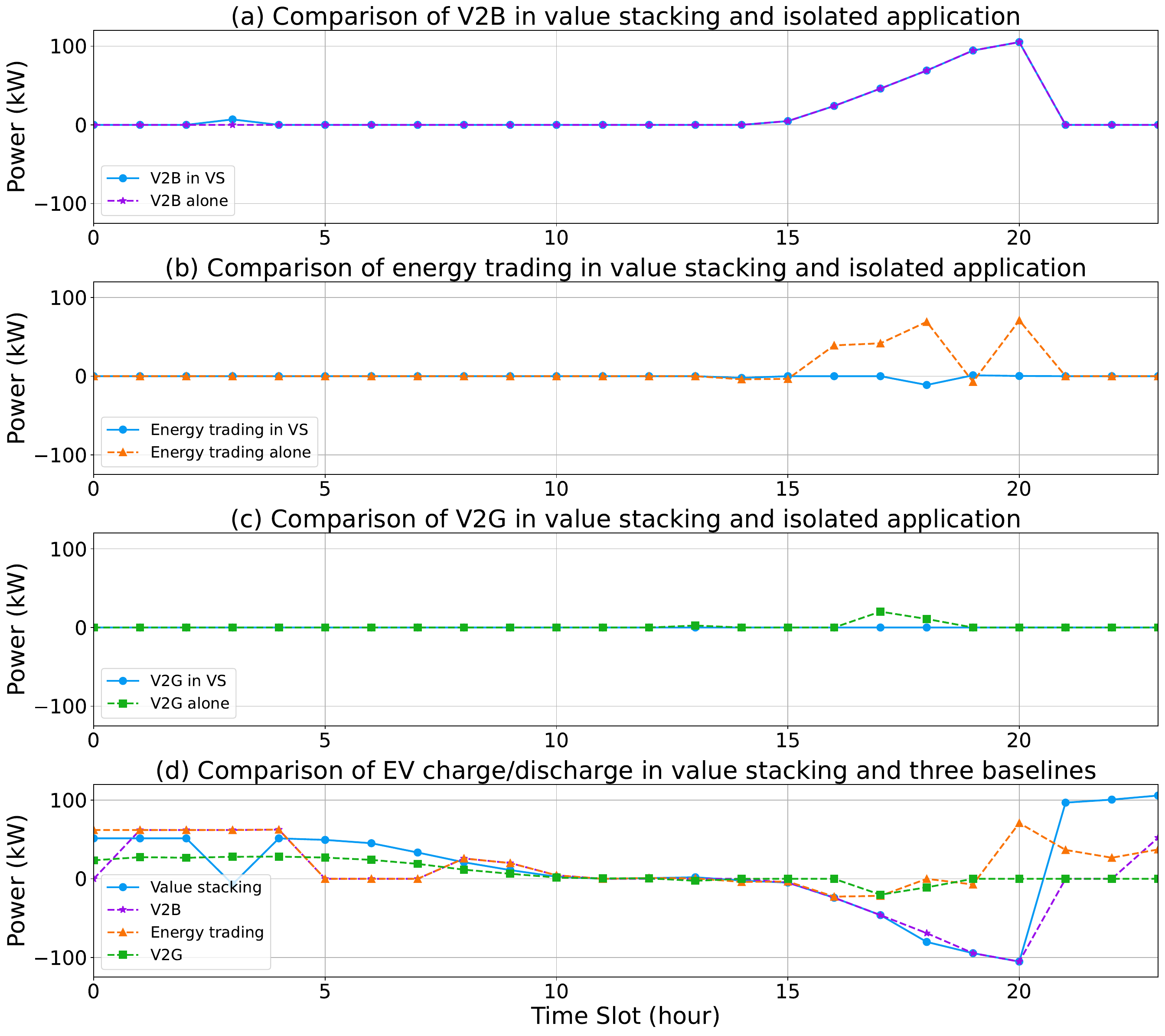}
        \caption{EV batterys' decision of Community 6 under TOU tariff.}
        \label{fig:Community 6 Decision TOU}
    \end{subfigure}
\caption{Decision-making process of EV batteries under TOU across six communities.}
\label{fig: six communities Decision TOU}
\end{figure*}

This work performs a comparative analysis of the decision-making process for EV batteries under two distinct tariffs, TOU and TPT, across six residential communities. Initially, this work examines the decision-making strategies involving V2X value stacking compared to three baselines within each community under the TOU tariff, as illustrated in Fig.~\ref{fig: six communities Decision TOU}. In these figures, solid lines denote the components of V2B, energy trading, and V2G in value stacking, while dashed lines represent these components individually as baselines. Additionally, it is noted that peak hours under the TOU tariff occur between time slots 15 and 21. 

In Community 1, the V2B in value stacking contributes similarly to V2B alone during all time slots besides time slot 18, presented in Fig.~\ref{fig:Community 1 Decision TOU} (a). V2B in value stacking achieves 82.56 kW in time slot 18, which is slightly lower than V2B alone ($95.56$ kW). The comparison of energy trading in value stacking and isolated application is illustrated in Fig.~\ref{fig:Community 1 Decision TOU} (b), where positive and negative values denote power purchases and sales to and from the local energy market respectively. Different from energy trading alone, energy trading in value stacking sells $30.53$ kW of power to local energy market to make profits in time slot 20. Due to the wholesale market price being lower than the peak-hour tariff, encouraging the use of EVs for V2B and energy trading to minimize costs. Therefore, V2G in value stacking only achieves $1.13$ kW in time slot 17, which is much lower than V2G alone ($55.82$ kW), presented in Fig.~\ref{fig:Community 1 Decision TOU} (c). As Fig.~\ref{fig:Community 1 Decision TOU} (d) shows, EV charge and discharge power from value stacking and three baselines are compared, where positive and negative values represent charge and discharge power respectively. From time slots 1 to 10, both value stacking and three baselines will charge to meet the charging demand before departure. From time slot 15, value stacking, V2B, and energy trading perform discharge, where value stacking achieves the highest amount of discharge of $105.96$ kW and $108.33$ kW in time slots 19 and 20, since EVs discharge to perform V2B and sell power to the local energy market. This demonstrates that a single value stream is insufficient to maximize EV battery utility, whereas value stacking optimally uses surplus power to cut overall system costs.

As Fig.~\ref{fig:Community 2 Decision TOU} (a) shows, similar to V2B in value stacking in Community 1, in Community 2, V2B in value stacking achieves $43.56$ kW and $68.32$ kW in time slots 16 and 17, which is slightly lower than V2B alone ($57.22$ kW and $83.31$ kW). Fig.\ref{fig:Community 2 Decision TOU} (b) indicates that, different from energy trading alone, value stacking refrains from trading power during peak hours, thereby balancing individual and total system costs. Similarly to Community 1, V2G in value stacking in Community 2 is limited to $1.32$ kW in time slot 18, which is much lower than V2G alone ($67.32$ kW), as displayed in Fig.\ref{fig:Community 2 Decision TOU} (c). Fig.\ref{fig:Community 2 Decision TOU} (d) shows that from time slots 1 to 10, all EVs charge to satisfy departure demands. From time slot 15 to 20, value stacking, V2B, and energy trading start to perform discharge, where energy trading alone achieves the highest of $97.86$ kW in time slot 18. During times slot 21 and 24, EVs start to charge to the desired energy, where value stacking charges $22.12$ kW power less than V2B alone and energy trading alone in time slot 21. 

Unlike communities 1 and 2, in Community 3, V2B in value stacking contributes 21.02 kW in time slot 3 before EVs departure. From time slot 15 to 20, V2B in value stacking is nearly mirroring V2B alone but lower than it, illustrated in Fig.~\ref{fig:Community 3 Decision TOU} (a). Different from energy trading alone, energy trading in value stacking buys 38.22 kW of power from the local energy market in time slot 20 to reduce the power usage in peak hour, illustrated in Fig.~\ref{fig:Community 3 Decision TOU} (b). As Fig.~\ref{fig:Community 3 Decision TOU} (c) illustrates, V2G in value stacking achieves 0, due to the unfavorable market pricing relative to peak tariffs. In Fig.~\ref{fig:Community 3 Decision TOU} (d), both value stacking and three baselines charge to meet the charging demand, where value stacking discharges to perform V2B in time slot 3, and then charges to the desired energy before departure. From time slot 15 to 20, value stacking, V2B, and energy trading start to perform discharge, where V2B alone achieves the highest of 101.26 kW in time slot 19. Different from V2B alone, value stacking purchases power from the local energy market in time slot 20, which reduce the discharge power from time 15 to 20 and charge power from time 21 to 24. 

In Community 4, V2B in value stacking is nearly mirroring V2B alone but lower than it from time slot 16 to 20, depicted in Fig.~\ref{fig:Community 4 Decision TOU} (a). Similar to energy trading in value stacking in Community 3, it purchases 19.31 kW of power from the local energy market in time slot 20, which is shown in Fig.~\ref{fig:Community 4 Decision TOU} (b). Similar to communities 1, 2, and 3, V2G in value stacking only achieves 0.95 kW in time slot 18, which is much lower than V2G alone (52.21 kW), presented in Fig.~\ref{fig:Community 4 Decision TOU} (c). For the comparison of EV charge and discharge power from value stacking and three baselines, both value stacking and three baselines charge to meet the charging demand from time slot 1 to 10 in Fig.~\ref{fig:Community 4 Decision TOU} (d). From time slot 15 to 20, value stacking, V2B, and energy trading perform discharge, where V2B alone achieves the highest of 101.26 kW in time slot 18. Unlike V2B alone, value stacking purchases power from the local energy market in time slot 20, which reduces its charge power from time 21 to 24. 

In Community 5, the V2B in value stacking contributes similarly to V2B alone during all time slots, presented in Fig.~\ref{fig:Community 5 Decision TOU} (a). The comparison of energy trading in value stacking and isolated application is illustrated in Fig.~\ref{fig:Community 5 Decision TOU} (b), where energy trading in value stacking sells 20.23 kW of power to the local energy market to make profits in time slot 20. As Fig.~\ref{fig:Community 5 Decision TOU} (c) illustrates, there is no V2G in value stacking, due to the wholesale market price being lower than the peak-hour tariff. Fig.~\ref{fig:Community 5 Decision TOU} (d) shows the comparison of EV charge and discharge power from value stacking and three baselines, where both value stacking and three baselines charge to meet the charging demand before departure. From time slot 15 to 20, value stacking, V2B, and energy trading perform discharge, where value stacking achieves the highest amount of discharge of 101.22 kW in time slot 20 since value stacking discharges to perform V2B and sell power to the local energy market. This suggests that value stacking can contribute greater value to our system.  

In Community 6, the V2B in value stacking contributes similarly to V2B alone during all time slots, presented in Fig.~\ref{fig:Community 6 Decision TOU} (a). The comparison of energy trading in value stacking and isolated application is illustrated in Fig.~\ref{fig:Community 6 Decision TOU} (b), where energy trading alone purchases power from the local energy market in time slots 16, 17, 18, and 20, respectively. Different from energy trading alone, energy trading in value stacking sells 9.37 kW of power to the local energy market in time slot 18 to reduce the power usage of the entire system, illustrated in Fig.~\ref{fig:Community 6 Decision TOU} (b). As Fig.~\ref{fig:Community 6 Decision TOU} (c) illustrates, value stacking does not perform V2G across all time slots. Fig.~\ref{fig:Community 6 Decision TOU} (d) presents, from time slot 1 to 10, both value stacking and three baselines charge to meet the charging demand before departure, where value stacking discharges to perform V2B in time slot 3, and then charge to the desired energy before departure. From time slot 15 to 20, value stacking, V2B, and energy trading start to perform discharge, where value stacking achieves the highest of 100.25 kW in time slot 20. This indicates that a single value stream cannot maximize the value of EV batteries, while value stacking leverages surplus power to reduce overall system costs.

\begin{figure*}
\centering
\begin{subfigure}{0.45\textwidth}
    \includegraphics[width=\textwidth]{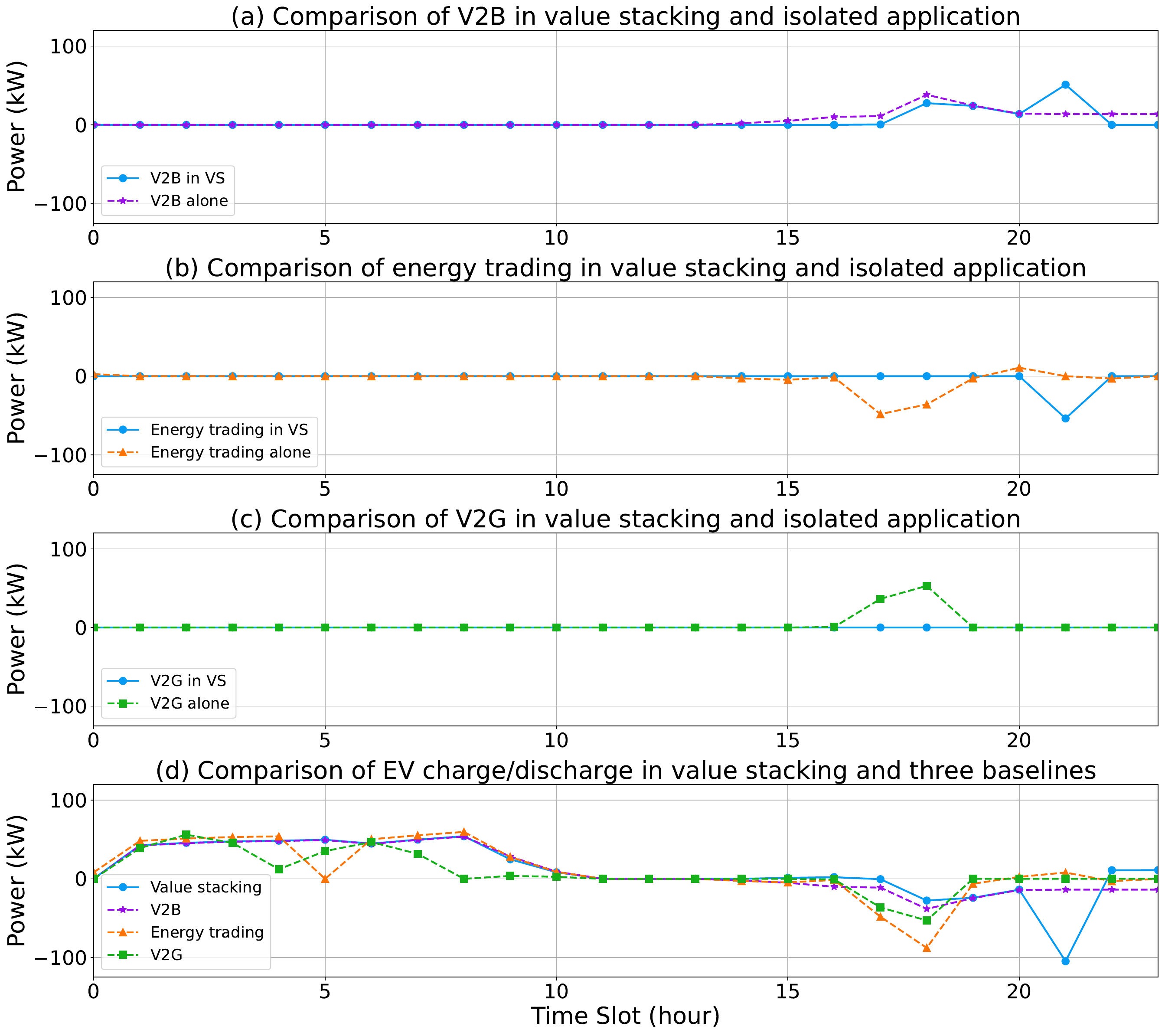}
    \caption{EV batterys' decision of Community 1 under TPT tariff.}
    \label{fig:Community 1 Decision TPT}
\end{subfigure}
\hfill
\begin{subfigure}{0.45\textwidth}
    \includegraphics[width=\textwidth]{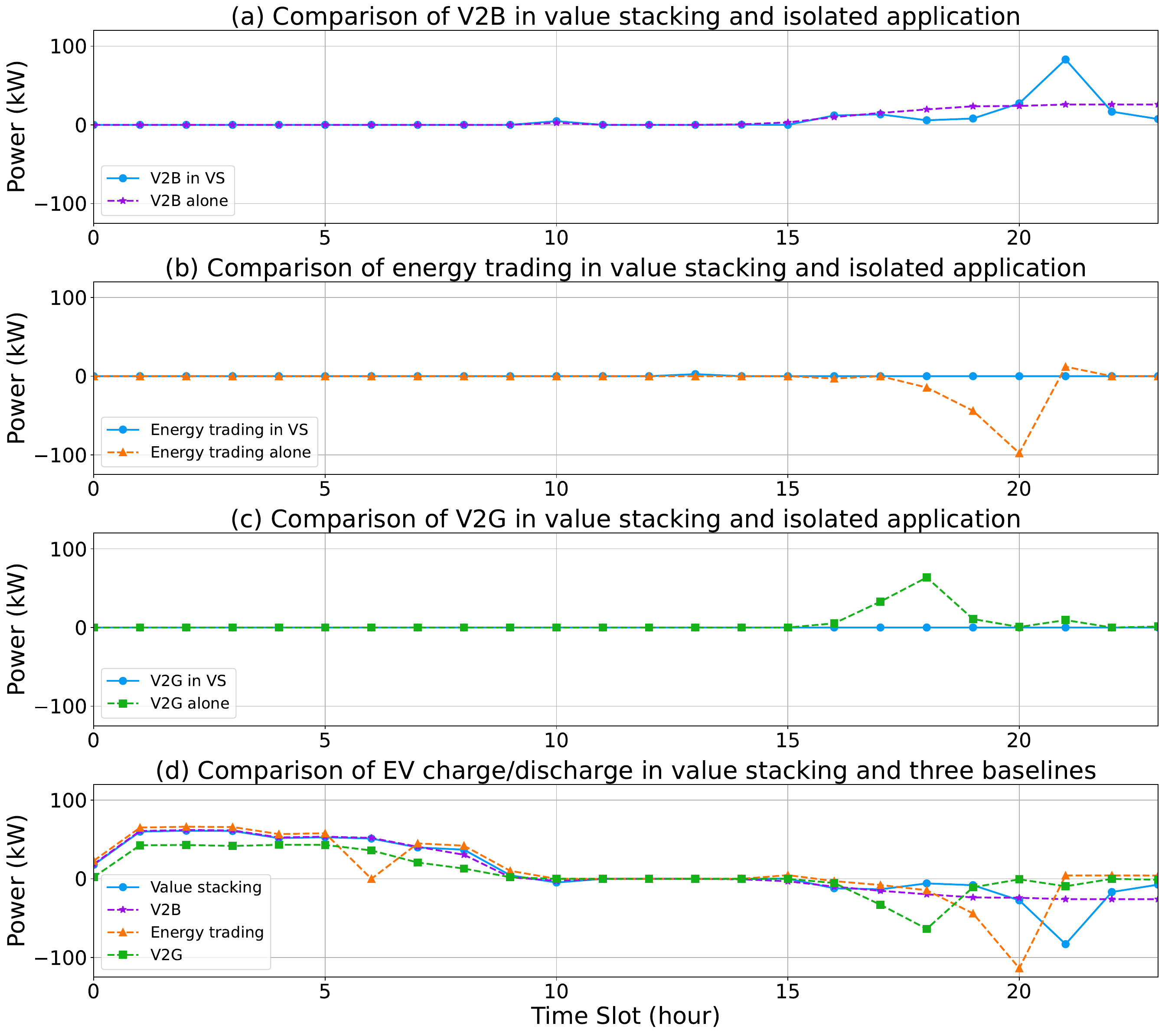}
    \caption{EV batterys' decision of Community 2 under TPT tariff.}
    \label{fig:Community 2 Decision TPT}
\end{subfigure}
\hfill
\begin{subfigure}{0.45\textwidth}
    \includegraphics[width=\textwidth]{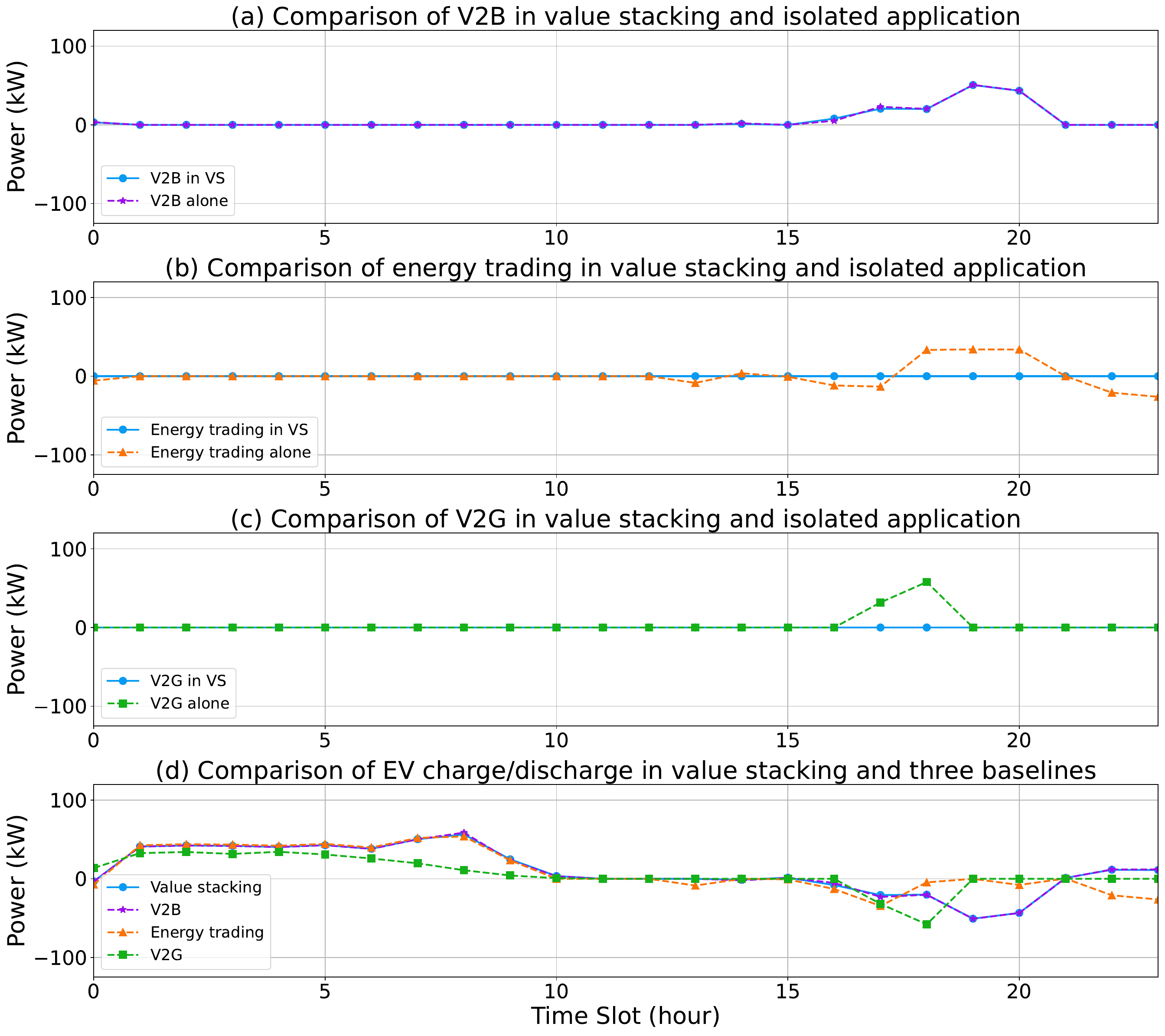}
    \caption{EV batterys' decision of Community 3 under TPT tariff.}
    \label{fig:Community 3 Decision TPT}
\end{subfigure}
\hfill
    \begin{subfigure}{0.45\textwidth}
        \centering
        \includegraphics[width=\textwidth]{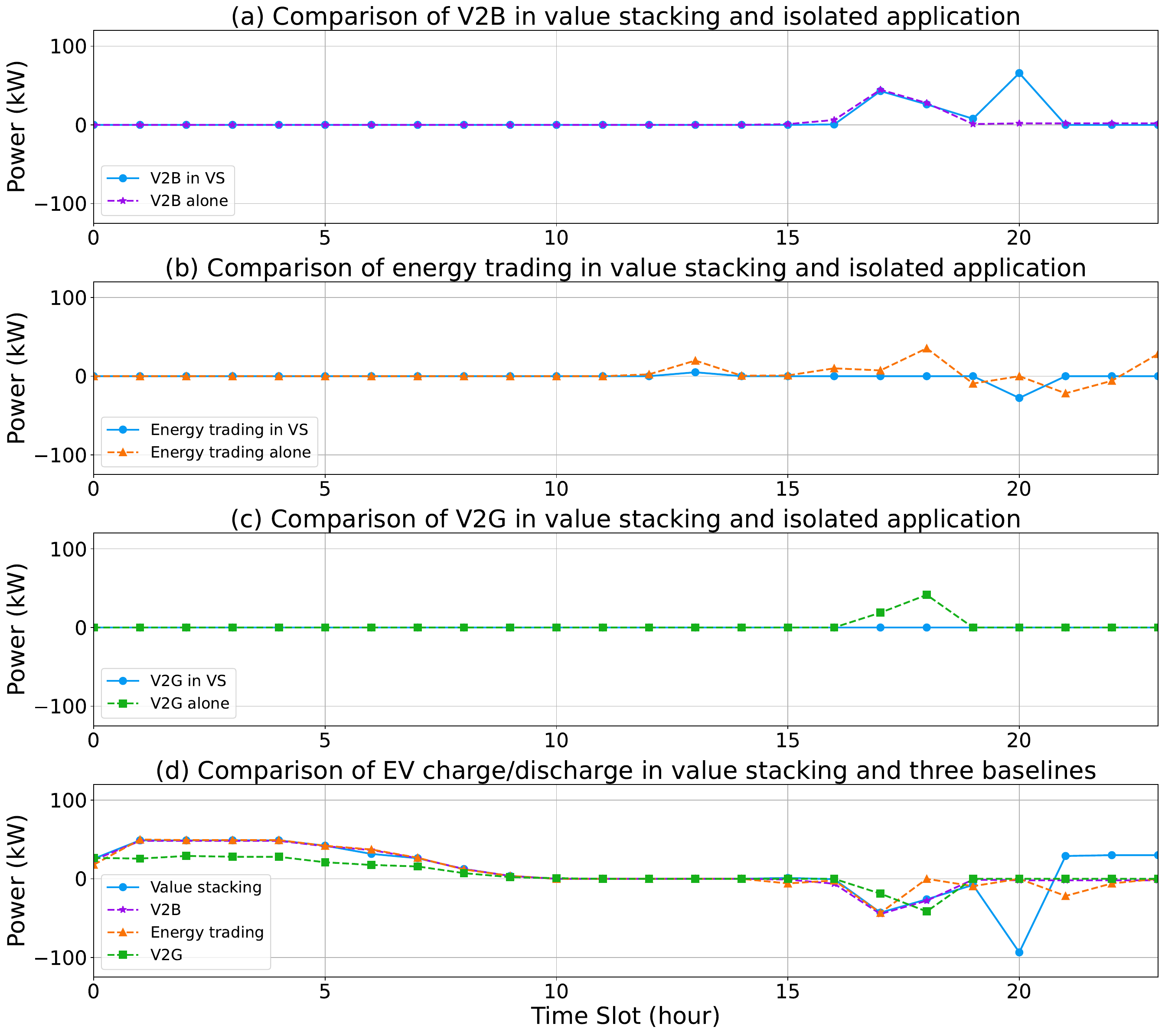}
        \caption{EV batterys' decision of Community 4 under TPT tariff}
        \label{fig:Community 4 Decision TPT}
    \end{subfigure}
\hfill
\begin{subfigure}{0.45\textwidth}
    \includegraphics[width=\textwidth]{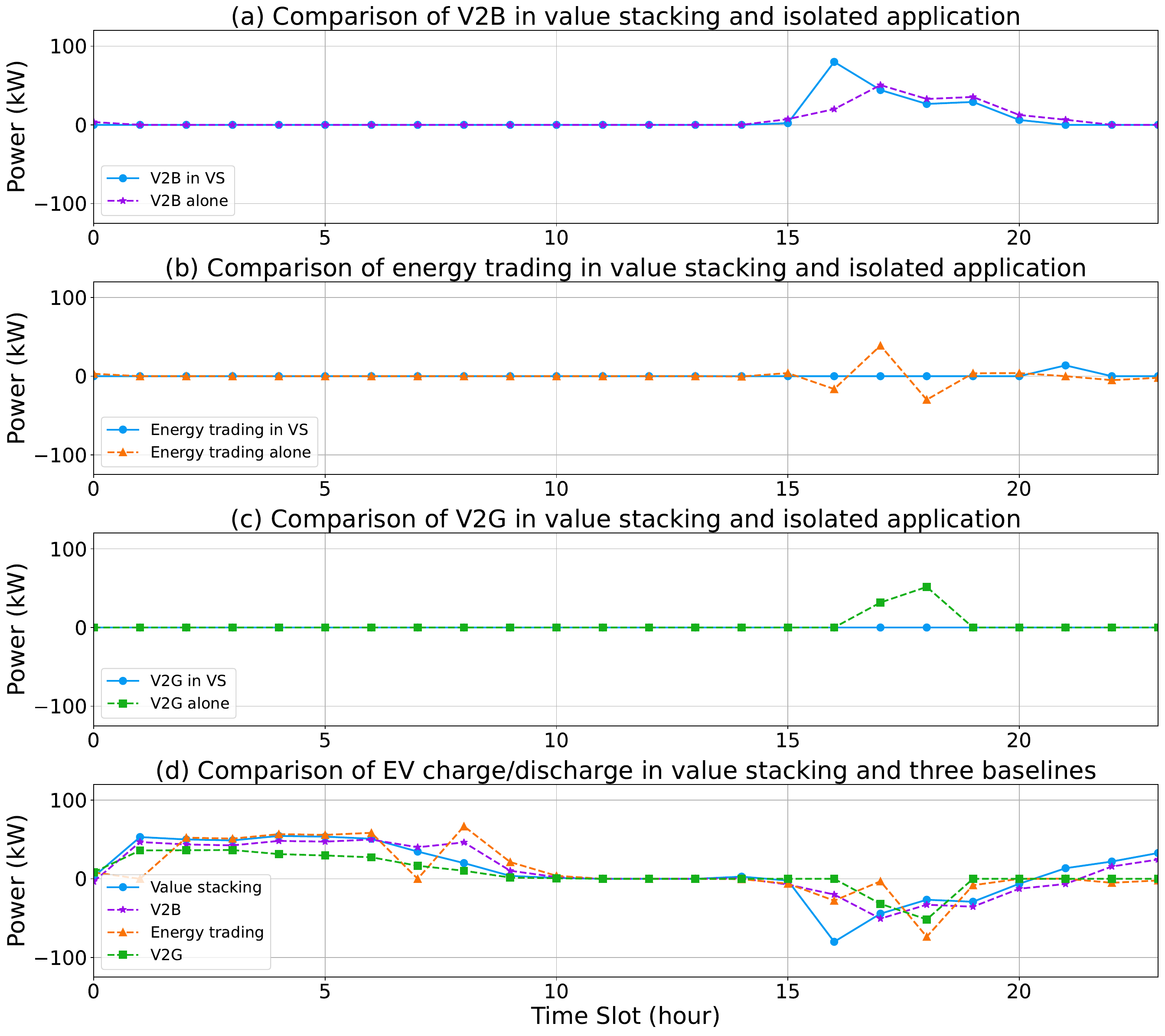}
    \caption{EV batterys' decision of Community 5 under TPT tariff.}
    \label{fig:Community 5 Decision TPT}
\end{subfigure}
\hfill
    \begin{subfigure}{0.45\textwidth}
        \centering
        \includegraphics[width=\textwidth]{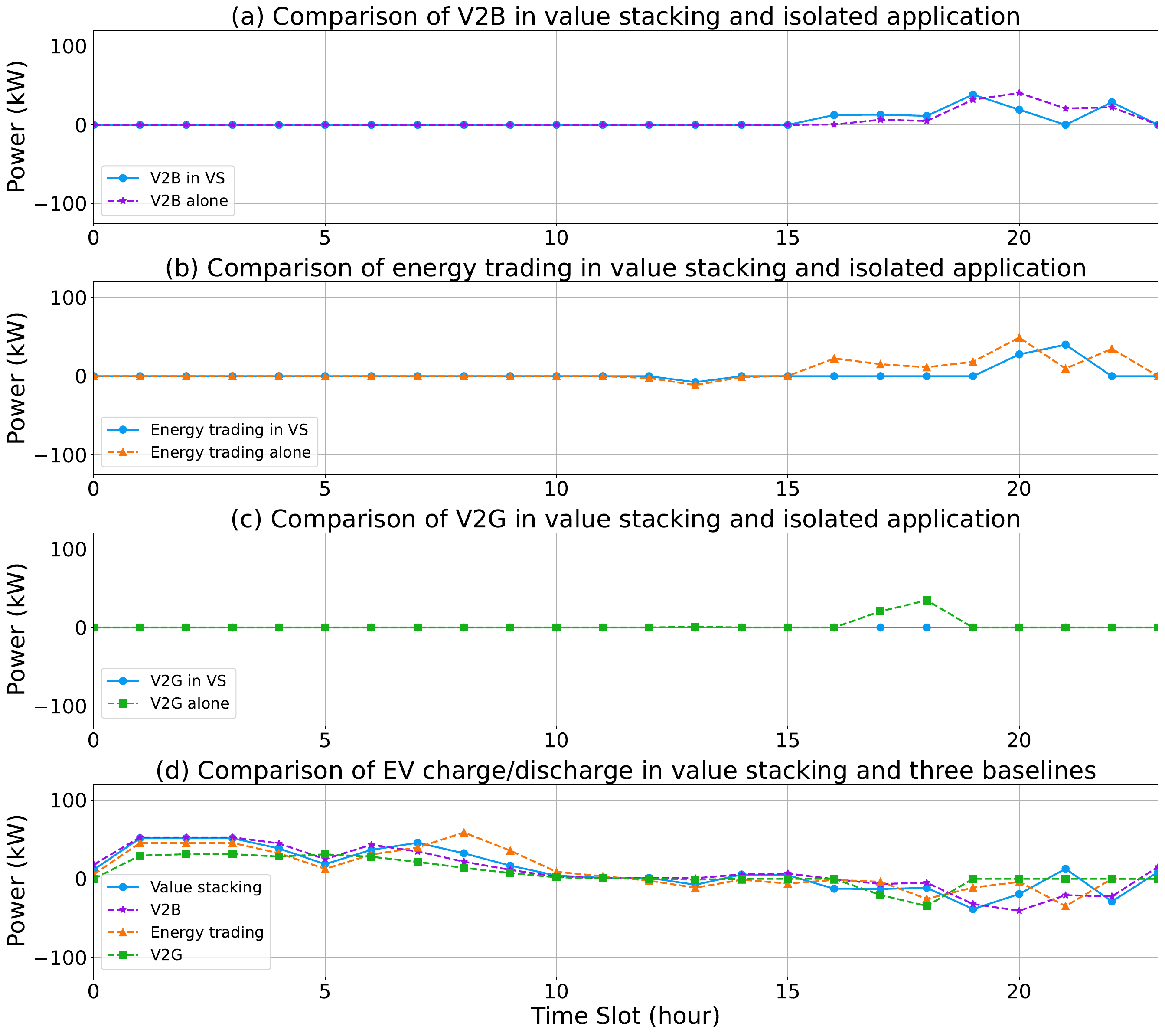}
        \caption{EV batterys' decision of Community 6 under TPT tariff.}
        \label{fig:Community 6 Decision TPT}
    \end{subfigure}
\caption{Decision-making process of EV batteries under TPT across six communities.}
\label{fig: six communities Decision TPT}
\end{figure*}

In the second part of our analysis, this work explores the decision-making process for EV batteries in six communities under the TPT tariff, as depicted in Fig.~\ref{fig: six communities Decision TPT}. Contrary to the TOU tariff, TPT calculates the total cost, consisting of energy price plus peak price. 

In Community 1, presented in Fig.~\ref{fig:Community 1 Decision TPT} (a), the V2B in value stacking contributes a total of 139.16 kW, which is less than V2B alone (146.77 kW) during all time slots. However, in time slot 21, V2B in value stacking achieves 39.25 kW, which is higher than that of V2B alone. As Fig.~\ref{fig:Community 1 Decision TPT} (b) illustrates, energy trading in value stacking sells 38.93 kW of power to the local energy market to make profits in time slot 21. Due to the wholesale market price being lower than the peak-hour tariff, encouraging the use of EVs for V2B and energy trading to minimize costs. Therefore, V2G in value stacking only achieves 1.13 kW in time slot 18, which is much lower than V2G alone (51.32 kW), as presented in Fig.~\ref{fig:Community 1 Decision TPT} (c). In Fig.~\ref{fig:Community 1 Decision TPT} (d), EV charge and discharge power from value stacking and three baselines are compared, where both value stacking and three baselines charge to meet the charging demand before departure. From time slots 15 to 21, value stacking, V2B, and energy trading perform discharge, where value stacking achieves the highest amount of discharge of 102.32 kW and in time slot 21. In time slot 21, EVs discharge to perform V2B and sell power to the local energy market. This suggests that a single value stream cannot maximize the value of EV batteries, whereas value stacking leverages surplus power to diminish entire system costs. 

Similar to V2B in value stacking in Community 1, in Community 2, the V2B in value stacking contributes a total of 137.34 kW, which is less than V2B alone (175.40 kW) during all time slots, which is depicted in Fig.~\ref{fig:Community 2 Decision TPT} (a). In time slot 21, V2B in value stacking achieves 79.15 kW, which is higher than V2B alone. As Fig.~\ref{fig:Community 2 Decision TPT} (b) illustrates, there is no energy trading in value stacking across all time slots. In Fig.~\ref{fig:Community 2 Decision TPT} (c), there is no V2G in value stacking across all time slots, while V2G alone sells 99.32 kW to the main grid to make profits. For the comparison of EV charge and discharge power from value stacking and three baselines, both value stacking and three baselines charge to meet the charging demand before departure, as depicted in Fig.~\ref{fig:Community 2 Decision TPT} (d). From time slot 15 to 20, value stacking, V2B, and energy trading start to perform discharge, where energy trading alone achieves the highest of 117.12 kW in time slot 20. During times slots 23 and 24, EVs start to charge to the desired energy, where value stacking charges a total of 32.11 kW power, which is less than V2B alone and energy trading alone. It means the entire energy system trades off the total cost and individual cost and trades energy among other communities to fulfill total cost reduction. 

Unlike communities 1 and 2, in Community 3, the V2B in value stacking contributes similarly to V2B alone during all time slots, illustrated in Fig.~\ref{fig:Community 3 Decision TPT} (a). for the comparison of energy trading, there is no energy trading in value stacking across all time slots, which shows in Fig.~\ref{fig:Community 3 Decision TPT} (b). Similar to Community 1, V2G in value stacking in community 3 only achieves 0.91 kW in time slot 19, which is much lower than V2G alone (62.70 kW), presented in Fig.~\ref{fig:Community 3 Decision TOU} (c). For the comparison of EV charge and discharge power from value stacking and three baselines, both value stacking and three baselines charge to meet the charging demand before departure depicted in Fig.~\ref{fig:Community 3 Decision TPT} (d). From time slot 1 to 10, both value stacking and three baselines charge to meet the charging demand before departure. From time slot 15 to 24, value stacking and V2B discharge the same power across all times. Different from V2B alone, value stacking purchases power from the local energy market in time slot 20, which reduces the discharge power from time 15 to 20 and charge power from time 21 to 24. 

In Community 4, V2B in value stacking is nearly mirroring V2B alone but is lower than it from time slot 15 to 19, as depicted in Fig.~\ref{fig:Community 4 Decision TPT} (a). Similar to energy trading in value stacking in Community 3, it purchases 19.31 kW of power from the local energy market in time slot 20, as shown in Fig.~\ref{fig:Community 4 Decision TPT} (b). Similar to communities 1, 2, and 3, V2G in value stacking only achieves 0.95 kW in time slot 18, which is much lower than V2G alone (52.21 kW), presented in Fig.~\ref{fig:Community 4 Decision TPT} (c). For the comparison of EV charge and discharge power from value stacking and three baselines, both value stacking and three baselines charge to meet the charging demand before departure depicted in Fig.~\ref{fig:Community 4 Decision TPT} (d). From time slot 1 to 10, both value stacking and three baselines charge to meet the charging demand before departure. From time slot 15 to 20, value stacking, V2B, and energy trading perform discharge, where V2B alone achieves a maximum of 98.76 kW in time slot 20. This indicates that a single value stream cannot maximize the value of EV batteries, while value stacking leverages surplus power to reduce overall system costs.   

In Community 5, the V2B in value stacking contributes a total of 168.92 kW, which is slightly lower than V2B alone (169.54 kW) during all time slots, presented in Fig.~\ref{fig:Community 5 Decision TPT} (a). The comparison of energy trading in value stacking and isolated application is illustrated in Fig.~\ref{fig:Community 5 Decision TPT} (b), where energy trading in value stacking purchases 10.13 kW of power from the local energy market in time slot 21. As Fig.~\ref{fig:Community 5 Decision TPT} (c) illustrates, value stacking does not discharge to perform V2G, due to the wholesale market price being lower than the tariff. Fig.~\ref{fig:Community 5 Decision TPT} (d) illustrates the comparison of EV charge and discharge power from value stacking and three baselines, where both value stacking and three baselines charge to meet the charging demand from time slots 1 to 10. From time slot 15 to 20, value stacking, V2B, and energy trading discharge power, where value stacking achieves the highest amount of discharge of 81.02 kW in time slot 16 since value stacking discharges to perform V2B. It suggests that V2B in value stacking contributes most to three value streams in Community 5. 

In Community 6, presented in Fig.~\ref{fig:Community 6 Decision TPT} (a), the V2B in value stacking achieves the highest power (57.11 kW) in time slot 19, while V2B alone reaches the same power in time slot 20. The reason is that in Fig.~\ref{fig:Community 6 Decision TPT} (b) energy trading in value stacking purchases power from the local energy market in time slots 20 and 21 to reduce the electricity cost, which allows V2B in value stacking to contribute early than V2B alone. As Fig.~\ref{fig:Community 6 Decision TPT} (c) illustrates, V2G in value stacking only contributes 1.03 kW in time slot 19, which is much lower than V2G alone. Fig.~\ref{fig:Community 6 Decision TPT} (d) presents, from time slot 1 to 10, both value stacking and three baselines charge to meet the charging demand before departure. From time slot 15 to 20, value stacking, V2B, and energy trading perform discharge, where V2B achieves the highest of 29.25 kW in time slot 19. Different from V2B and energy trading alone, value stacking purchases power from the local energy market in time slots 19 and 20, which reduces the discharge power from time 15 to 20.

\subsection{Comparison of HVAC Usage with and without V2B}\label{HVAC with and without V2B}
\begin{figure}[!t]
    \centering
    \includegraphics[width=0.8\linewidth]{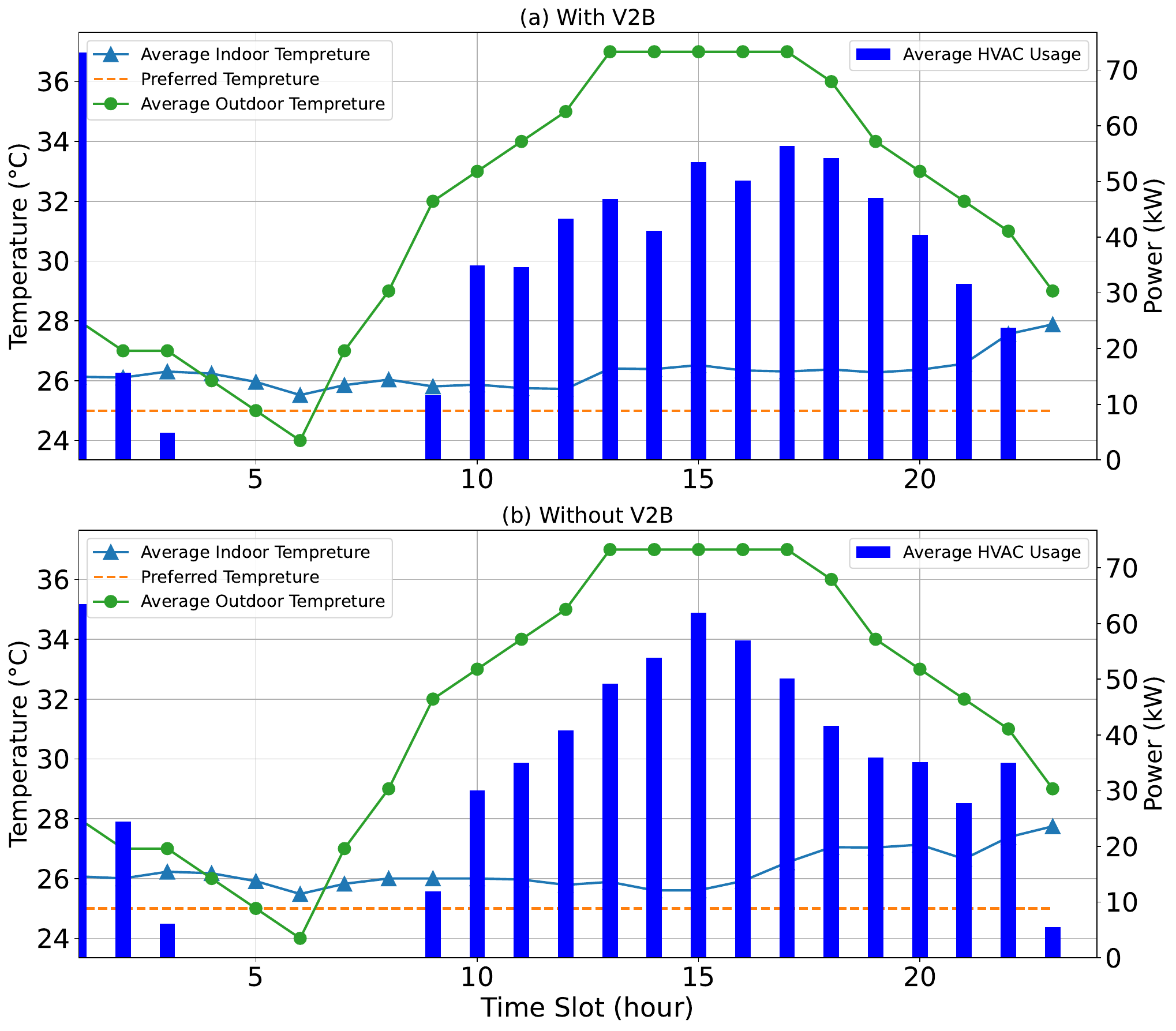}
    \caption{Comparison of HVAC usage with and without V2B under TPT.}
\label{fig:Comparison of HVAC under TPT}
\end{figure}

To evaluate the contribution of V2B in HVAC usage, this work compares the HVAC usage in scenarios with and without V2B in the proposed value-stacking optimization problem under TPT, as depicted in Fig.~\ref{fig:Comparison of HVAC under TPT}, where the left and right Y-axes represent indoor temperature ($^{\circ}$C) and HVAC usage (kW). The bar chart showcases the average HVAC usage (kW) across six residential communities, with the dashed line indicating the preferred indoor temperature, the blue solid line representing the indoor temperature, and the green line denoting the outdoor temperature, respectively. Notably, there is a peak load of HVAC usage in time slot 15 in Fig.~\ref{fig:Comparison of HVAC under TPT} (b), and V2B can help in peak load shaving and make the HVAC usage more evenly, as depicted in Fig.~\ref{fig:Comparison of HVAC under TPT} (a). Although the difference between indoor temperature and preferred indoor temperature in Fig.~\ref{fig:Comparison of HVAC under TPT} (a) is larger than that in Fig.~\ref{fig:Comparison of HVAC under TPT} (b), HVAC usage in time slots 20 and 21 in Fig.~\ref{fig:Comparison of HVAC under TPT} (a) is lower than that in Fig.~\ref{fig:Comparison of HVAC under TPT} (b). This variance could be attributed to the fact that V2B supplies most power to inflexible building load usage, sacrificing the thermal comfort to trade off between the cost of electricity price and dissatisfaction with the thermal environment. Therefore, V2B can assist in peak load shaving in HVAC usage, resulting in a more evenly distributed HVAC usage under TPT.

\begin{figure}[!t]
    \centering
    \includegraphics[width=0.8\linewidth]{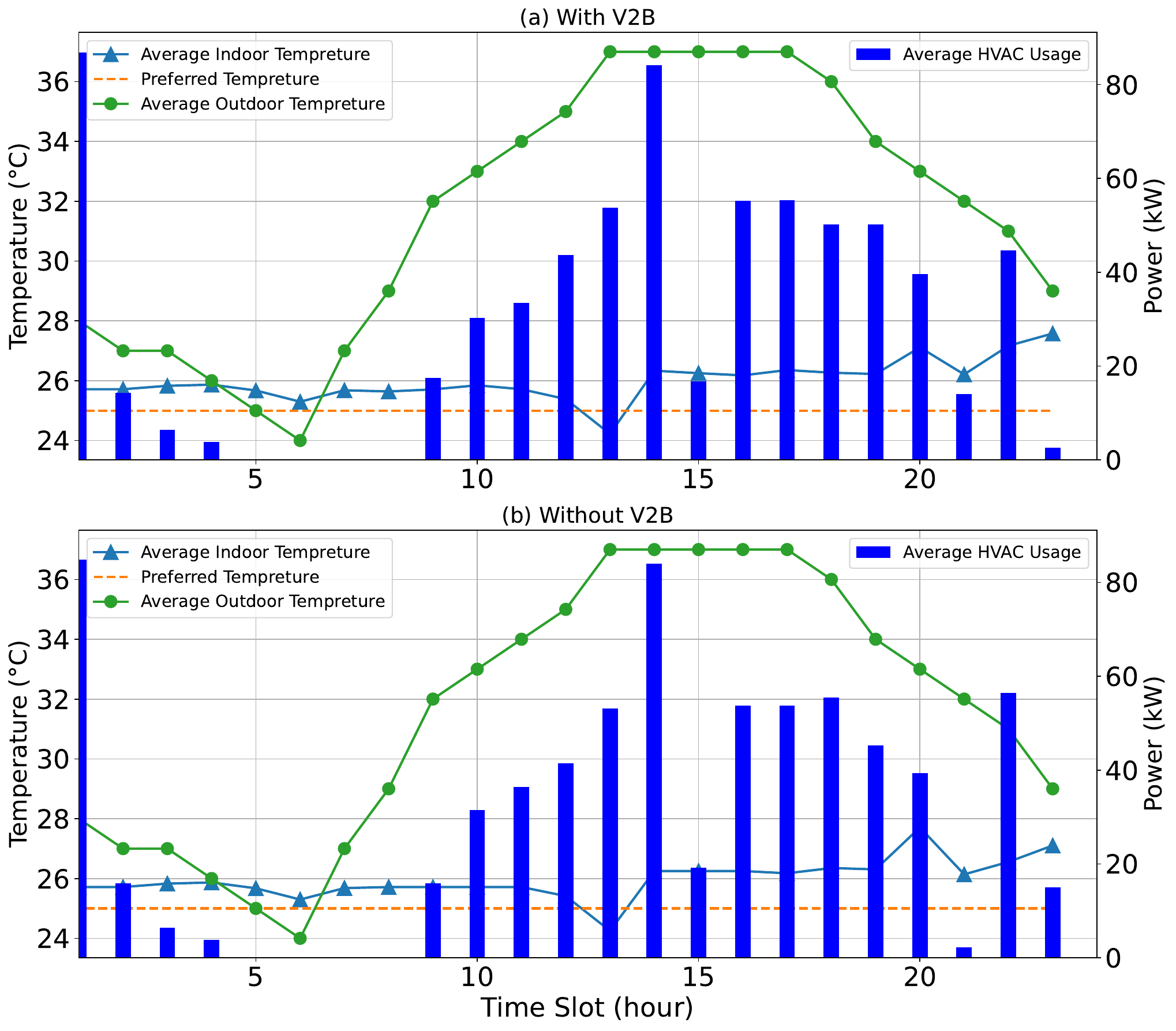}
    \caption{Comparison of HVAC usage with and without V2B under TOU.}
\label{fig:Comparison of HVAC under TOU}
\end{figure}

Fig.~\ref{fig:Comparison of HVAC under TOU} shows HVAC usage with and without V2B integration within the proposed value-stacking optimization problem considering TOU pricing. In contrast to the TPT tariff, the inclusion of V2B in TOU does not facilitate peak load shaving for HVAC usage. However, during peak load hours (from time slots 15 to 21), V2B exhibits the potential to curtail HVAC usage in time slots 16, 17, and 19, as shown in Fig.\ref{fig:Comparison of HVAC under TOU} (a). Notably, while HVAC usage in time slots 19 and 20 in Fig.\ref{fig:Comparison of HVAC under TOU} (a) exceeds that of Fig.~\ref{fig:Comparison of HVAC under TOU} (b), the introduction of V2B reduces dissatisfaction with the thermal environment. This is plausibly attributed to V2B's ability to accommodate all inflexible load usage during time slots 19 and 20, thereby enabling the residential communities to maintain thermal comfort without sacrificing it to reduce overall costs. Consequently, V2B integration under TOU demonstrates the potential to decrease HVAC usage during peak hours.

\subsection{Evaluation of Forecast Accuracy} \label{forecast assessment}
To predict building demand, PV generation, and EV arrivals for residential communities, the GRU-EN-TFD model is employed as described in Section \ref{Transformer}. This model forecasts hourly building load, PV generation, and EV arrivals across six communities over a period of 200 days. Specifically, this work uses data from September 1, 2019, to September 1, 2020, as the training dataset, data from September 2, 2020, to December 31, 2020, as the validation dataset, and data from January 1, 2021, to July 19, 2021, as the testing dataset.

\subsubsection{Model Training Process}
\begin{figure}[!t]
\centering
\includegraphics[width=1.0\linewidth]{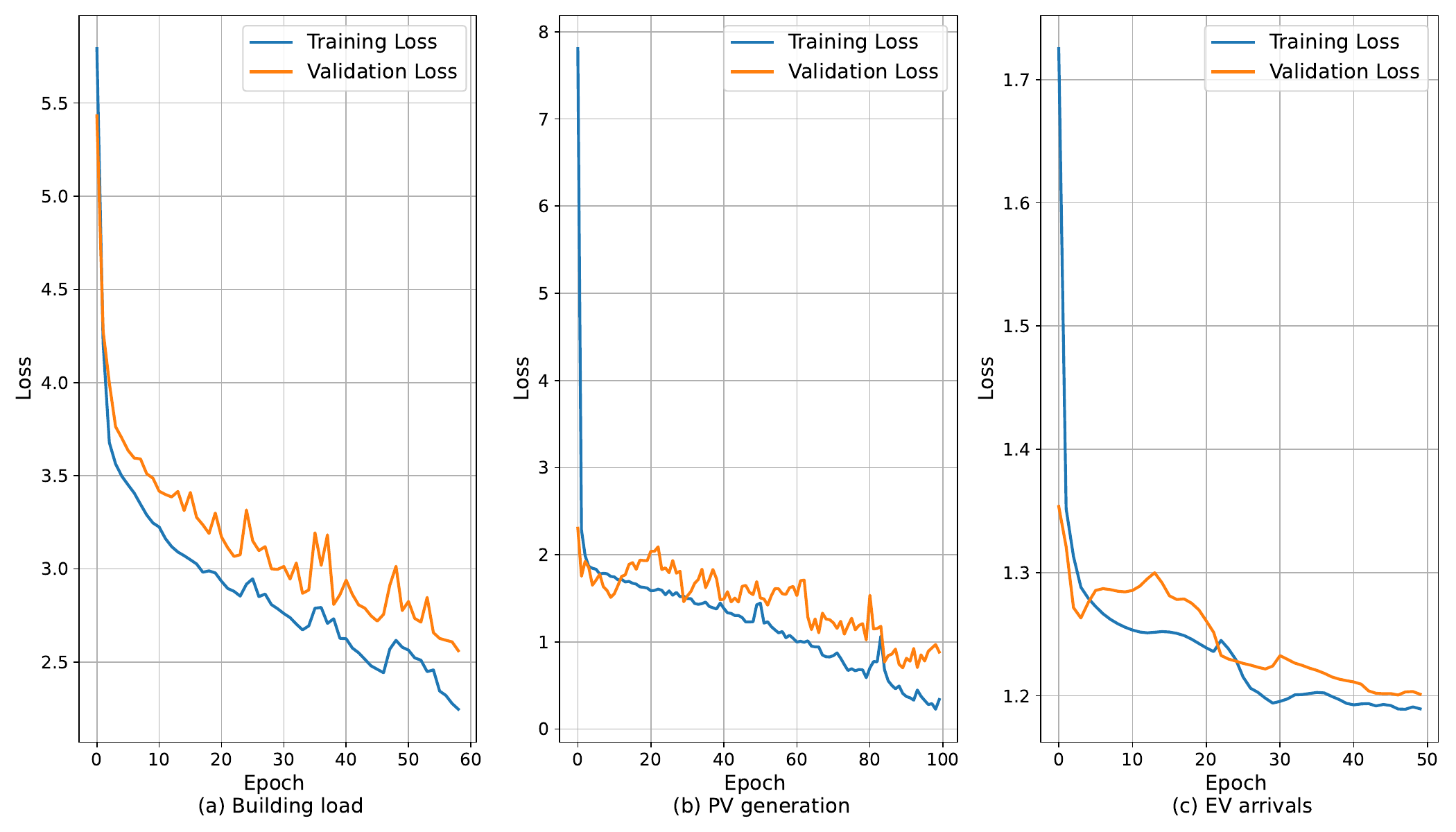} 
\caption{Training process of the GRU-EN-TFD model for forecasting building demand, PV generation, and EV arrivals, respectively.}
\label{fig: Training process}
\end{figure}

Fig.~\ref{fig: Training process} illustrates the training and validation losses of the proposed GRU-EN-TFD model, with each subplot corresponding to forecasts for building load, PV generation, and EV arrivals, respectively. As shown in Fig.~\ref{fig: Training process} (a), both training and validation losses start relatively high and decrease rapidly within the first 10 epochs, indicating quick learning for building load data. The training loss shows a gradual but steady decrease. The validation loss decreases alongside the training loss but exhibits some variability, suggesting occasional fluctuations in model performance on validation data.
Fig.~\ref{fig: Training process} (b) presents a similar trend in training and validation of PV generation forecasting. The losses are lower than those in the building load subplot in Fig.~\ref{fig: Training process} (a), potentially suggesting simpler patterns or more effective learning for PV generation forecast.
Fig.~\ref{fig: Training process} (c) illustrates the trend in losses during training and validation for EV arrival forecasting. Both the training and validation losses level off more quickly (around epoch 10) and exhibit less variability and converge more closely, indicating a robust model performance for EV arrivals forecasting. 

\subsubsection{Comparison of Model Performance}
This work further analyzes the RMSE of the GRU-EN-TFD model in the testing data set across six communities, in comparison with the LSTM model. As shown in \ref{Tab: 1}, the GRU-EN-TFD model consistently outperforms the LSTM model in the testing data set across the six residential communities. For example, the GRU-EN-TFD model demonstrates a higher predictive accuracy in Community 3 in forecasting EV arrivals at EVPL with an RMSE of 1.16. This performance underscores the model's ability in handling the variability in EV parking patterns in Community 3. Additionally, the model excels in predicting PV generation in Community 2, achieving an RMSE of 0.85, compared to an RMSE of 1.09 by the LSTM. 

Overall, the GRU-EN-TFD model consistently demonstrates strong performance across diverse forecasting tasks, including building load demand, PV generation, and EV arrivals, in all communities. Its notable accuracy, especially in the complex scenarios of EV arrival and PV generation forecasting, highlights its ability to effectively capture temporal variabilities, thus improving decision-making in residential communities.

\begin{table*}
\centering
\caption{Comparative RMSE performance metrics for GRU-EN-TFD and LSTM models on test data sets across six communities.}
\label{Tab: 1}
\begin{tblr}{
  cell{1}{1} = {r=2}{},
  cell{1}{2} = {c=3}{c},
  cell{1}{5} = {c=3}{c},
  cell{1}{8} = {c=3}{c},
  cell{5}{1} = {r=2}{},
  cell{5}{2} = {c=3}{c},
  cell{5}{5} = {c=3}{c},
  cell{5}{8} = {c=3}{c},
  vlines,
  hline{1,3-5,7-9} = {-}{},
  hline{2,6} = {2-10}{},
}
Model      & Community 1~ &      &        & Community 2 &      &        & Community 3 &      &        \\
           & Building 1   & PV 1 & EVPL 1 & Building 2  & PV 2 & EVPL 2 & Building 3  & PV 3 & EVPL 3 \\
LSTM       & 3.31         & 1.15 & 1.33   & 3.56        & 1.09 & 1.36   & 3.42        & 1.13 & 1.47   \\
GRU-EN-TFD & 2.67         & 0.86 & 1.21   & 2.73        & 0.85 & 1.25   & 2.56        & 0.93 & 1.16   \\
Model      & Community 4  &      &        & Community 5 &      &        & Community 6 &      &        \\
           & Building 4   & PV 4 & EVPL 4 & Building 5  & PV 5 & EVPL 5 & Building 6  & PV 6 & EVPL 6 \\
LSTM       & 3.03         & 1.05 & 1.22   & 3.23        & 1.11 & 1.43   & 3.38        & 1.21 & 1.32   \\
GRU-EN-TFD & 2.83         & 0.91 & 1.19   & 2.71        & 0.89 & 1.28   & 2.66        & 0.92 & 1.17   
\end{tblr}
\end{table*}

\begin{figure}[!t]
\centering
\includegraphics[width=1.0\linewidth]{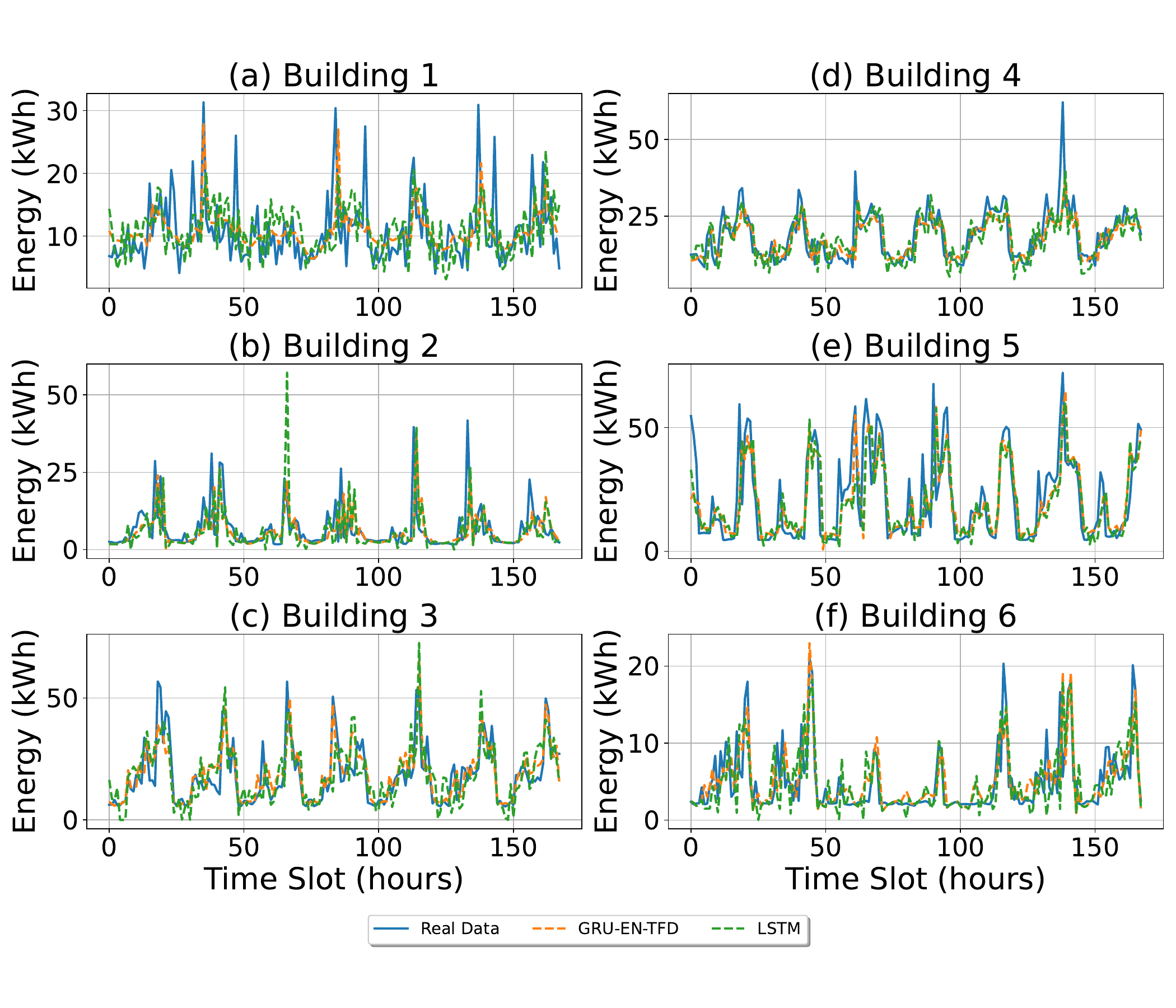} 
\caption{Comparison of the forecasted load and actual load for six communities.}
\label{fig:load forecast}
\end{figure}

\begin{figure}[!t]
\centering
\includegraphics[width=1.0\linewidth]{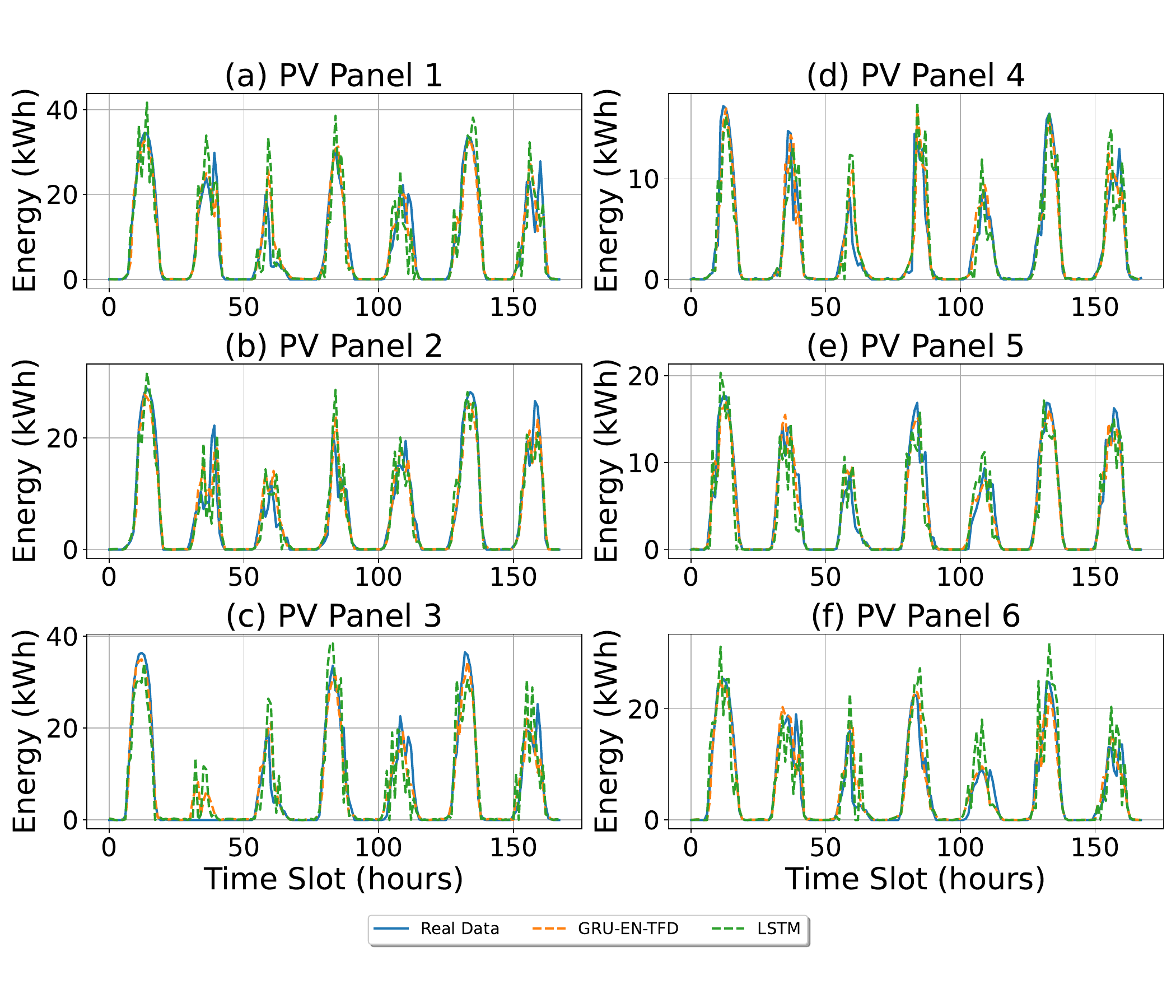} 
\caption{Comparison of the forecasted PV generation and actual PV generation for six communities.}
\label{fig:PV forecast}
\end{figure}

\begin{figure}[!t]
\centering
\includegraphics[width=1.0\linewidth]{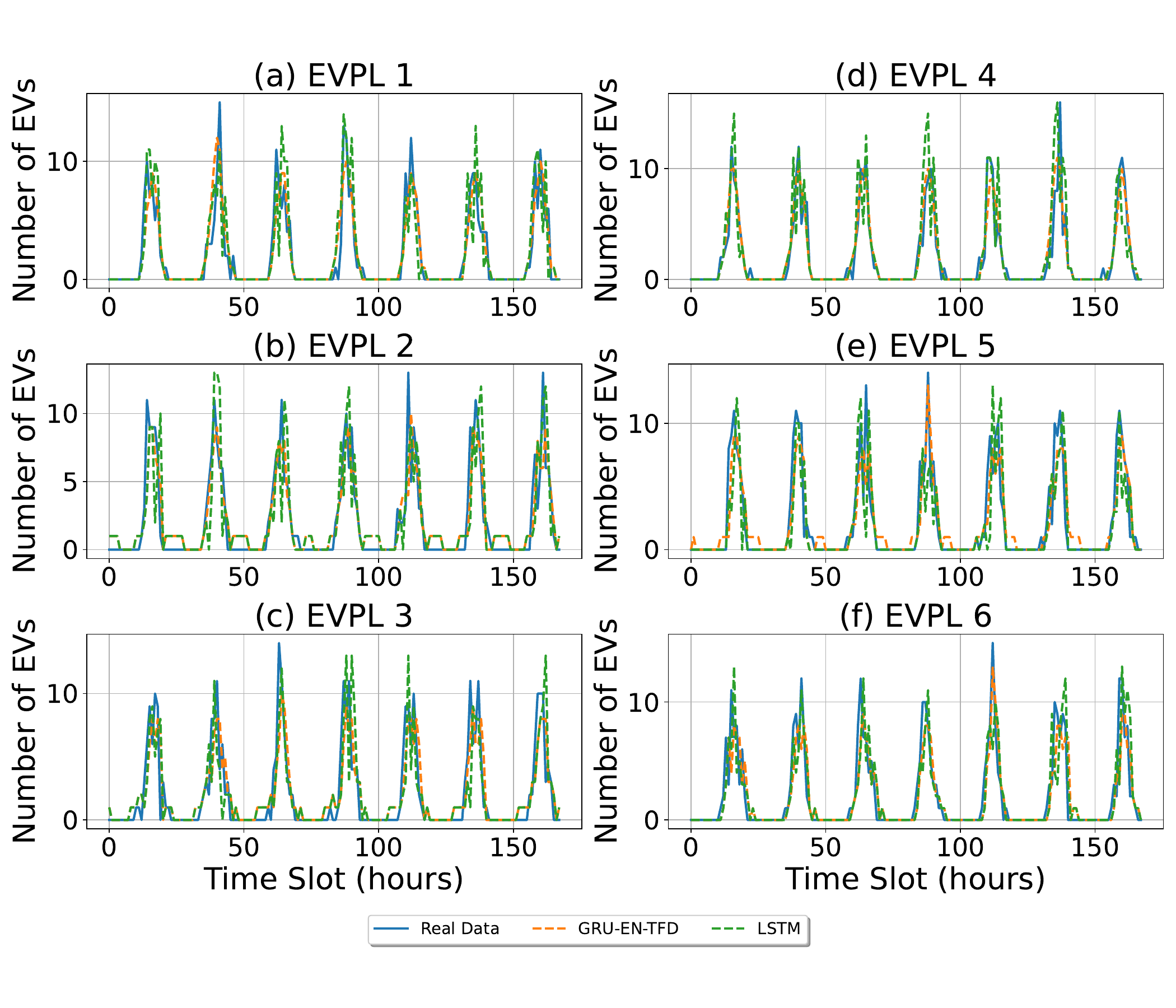} 
\caption{Comparison of the forecasted EV arrivals and actual arrivals for six communities.}
\label{fig:EV forecast}
\end{figure}

Additionally, this work compares the forecasted building demand generated by the GRU-EN-TFD model and an LSTM model to the actual load data across six communities over one week (168 hours). The comparison results are illustrated in Fig.~\ref{fig:load forecast}, demonstrating that the GRU-EN-TFD model consistently outperforms the LSTM model, particularly during peak hours.

Similarly, this work compares predictions for PV generation across these communities using both models, as shown in Fig.~\ref{fig:PV forecast}. The GRU-EN-TFD model again demonstrates higher accuracy than the LSTM model, particularly during time slots, such as 83 and 135 for PV panel 1, time slot 83 for PV panel 2, time slot 84 for PV panel 3, and time slot 133 for PV panel 6. Compared to the building load forecast, PV generation predictions exhibit better performance, suggesting that modeling PV generation patterns is relatively simpler than capturing residential electricity consumption behaviors. 

Additionally, Fig.~\ref{fig:EV forecast} illustrates the forecasts for hourly EV arrivals in the same six communities using both models. As with the building load and PV generation forecasts, the GRU-EN-TFD model exhibits superior performance over the LSTM model, also in particular during peak arrival times.

In summary, the GRU-EN-TFD model outperforms the LSTM model across all forecasting tasks. The accuracy of EV arrivals and PV generation forecasts surpasses that of building load predictions, indicating GRU-EN-TFD's robust performance across different forecasting tasks.

\subsection{Influence of Forecasting Errors on Value-Stacking Performance} \label{impact of prediction error on VS}
\begin{figure}[!t]
    \centering
    \includegraphics[width=0.8\linewidth]{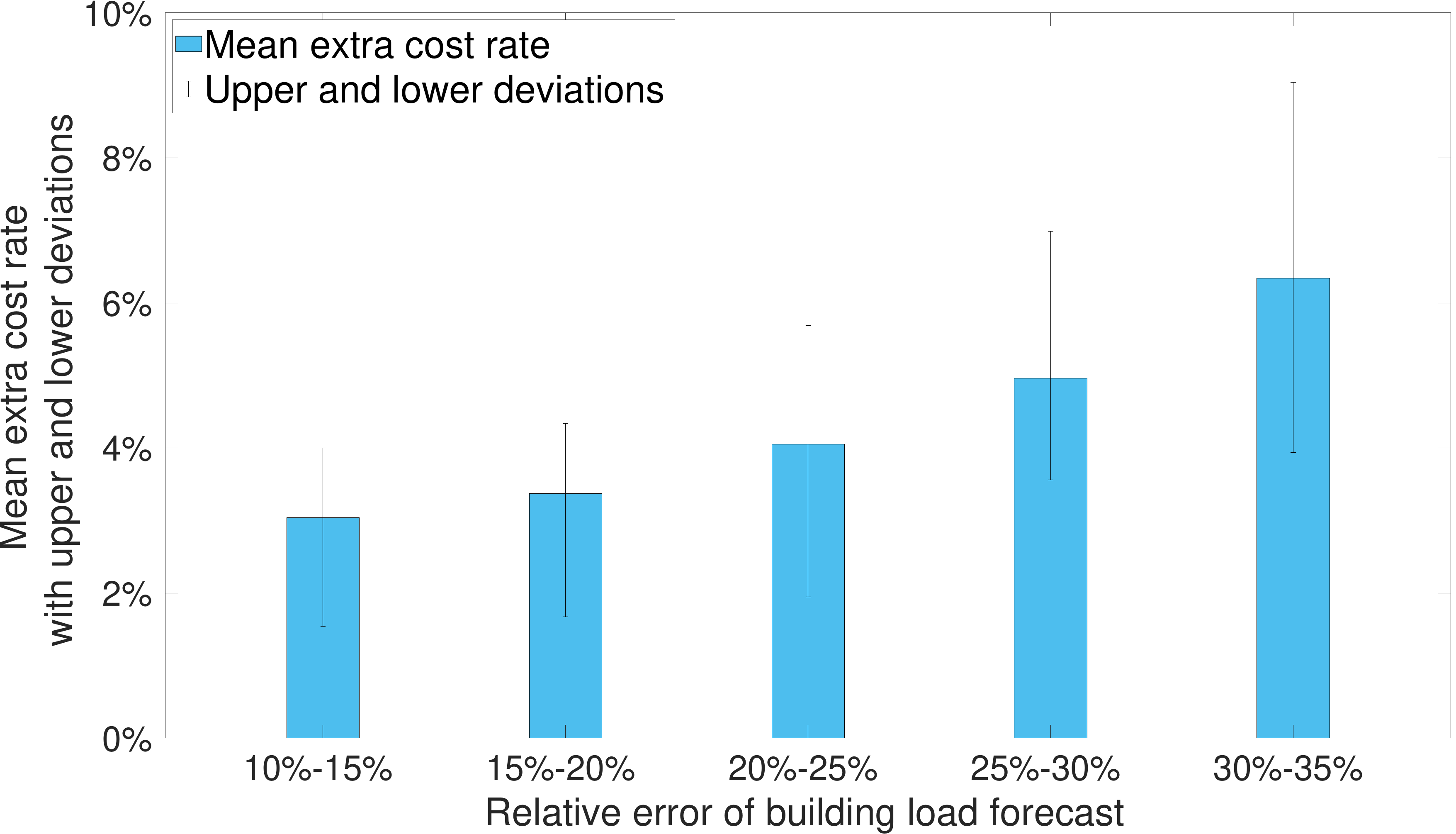}
    \caption{Correlation between mean extra cost rate and relative error in building load prediction.}
\label{fig: extra cost on load}
\end{figure}

To assess the impact of prediction errors, the extra costs are calculated using the predicted values of building demand, PV generation, and EV arrivals for each time slot, as defined in Eq.(\ref{eq26}). Fig.\ref{fig: extra cost on load} illustrates the relationship between the relative error in load prediction and the mean extra cost rate (i.e., the average extra cost compared to the minimum cost), including upper and lower standard deviations. As the error in load prediction increases, the mean extra cost rate also increases, reaching an average of 6.34\% within a range of [3.91, 9.1\%] when the relative error in load prediction is between 30 and 35\%.

\begin{figure}[!t]
    \centering
    \includegraphics[width=0.8\linewidth]{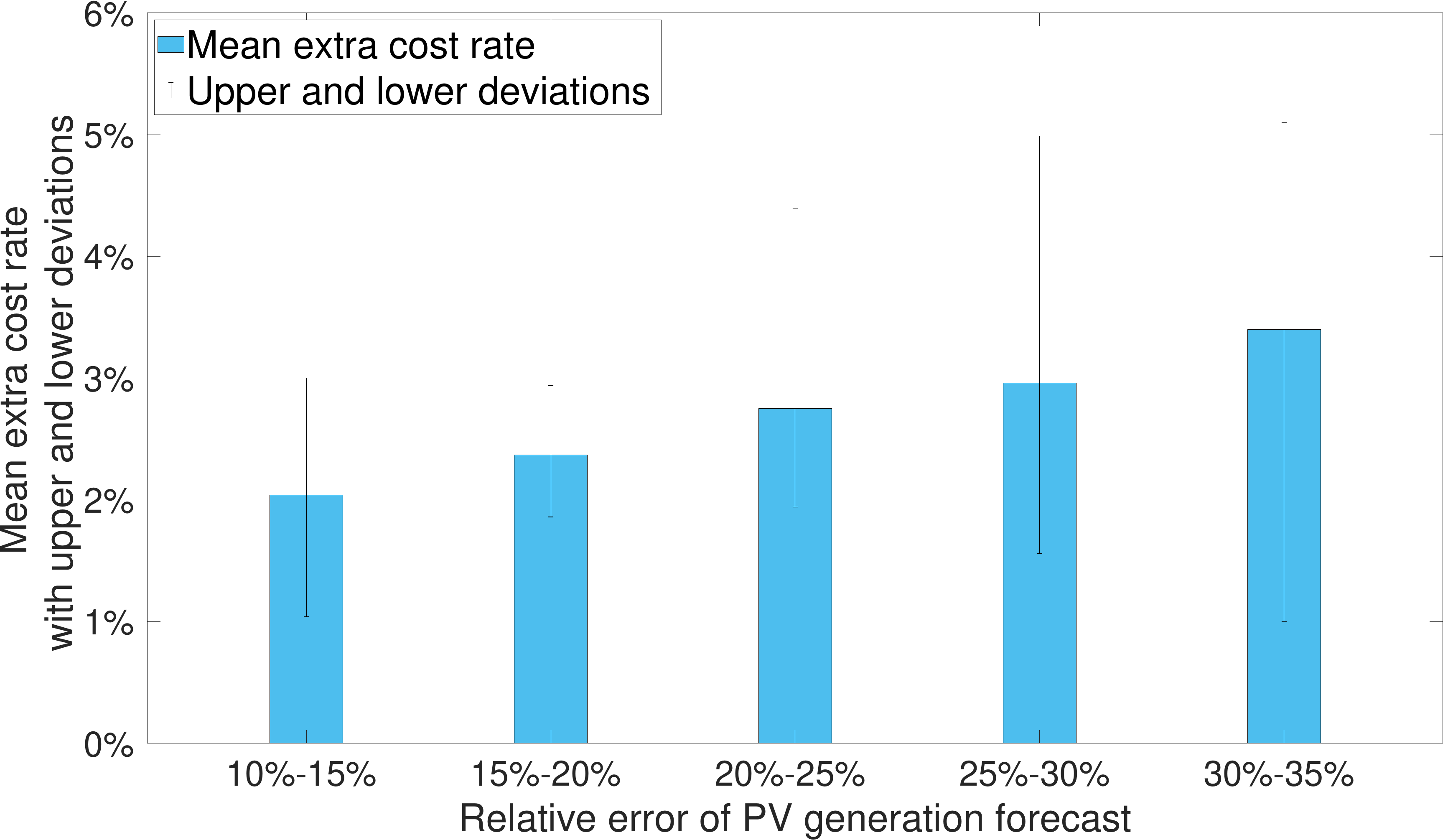}
    \caption{Correlation between the mean extra cost rate and the relative error in PV generation prediction.}
\label{fig: extra cost on PV}
\end{figure}

The relationship between the relative error in PV generation forecasts and the mean extra cost rate is depicted in Fig.~\ref{fig: extra cost on PV}. It shows that the mean extra cost rate marginally increases with the increase in the relative error of the predicted PV generation. The influence of errors in PV generation predictions is less significant compared to the impact of load prediction errors. These results are derived from our simulation settings and can be explained by the fact that, EVs are often not present at home during daytime hours when PV systems are active, thereby lessening the impact of PV generation prediction errors on the overall effectiveness of EV value-stacking.

\begin{figure}[!t]
    \centering
    \includegraphics[width=0.8\linewidth]{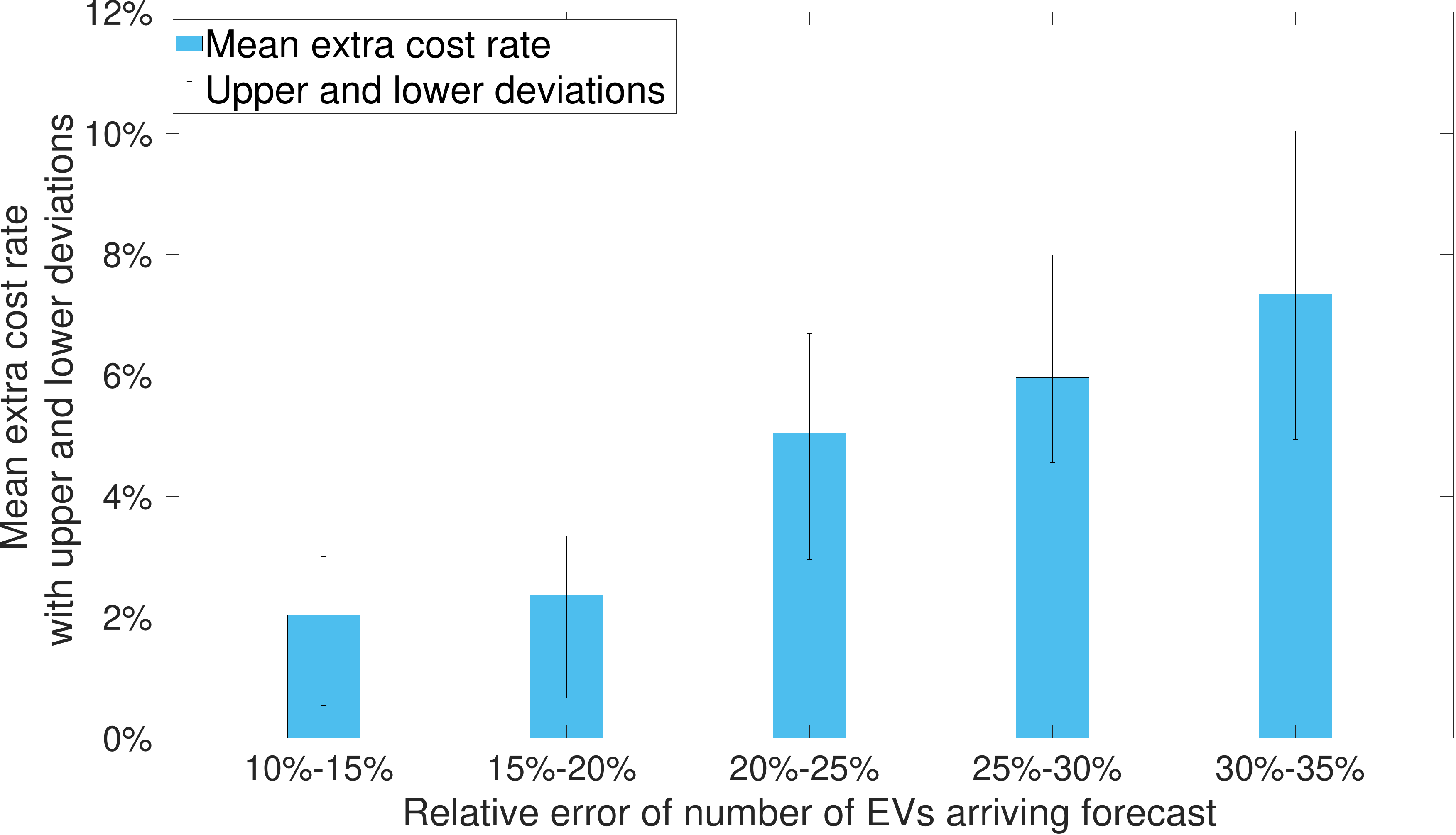}
    \caption{Correlation between the mean extra cost rate and the relative error in EV arrivals prediction.}
\label{fig: extra cost on EV}
\end{figure}

Fig. \ref{fig: extra cost on EV} depicts the correlation between the relative error in forecasted EV arrivals and the average extra cost rate. Specifically, when the relative prediction errors range between 10-15\% and 15-20\%, the mean extra cost rates decrease to 2.04\% and 2.37\%, respectively. However, as the relative prediction error increases to a range between 20 and 25\%, the mean extra cost rate sharply rises to 5.05\%. This significant increase is primarily due to the impact of the prediction error on EV scheduling, which can cause some vehicles to charge prematurely to meet their needs before departure or lead to increased charging demand during peak load hours. Additionally, when the relative load prediction error lies between 30-35\%, the mean extra cost rate climbs to 7.34\%, with a range of [4.6, 9.8\%]. Thus, the average extra cost is more significantly  affected by the prediction error in EV arrivals compared to errors in load demand and PV generation forecasts.

\subsection{Summary and Discussion}
Firstly, the value-stacking performance is evaluated through an analysis based on perfect predictions, under varying tariffs: TOU and TPT. Under TOU, significant cost reductions are achieved. TPT shows lower cost reductions but highlights the influence of peak load variability. Additionally, this work illustrates the decision-making processes in EV value stacking against three baselines under TOU and TPT tariffs. The results indicate that none of the single value stream can dominate the benefit, and value stacking can leverage the synergy to reduce overall system costs under both tariffs. Furthermore, the impact of V2B on HVAC usage is elucidated through a comparison between scenarios with and without V2B integration. The findings show that V2B significantly benefits the HVAC load management and optimizes energy use during peak periods. Moreover, the accuracy of the GRU-EN-TFD model in forecasting residential building demand, PV generation, and EV arrivals is evaluated. The findings show that the GRU-EN-TFD model consistently outperforms the LSTM in all forecast tasks for PV generation, residential load demand, and EV arrivals. This work analyzes the effects of prediction errors on the performance of value stacking. The simulation results show the mean additional costs due to uncertainties are notably influenced by errors in EV arrival predictions, underscoring the importance of accurate EV behavior predictions in the V2X value stacking.  

Secondly, this work considers each residential community purchasing energy based on the TOU or TPT and selling energy to provide V2G at the whole-market price. Specifically, in markets like NYISO and ISO-NE, whole-market prices are often much lower than off-peak rates in TOU and energy price in TPT. Therefore, communities are inclined to store energy in EV batteries during off-peak periods and then utilize it locally during peak times when tariffs are higher, which significantly limits the contribution of V2G. In contrast, V2B allows communities to leverage the stored energy in EV batteries to offset the higher tariffs during peak times, thus providing substantial cost savings. The profitability of performing V2G is not sufficient to compensate for the cost of purchasing or competitive to the benefit of offseting the higher tariffs via V2B. Therefore, V2B provides greater benefits to communities in these three energy markets compared to V2G. Looking ahead, the future of V2X presents promising opportunities. For V2G, there are expectations for enhanced services that contribute to grid stability and resilience enhancement with greater economic incentives \citep{inci2022integrating}. Regarding V2B, continuous innovation in smart building technologies and energy management systems is expected to broaden the applications and benefits of V2B, positioning it as a key component of urban energy solutions \citep{islam2022state}. 

Thirdly, despite the advantage and suitability of the dynamic RHO in the V2X value-stacking problem, such centralized optimization may face potential limitations when privacy concerns arise or when scaled to a significantly larger system. Specifically, the centralized nature of dynamic RHO requires access to detailed information about the distribution network and each residential community. As such, a strong commitment to collaboration is necessary to ensure that all required information is available for problem-solving. In the context of studying the benefits of V2X value stacking in our work, centralized problem-solving is effective. However, when concerns about privacy leakage and vulnerability to cyber-attacks emerge, centralized optimization may no longer be suitable. In these cases, privacy-preserving solutions are needed, often in decentralized manners, e.g., using decentralized optimization and blockchain technologies. For example, decentralized methods can decompose centralized optimization problems into smaller sub-problems that can be solved in parallel, thus helping preserve the privacy of each community. 

Another potential shortcoming of centralized RHO is scalability. While the dynamic nature of the RHO allows for real-time adaptation and updates, the method may face scalability challenges as the application scale extends with a much greater number of communities. The computational complexity associated with handling a large number of entities with integer variables may strain computational resources. This could potentially limit the method's feasibility for large-scale applications. A promising solution to address these scalability challenges involves decomposing centralized optimization problems into smaller sub-problems using algorithms such as the Alternating Direction Method of Multipliers (ADMM), which can be solved in parallel more efficiently, as shown in \citep{Zhou2020admm}, \citep{yang2021privacy}, and \citep{yang2022fully}. Our future work aims to consider decentralized optimization methods for V2X value stacking with privacy-preserving at scale.

\color{black}
\section{Conclusion and Future Work}\label{Conclusion and Future Work}
This paper introduced a V2X value-stacking optimization problem that incorporates local network constraints to optimize the economic benefits of coordinating EV charging and discharging. This problem is tackled by formulating a dynamic rolling-horizon optimization problem to solve the optimal EV charge/discharge decisions. To address energy uncertainties, including local building load, PV generation, and EV arrivals, this work developed a Transformer-based forecasting model named GRU-EN-TFD. The accuracy of these predictions and their impact on the efficiency of value stacking are assessed. Our simulations reveal that, under both TOU and TPT tariffs, V2X value-stacking outperforms three individual value streams in terms of cost reductions. Additional analysis on the marginal contributions of three value streams across three markets-the Australian's NEM, ISO-NE, and NYISO in the US-showed that V2B is particularly effective in reducing costs. Our study also revealed that the proposed GRU-EN-TFD model outperforms benchmark forecast models. The mean additional costs due to uncertainties are notably influenced by errors in EV arrival predictions, whereas the impact of errors in PV generation predictions is less significant, underscoring the importance of accurate EV behavior predictions in V2X value stacking. 

The dynamic RHO method, while effective for the V2X value-stacking problem, faces limitations related to privacy concerns and scalability. The centralized nature of the approach requires access to detailed information, which poses privacy risks and vulnerability to cyber-attacks. Additionally, as the number of communities increases, the computational complexity rises significantly due to the presence of binary variables, making scalability a challenge. Our future work aims to address these challenges by integrating decentralized methods with our dynamic RHO approach, ensuring privacy preservation for a large number of residential communities and enhancing scalability for larger applications.

\section*{Acknowledgement}
This work was supported in part by the Australian Research Council (ARC) Discovery Early Career Researcher Award (DECRA) under Grant DE230100046.

We dedicate this paper to the memory of our esteemed colleague, Professor Ariel Liebman, whose contributions will always be remembered.

\bibliography{reference.bib} 

\end{document}